\documentclass[a4paper,francais]{article}

\usepackage{amsmath,amscd,latexsym}
\usepackage[psamsfonts]{amssymb}
\usepackage[ansinew]{inputenc}
\usepackage[frenchb]{babel}
\usepackage[all]{xy}
\usepackage{lscape}

\newtheorem{theorem}{Th\'eor\`eme}[subsection]
\newtheorem{prop}[theorem]{Proposition}
\newtheorem{lem}[theorem]{Lemme}
\newtheorem{corol}[theorem]{Corollaire}

\newcommand{\compt}{\stepcounter{theorem}\thetheorem}

\newenvironment{defi}{\bigskip\noindent\textbf{D\'efinition \compt \,}}{\bigskip}

\newenvironment{rem}{\noindent\textbf{Remarque \compt\/ \;}}{}
\newcommand{\debrem}{\begin{rem}}
\newcommand{\finrem}{\end{rem}}

\newenvironment{demo}{\noindent\textit{D\'emonstration. }}{\hspace*{\fill}$\square$\bigskip}
\newcommand{\dem}{\begin{demo}}
\newcommand{\med}{\end{demo}}
\newcommand{\carre}{\hspace*{\fill}$\square$}
\newenvironment{demoref}[1]{\noindent\textit{D\'emonstration (#1). }}{\hspace*{\fill}$\square$\bigskip}
\newcommand{\demr}{\begin{demoref}}
\newcommand{\rmed}{\end{demoref}}

\newfont\pesp{rsfs10 at 8pt}
\newfont\esp{rsfs10 at 12pt}

\newcommand{\cc}{\mathcal{C}}

\newcommand{\vv}{\mathcal{V}}
\newcommand{\uu}{\mathcal{U}}

\newcommand{\mm}{\mathcal{M}}
\newcommand{\nn}{\mathcal{N}}

\newcommand{\pp}{\mathcal{P}}
\newcommand{\qq}{\mathcal{Q}}
\newcommand{\rr}{\mathcal{R}}

\newcommand{\gcc}{\mathfrak{C}}
\newcommand{\gdd}{\mathfrak{D}}
\newcommand{\gff}{\mathfrak{F}}
\newcommand{\gww}{\mathfrak{W}}
\newcommand{\gqq}{\mathfrak{Q}}
\newcommand{\gll}{\mathfrak{L}}

\newcommand{\coc}{{\cc}ocont}
\newcommand{\cres}{c{\rr}es}
\newcommand{\cresc}{c{\rr}es_{c}}

\newcommand{\modq}{\mathfrak{ModQ}}
\newcommand{\modqc}{\modq^{c}}
\newcommand{\modqd}{\modq^{\dag}}
\newcommand{\modqpt}{\modq_{\ast}}
\newcommand{\modqcpt}{\modqpt^{c}}
\newcommand{\modqdpt}{\modqpt^{\dag}}
\newcommand{\modqst}{\modq_{st}}
\newcommand{\modqcst}{\modqst^{c}}

\newcommand{\hmodq}{\mathfrak{THQ}}
\newcommand{\hmodqc}{\hmodq^{c}}
\newcommand{\hmodqd}{\hmodq^{\dag}}
\newcommand{\hmodqpt}{\hmodq_{\ast}}

\newcommand{\hmodqdpt}{\hmodqpt^{\dag}}
\newcommand{\hmodqst}{\hmodq_{st}}
\newcommand{\hmodqcst}{\hmodqst^{c}}

\newcommand{\pder}{\mathfrak{PDer}}

\newcommand{\der}{\mathfrak{Der}}
\newcommand{\derd}{\der^{d}}
\newcommand{\derad}{\der_{ad}}
\newcommand{\deradpp}{\der_{ad}^{pp}}

\newcommand{\cat}{\mathfrak{Cat}}
\newcommand{\Cat}{\mathfrak{CAT}}
\newcommand{\Catps}{\text{2-$\Cat_{ps}$}}
\newcommand{\Catad}{\Cat_{ad}}

\newcommand{\tron}{\pi_{0}}
\newcommand{\ptron}{\pi_{0}^{iso}}

\newcommand{\simp}{\mathcal{S}}
\newcommand{\ens}{\mathcal{E}ns}

\newcommand{\Sp}{S\!p}
\newcommand{\Sps}{\Sp^{\Sigma}}

\newcommand{\D}{\mathbb{D}}
\newcommand{\hot}{\mathbb{H}}

\newcommand{\ch}{{\cc}h}
{\newcommand{\cyl}{{\cc}yl}

\newcommand{\Ho}{\mathrm{Ho}}
\newcommand{\Lb}{\textbf{L}}
\newcommand{\Rb}{\textbf{R}}

\newcommand{\com}{\text{\large{m}}}
\newcommand{\un}{\text{\large{u}}}
\newcommand{\tpn}{\text{\large{$\lambda$}}}
\newcommand{\uni}{\text{\Large{$\eta$}}}
\newcommand{\cou}{\text{\Large{$\varepsilon$}}}
\newcommand{\loc}{\text{\large{$\gamma$}}}

\newcommand{\pmm}{\text{\tiny{$\mm$}}}
\newcommand{\pnn}{\text{\tiny{$\nn$}}}

\begin{document}

\title{\textbf{Théories homotopiques de Quillen combinatoires et dérivateurs de Grothendieck}}

\author{Olivier Renaudin}
\date{}
\maketitle

\footnotetext{MSC2000~: 18G55, 55U35. \; Mots cl\'es~: modèles de Quillen, dérivateurs de Grothendieck.}

\begin{abstract}
On construit une pseudo-localisation de la 2-catégorie des catégories modèles de Quillen combinatoires
relativement aux équivalences de Quillen, puis on vérifie que celle-ci se plonge dans une 2-catégorie
de dérivateurs de Grothendieck.
\end{abstract}

\bigskip

\bigskip

L'objectif de ce papier est de comparer la 2-catégorie des théories homotopiques de Quillen combinatoires avec
une 2-catégorie de dérivateurs de Grothendieck. La 2-catégorie $\hmodqc$ des théories homotopiques de Quillen
combinatoires est la pseudo-localisation de la 2-catégorie $\modqc$ des catégories modèles de Quillen combinatoires
relativement aux équivalences de Quillen. Dans un premier temps, on utilise des résultats de D. Dugger pour
produire une construction de la 2-catégorie $\hmodqc$.
Dans un second temps, on utilise des résultats de D.-C. Cisinski pour obtenir une équivalence locale
$\hmodqc \rightarrow \derad$ où $\derad$ désigne la 2-catégorie des dérivateurs à droite et à gauche avec les
adjonctions pour 1-morphismes.

\medskip

La première section est consacrée à quelques rappels concernant le cadre 2-catégorique du papier et à la notion
de pseudo-localisation d'une 2-catégorie $\gcc$ relativement à une classe $\gww$ de ces 1-morphismes~: il s'agit
essentiellement d'un pseudo-foncteur $\Gamma : \gcc \longrightarrow \gcc[\gww^{-1}]$ ``2-universel" parmi ceux
envoyant $\gww$ dans les équivalences.

\smallskip

Le but de la deuxième section est de construire une pseudo-localisation de la 2-catégorie $\modqc$ des modèles
de Quillen combinatoires relativement aux équivalences de Quillen.

On introduit d'abord des objets cylindre et chemin dans la 2-catégorie $\modq$ des modèles de Quillen, qui mènent
naturellement aux homotopies de Quillen, i.e. aux transformations naturelles de $\modq$ qui sont des équivalences
faibles sur les cofibrants. L'intérêt principal de ces modèles cylindres pour la suite est qu'ils assurent qu'un
pseudo-foncteur de source $\modq$ qui envoie les équivalences de Quillen dans les équivalences, envoie également
les homotopies de Quillen dans les isomorphismes.

On se place ensuite dans la 2-catégorie $\modqc$ des modèles de Quillen combinatoires. On rappelle la notion de
modèle présentable et le théorème de ``résolution" de D. Dugger \cite{D2}, selon lequel tout modèle combinatoire
est but d'une équivalence de Quillen de source un modèle présentable. Suivant toujours des observations de
D. Dugger \cite{D1}, on énonce un résultat d'invariance homotopique qui reflète le caractère ``cofibrant" des
modèles présentables.

On utilise alors ces propriétés pour construire une 2-catégorie $\hmodqc$ et un pseudo-foncteur
$\Gamma : \modqc \rightarrow \hmodqc$, puis on vérifie que celui-ci possède la propriété universelle caractérisant
une pseudo-localisation de $\modqc$ relativement aux équivalences de Quillen. La 2-catégorie $\hmodqc$ a pour objets
les modèles de Quillen combinatoires et pour catégories de morphismes les localisations relativement aux homotopies
de Quillen des catégories de morphismes entre modèles présentables de $\modqc$.
On conclue cette partie par quelques observations, concernant notamment les troncations de la
pseudo-localisation $\Gamma : \modqc \rightarrow \hmodqc$.

Enfin on termine la deuxième section en indiquant comment les considérations qui précèdent s'appliquent également,
d'une part aux modèles pointés et aux modèles stables, et d'autre part aux modèles simpliciaux et aux modèles
spectraux, en utilisant à nouveau des observations de D. Dugger.

\smallskip

La troisième section débute par le rappel de définitions concernant les dérivateurs. On rappelle ensuite sommairement
la construction, dû à D.-C. Cisinski \cite{C1}, d'un pseudo-foncteur de la 2-catégorie des modèles de Quillen dans
celle des dérivateurs, dont on déduit un pseudo-foncteur $\hmodqc \rightarrow \derad$.
En utilisant un théorème de représentation également de D.-C. Cisinski \cite{C2}, on vérifie que ce dernier
pseudo-foncteur est une équivalence locale. Pour clore cette section, on ébauche l'étude de l'image essentielle de ce
même pseudo-foncteur en introduisant les dérivateurs de petite présentation.

\smallskip

Enfin, un appendice recueille divers diagrammes dont la vérification de la commutativité est requise dans certaines
démonstrations de la partie principale du texte.

\tableofcontents

\newpage

\section{Préliminaires 2-catégoriques}

On note $\cat$ la 2-catégorie des petites catégories, et $\Cat$ la 2-catégorie de toutes les catégories.

\subsection{Rappels}

On rappelle dans cette section quelques définitions concernant les 2-catégories, pseudo-foncteurs,
transformations pseudo-naturelles\dots (cf. e.g. \cite{L1,L2}), qui constituent le cadre dans lequel s'inscrivent
les sections suivantes.

\medskip

Une 2-catégorie est une $\Cat$-catégorie, i.e. une catégorie enrichie sur la catégorie des catégories.
On note $e$ la 2-catégorie ponctuelle.

\begin{defi} \label{def-ps-fct}
Soient $\gcc$ et $\gdd$ deux 2-catégories. Un \emph{pseudo-foncteur} $\Phi : \gcc \rightarrow \gdd$ est la donnée~:
\begin{itemize}
    \item d'une fonction $\Phi : Ob(\gcc) \rightarrow Ob(\gdd)$,
    \item pour tout $c,c' \in Ob(\gcc)$, d'un foncteur $\Phi_{c,c'} : \gcc(c,c') \longrightarrow \gdd(\Phi(c),\Phi(c'))$,
    \item pour tout $c,c',c'' \in Ob(\gcc)$, des isomorphismes naturels
$$\xymatrix{
\gcc(c,c') \times \gcc(c',c'') \ar[r]^(.6){\circ} \ar[d]_{\Phi_{c,c'} \times \Phi_{c',c''}} & \gcc(c,c'')
\ar[d]^{\Phi_{c,c''}} & & e \ar[r]^(.4){Id_{c}} \ar@{=}[d] &  \gcc(c,c) \ar[d]^{\Phi_{c,c'}}   \\
\gdd(\Phi(c),\Phi(c')) \times \gdd(\Phi(c'),\Phi(c'')) \ar[r]_(.65){\circ} \ar@{=>}[ur]_(.65){\com^{\Phi}_{cc'c''}}
& \gdd(\Phi(c),\Phi(c''))   & &   e \ar[r]_(.3){Id_{\Phi(c)}}  \ar@{=>}[ur]^(.45){\un^{\Phi}_{c}} & \gdd(\Phi(c),\Phi(c))  }$$
    soit, pour tout $c \in Ob(\gcc)$, d'un 2-isomorphisme
    $\un^{\Phi}_{c} : Id_{\Phi(c)} \stackrel{\sim}{\longrightarrow} \Phi(Id_{c})$,
    et pour toute paire $f, g$ de 1-morphismes composables de $\gcc$, d'un 2-isomorphisme
    $\com^{\Phi}_{g,f} : \Phi(g) \circ \Phi(f) \stackrel{\sim}{\longrightarrow} \Phi(g \circ f)$
    naturel en $f$ et en $g$,
\end{itemize}
telle que les diagrammes suivants commutent~:

\noindent - Axiome d'associativité~:
$$\xymatrix{
\Phi(h) \circ \Phi(g) \circ \Phi(f) \ar[d]_{\com^{\Phi}_{h,g} \circ \Phi(f)} \ar[rrr]^{\Phi(h) \circ \com^{\Phi}_{g,f}}
& & & \Phi(h) \circ \Phi(g \circ f) \ar[d]^{\com^{\Phi}_{h,g \circ f}} \\
\Phi(h \circ g) \circ \Phi(f) \ar[rrr]_{\com^{\Phi}_{h \circ g,f}} & & & \Phi(h \circ g \circ f)  }$$
- Axiomes d'unité~:
$$\xymatrix{
Id_{\Phi(c)} \circ \Phi(f) \ar@{=}[d] \ar[rr]^{\un^{\Phi}_{c} \circ \Phi(f)} & & \Phi(Id_{c}) \circ \Phi(f) \ar[d]^{\com^{\Phi}_{Id_{c}, f}}
 &  & \Phi(f) \circ  Id_{\Phi(c)} \ar@{=}[d] \ar[rr]^{\Phi(f) \circ \un^{\Phi}_{c}} &  & \Phi(f) \circ \Phi(Id_{c}) \ar[d]^{\com^{\Phi}_{f, Id_{c}}} \\
\Phi(f) \ar@{=}[rr]  & & \Phi(Id_{c} \circ f) & &  \Phi(f) \ar@{=}[rr]  & & \Phi(f \circ Id_{c})  }$$
Un 2-foncteur est un pseudo-foncteur dont les isomorphismes structuraux $\un$ et $\com$
sont des morphismes identités.
\end{defi}

\begin{defi} \label{def-ps-nat}
Soient deux 2-catégories $\gcc$ et $\gdd$ et deux pseudo-foncteurs $\Phi, \Psi : \gcc \rightarrow \gdd$.
Une \emph{transformation pseudo-naturelle} $\alpha : \Phi \Rightarrow \Psi$ est la donnée~:
\begin{itemize}
    \item pour tout $c \in Ob(\gcc)$, d'un 1-morphisme $\alpha_{c} : \Phi(c) \rightarrow \Psi(c)$,
    \item pour tout $c,c' \in Ob(\gcc)$, d'un isomorphisme naturel
    $$\xymatrix{
\gcc(c,c')  \ar[d]_{\Phi_{c,c'}} \ar[r]^(.4){\Psi_{c,c'}} & \gdd(\Psi(c),\Psi(c')) \ar[d]^{(\alpha_{c})^{*}} \\
\gdd(\Phi(c),\Phi(c')) \ar[r]_{(\alpha_{c'})_{*}} \ar@{=>}[ur]^{\tpn^{\alpha}_{cc'}} & \gdd(\Phi(c),\Psi(c')) }$$
    soit, pour tout 1-morphisme $f : c \rightarrow c'$ de $\gcc$, d'un 2-isomorphisme
    $\tpn^{\alpha}_{f} : \alpha_{c'} \circ \Phi(f) \stackrel{\sim}{\longrightarrow} \Psi(f) \circ \alpha_{c}$
    naturel en $f$,
\end{itemize}
telle que les diagrammes suivants commutent~:

 Axiome de composition~: \hspace{7cm} Axiome d'unité~:
$$\xymatrix{
\alpha_{c''} \circ  \Phi(g) \circ \Phi(f)  \ar[d]_{Id \ast \com^{\Phi}_{g,f}} \ar[r]^{\tpn^{\alpha}_{g} \ast Id} &
\Psi(g) \circ \alpha_{c'} \circ \Phi(f)  \ar[r]^{Id \ast \tpn^{\alpha}_{f}} &
\Psi(g) \circ \Psi(f) \circ \alpha_{c} \ar[d]^{\com^{\Psi}_{g,f} \ast Id} &
\alpha_{c} \circ Id_{\Phi(c)}  \ar@{=}[r] \ar[d]_{Id \ast \un^{\Phi}_{c}} &
Id_{\Psi(c)} \circ \alpha_{c} \ar[d]^{\un^{\Psi}_{c} \ast Id}  \\
\alpha_{c''} \circ \Phi(g \circ f) \ar[rr]_{\tpn^{\alpha}_{g \circ f}} & & \Psi(g \circ f) \circ \alpha_{c}
&  \alpha_{c} \circ \Phi(Id_{c})  \ar[r]_{\tpn^{\alpha}_{Id}} &  \Psi(Id_{c}) \circ \alpha_{c}  }$$
Une transformation naturelle est une transformation pseudo-naturelle dont les isomorphismes structuraux $\tpn$
sont des morphismes identités.
\end{defi}

\begin{defi} \label{def-modif}
Soient deux 2-catégories $\gcc$ et $\gdd$, deux pseudo-foncteurs $\Phi, \Psi : \gcc \rightarrow \gdd$
et deux transformations pseudo-naturelles $\alpha, \beta : \Phi \Rightarrow \Psi$.
Une \emph{modification} $\theta : \alpha \Rrightarrow \beta$ est la donnée, pour tout $c \in Ob(\gcc)$,
d'un 2-morphisme $\theta_{c} : \alpha_{c} \rightarrow \beta_{c}$ tels que, pour tout 1-morphisme
$f : c \rightarrow c' \in \gcc$, le diagramme suivant commute~:
$$\xymatrix{
\alpha_{c'} \circ \Phi(f) \ar[d]_{\theta_{c'} \ast Id} \ar[r]^{\tpn^{\alpha}_{f}} &
\Psi(f) \circ \alpha_{c} \ar[d]^{Id \ast \theta_{c}} \\
\beta_{c'} \circ \Phi(f) \ar[r]_{\tpn^{\beta}_{f}} & \Psi(f) \circ \beta_{c} }$$
\end{defi}

Avec des compositions évidentes, on obtient une 3-catégorie $\Catps$ dont les objets sont les 2-catégories,
dont les 1-morphismes sont les pseudo-foncteurs, dont les 2-morphismes sont les transformations
pseudo-naturelles, et dont les 3-morphismes sont les modifications.

\begin{defi} \label{def-eq-loc}
Un pseudo-foncteur $\Phi : \gcc \rightarrow \gdd$ est une \emph{équilvalence locale} (resp. un \emph{plongement}) si,
pour tout $c,c' \in Ob(\gcc)$, le foncteur $\Phi_{c,c'} : \gcc(c,c') \longrightarrow \gdd(\Phi(c),\Phi(c'))$
est une équivalence de catégorie (resp. un isomorphisme).

Un pseudo-foncteur $\Phi : \gcc \rightarrow \gdd$ est \emph{2-essentiellement surjectif}
(resp. \emph{essentiellement surjectif}) si, pour tout $d \in Ob(\gdd)$, il existe $c \in Ob(\gcc)$
et une équivalence (resp. un isomorphisme) $\Phi(c) \rightarrow d$.

Un pseudo-foncteur $\Phi : \gcc \rightarrow \gdd$ est une \emph{biéquivalence} (resp. une \emph{équivalence})
s'il existe un pseudo-foncteur $\Psi : \gcc \rightarrow \gdd$ et des équivalences (resp. des isomorphismes)
$\Phi \circ \Psi \longrightarrow Id$ et $Id \longrightarrow \Psi \circ \Phi$.
\end{defi}

\begin{prop}\cite{L2}
Un pseudo-foncteur $\Phi : \gcc \rightarrow \gdd$ est une biéquilvalence si et seulement s'il est
une équilvalence locale et 2-essentiellement surjectif.   \carre
\end{prop}

\medskip

Dans la suite, on utilisera la remarque suivante.

\begin{prop}
Soient deux 2-catégories $\gcc$ et $\gdd$, deux pseudo-foncteurs $\Phi, \Psi : \gcc \rightarrow \gdd$
et une transformation pseudo-naturelle $\alpha : \Phi \Rightarrow \Psi$.
L'axiome d'unité de $\alpha$ est redondant~: il découle des axiomes d'unité de $\Phi$ et $\Psi$
et de l'axiome de composition de $\alpha$.
\end{prop}

\dem
Des axiomes d'unité de $\Phi$ et $\Psi$, il découle, pour tout 1-morphisme $f : c \rightarrow c' \in \gcc$,
la commutativité du diagramme suivant~:
$$\xymatrix{
\alpha_{c'} \circ \Phi(f)
\ar[dd]_{Id \ast \un^{\Phi}_{c'} \ast Id}  \ar[rr]^(.45){\un^{\Psi}_{c'} \ast Id}
\ar[drr]_(.6){\tpn^{\alpha}_{f}}  \ar@{=}[ddrr]  & &
\Psi(Id_{c'}) \circ \alpha_{c'} \circ \Phi(f)  \ar[rr]^{Id \ast \tpn^{\alpha}_{f}} & &
\Psi(Id_{c'}) \circ \Psi(f) \circ \alpha_{c}  \ar[dd]^{\com^{\Psi}_{Id,f} \ast Id}      \\
 & & \Psi(f) \circ \alpha_{c} \ar@{=}[drr] \ar[urr]_{\un^{\Psi}_{c'} \ast Id}   & &   \\
\alpha_{c'} \circ \Phi(Id_{c'}) \circ \Phi(f)  \ar[rr]_(.55){Id \ast \com^{\Psi}_{Id,f}}  & &
\alpha_{c'} \circ \Phi(f) \ar[rr]_{\tpn^{\alpha}_{f}}  & &  \Psi(f) \circ \alpha_{c}     }$$
Par comparaison avec l'axiome de composition de $\alpha$ pour les 1-morphismes $Id_{c'}$ et $f$,
on obtient le diagramme commutatif~:
$$\xymatrix{
\alpha_{c'} \circ \Phi(f) \ar@{=}[rr] \ar[d]_{Id \ast \un^{\Phi}_{c'} \ast Id}  & &
\alpha_{c'} \circ \Phi(f) \ar[d]^{\un^{\Psi}_{c'} \ast Id}  \\
\alpha_{c'} \circ \Phi(Id_{c'}) \circ \Phi(f) \ar[rr]_{\tpn^{\alpha}_{Id_{c'}} \ast Id}  & &
\Psi(Id_{c'}) \circ \alpha_{c'} \circ \Phi(f)  }$$
En choisissant $f = Id_{c'}$ et en utilisant le 2-isomorphisme
$\un^{\Phi}_{c'} : Id_{\Phi(c')} \stackrel{\sim}{\longrightarrow} \Phi(Id_{c'})$,
on obtient l'axiome d'unité de $\alpha$.
\med

\subsection{Pseudo-localisation}

Dans toute cette section, on considère une 2-catégorie $\gcc$ et une classe $\gww$ de ses 1-morphismes.
Pour deux 2-catégories $\gcc$ et $\gdd$, on note $[\gcc,\gdd]$ pour $\Catps(\gcc,\gdd)$, la 2-catégorie des
1-morphismes de $\gcc$ dans $\gdd$.

\begin{defi} \label{def-p-loc}
Une \emph{pseudo-localisation} de $\gcc$ relativement à $\gww$ est la donnée d'une 2-catégorie $\gcc[\gww^{-1}]$ dotée
d'un pseudo-foncteur $\Gamma : \gcc \longrightarrow \gcc[\gww^{-1}]$ envoyant $\gww$ dans les équivalences, tels que,
pour toute 2-catégorie $\gdd$, le 2-foncteur induit
$$\Gamma^{*} : [\gcc[\gww^{-1}],\gdd] \longrightarrow [\gcc,\gdd]_{\gww}$$
est une équivalence de 2-catégories, où $[\gcc,\gdd]_{\gww}$ désigne la sous-2-catégorie pleine de $[\gcc,\gdd]$
constituée des pseudo-foncteurs envoyant $\gww$ dans les équivalences.
\end{defi}

Il découle immédiatement de la définition qu'une pseudo-localisation est unique à équivalence près, elle-même unique
à isomorphisme unique près, et que toute 2-catégorie équivalente à une pseudo-localisation détermine une une
pseudo-localisation.

De plus, si $\Gamma : \gcc \longrightarrow \gcc[\gww^{-1}]$ est une pseudo-localisation de $\gcc$ relativement
à $\gww$, alors pour tout pseudo-foncteur $\Phi : \gcc \longrightarrow \gdd$ envoyant $\gww$ dans les équivalences,
la 2-catégorie $\gff_{\Gamma}(\Phi)$ dont~:
\begin{itemize}
    \item les objets sont les paires $(\Psi,\alpha)$ constituée d'un pseudo-foncteur
    $\Psi : \gcc[\gww^{-1}] \longrightarrow \gdd$ et d'un isomorphisme pseudo-naturel
    $\alpha : \Psi \circ \Gamma \stackrel{\sim}{\Longrightarrow} \Phi$,
    \item les morphismes $(\Psi,\alpha) \rightarrow (\Psi',\alpha')$ sont les transformations pseudo-naturelles
    $\chi : \Psi \Rightarrow \Psi'$ telles que $\alpha' \circ (\chi \circ \Gamma) = \alpha$,
    \item les 2-morphismes $\chi \rightarrow \chi'$ sont les modifications
    $\theta : \chi \Rightarrow \chi'$ telles que $Id_{\alpha'} \ast (\theta \circ \Gamma) = Id_{\alpha}$,
\end{itemize}
est équivalente à la 2-catégorie ponctuelle.

\begin{prop} \label{p-str-ploc}
Soit $\Gamma : \gcc \longrightarrow \gcc[\gww^{-1}]$ une pseudo-localisation de $\gcc$ relativement à $\gww$.
Il existe une pseudo-localisation équivalente $\Gamma' : \gcc \longrightarrow \gll$ telle que le pseudo-foncteur
$\Gamma'$ soit l'identité sur les objets.
Si le pseudo-foncteur $\Gamma$ est l'identité sur les objets, alors, pour tout pseudo-foncteur
$\Phi : \gcc \longrightarrow \gdd$ envoyant $\gww$ dans les équivalences, il existe
$(\widetilde{\Phi},\alpha) \in \gff_{\Gamma}(\Phi)$ tel que, pour tout $c \in \gcc$, $\widetilde{\Phi}(c) = \Phi(c)$
et $\alpha_{c} = Id_{\Phi(c)}$.
\end{prop}

\dem
Soit $\Gamma^{*}(\gcc[\gww^{-1}])$ la 2-catégorie ayant pour objets ceux de $\gcc$ et pour catégorie
de morphismes $\Gamma^{*}(\gcc[\gww^{-1}])(c,c') = \gcc[\gww^{-1}](\Gamma(c),\Gamma(c'))$.
Le pseudo-foncteur $\Gamma : \gcc \longrightarrow \gcc[\gww^{-1}]$ est le composé de pseudo-foncteurs
$\hat{\Gamma} : \gcc \longrightarrow \Gamma^{*}(\gcc[\gww^{-1}])$ et
$\check{\Gamma} : \Gamma^{*}(\gcc[\gww^{-1}]) \longrightarrow \gcc[\gww^{-1}]$
où $\hat{\Gamma}$ est l'identité sur les objets et $\check{\Gamma}$ l'identité sur les catégories de morphisme.
Comme $\hat{\Gamma}$ envoie  $\gww$ dans les équivalences , il existe $(\Psi,\alpha) \in \gff_{\Gamma}(\hat{\Gamma})$,
et donc tel que $(\check{\Gamma} \circ \Psi, \check{\Gamma} \circ \alpha) \in \gff_{\Gamma}(\Gamma)$. Il existe donc
un isomorphisme pseudo-naturel $\chi : \check{\Gamma} \circ \Psi \Rightarrow Id$ qui montre que $\check{\Gamma}$
est essentiellement surjectif et donc une équivalence, ce qui assure que $\hat{\Gamma}$ est une pseudo-localisation
de $\gcc$ relativement à $\gww$.

Soient maintenant une pseudo-localisation $\Gamma : \gcc \longrightarrow \gcc[\gww^{-1}]$ telle que $\Gamma$ soit
l'identité sur les objets, et un pseudo-foncteur $\Phi : \gcc \longrightarrow \gdd$ envoyant $\gww$ dans les
équivalences. Soit $(\Psi,\beta) \in \gff_{\Gamma}(\Phi)$. On définit un pseudo-foncteur
$\widetilde{\Phi} : \gcc[\gww^{-1}]\longrightarrow \gdd$ en posant, pour tout $c, c' \in \gcc$,
$\widetilde{\Phi}(c) = \Phi(c)$ et
$\widetilde{\Phi}_{c,c'} : \gcc[\gww^{-1}](c,c') \longrightarrow \gdd(\widetilde{\Phi}(c),\widetilde{\Phi}(c'))$
égal à la composée
$$\xymatrix{
\gcc[\gww^{-1}](c,c') \ar[r]^{\Psi_{c,c'}} &   \gdd(\Psi(c),\Psi(c')) \ar[r]^{(\beta_{c'})^{*}}  &
\gdd(\Psi(c),\widetilde{\Phi}(c')) \ar[r]^{(\beta_{c}^{-1})^{*}} &   \gdd(\widetilde{\Phi}(c),\widetilde{\Phi}(c'))  }$$
Les 2-isomorphismes structuraux sont définis, pour tout $c \in \gcc$ et pour les 1-morphisme composables
$f : c \rightarrow c', g : c' \rightarrow c'' \in \gcc$, par
$\un^{\widetilde{\Phi}}_{c} = \beta_{c} \circ \un^{\Psi}_{c} \circ \beta_{c}^{-1} :
Id_{\widetilde{\Phi}(c)} \stackrel{\sim}{\longrightarrow} \widetilde{\Phi}(Id_{c})$ et
$$\com^{\widetilde{\Phi}}_{g,f} = \beta_{c''} \circ \com^{\Psi}_{g,f} \circ \beta_{c}^{-1} :
\widetilde{\Phi}(g) \circ \widetilde{\Phi}(f) =
\beta_{c''} \circ \Psi(g) \circ \beta_{c'} \circ \beta_{c'}^{-1} \circ \Psi(f) \circ \beta_{c}^{-1}
\stackrel{\sim}{\longrightarrow} \widetilde{\Phi}(g \circ f)$$
Les axiomes de pseudo-fonctorialité de $\widetilde{\Phi}$ se déduisent immédiatement de ceux de $\Psi$.
On définit un ismorphisme pseudo-naturel $\alpha : \widetilde{\Phi} \circ \Gamma \Rightarrow \Phi$ en posant,
pour tout $c \in \gcc$, $\alpha_{c} = Id_{c} : \widetilde{\Phi}(\Gamma(c)) \Rightarrow \Phi(c)$,
et pour tout 1-morphisme $f : c \rightarrow c' \in \gcc$
$$\tpn^{\alpha}_{f} = \tpn^{\beta}_{f} \circ \beta_{c}^{-1} :
\alpha_{c'} \circ \widetilde{\Phi}(\Gamma(f)) = \beta_{c'} \circ \Psi(\Gamma(f)) \circ \beta_{c}^{-1}
\Longrightarrow \Phi(f) \circ \alpha_{c}$$
Les axiomes de pseudo-naturalité de $\alpha$ se déduisent immédiatement de ceux de $\beta$.
\med

\bigskip

A toute 2-catégorie $\gcc$ sont associées trois 1-catégories possédant les mêmes objets.
La 1-catégorie sous-jacente $\gcc_{\circ}$ est obtenue en oubliant les 2-morphismes.
La 1-troncation $\tron(\gcc)$ est le quotient de $\gcc_{\circ}$ par la relation d'équivalence engendrée
par les 2-morphismes.
La pseudo-troncation $\ptron(\gcc)$ est le quotient de $\gcc_{\circ}$ par la relation d'équivalence engendrée
par les 2-isomorphismes.
Ces deux dernières catégories viennent donc avec deux foncteurs quotients $\gcc_{\circ} \longrightarrow \tron(\gcc)$
et $\gcc_{\circ} \longrightarrow \ptron(\gcc)$,
On a défini ainsi deux 2-foncteurs $\tron, \ptron : \Catps \longrightarrow \Cat$.

Les classes de morphisme de $\tron(\gcc)$ et $\ptron(\gcc)$ images de $\gww$ sont encore notées $\gww$.
On a alors un triangle commutatif dans $\Cat$~:
$$\xymatrix{  & \tron(\gcc) \ar[dl]_{} \ar[dr]^{\tron(\Gamma)} &  \\
\tron(\gcc)[\gww^{-1}]  \ar[rr]_{\overline{\tron(\Gamma)}} & & \tron(\gcc[\gww^{-1}])  }$$

\begin{prop} \label{p-tr-ploc}
Soit $\Gamma : \gcc \longrightarrow \gcc[\gww^{-1}]$ une pseudo-localisation de $\gcc$ relativement à $\gww$.
Le foncteur $\overline{\tron(\Gamma)} : \tron(\gcc)[\gww^{-1}] \longrightarrow \tron(\gcc[\gww^{-1}])$
est une équivalence. Il s'agit d'un isomorphisme si le pseudo-foncteur $\Gamma$ est l'identité sur les objets.
\end{prop}

\dem
Soit $D$ un catégorie. On note $D_{\delta}$ la 2-catégorie ``discrète'' associée (i.e. sans 2-morphisme
non-triviaux). On a une 2-adjonction $\tron : \Catps \rightleftarrows \Cat : (-)_{\delta}$, qui montre que
le foncteur $\tron(\Gamma)^{*} : [\tron(\gcc[\gww^{-1}]),D] \longrightarrow [\tron(\gcc),D]_{\gww}$
est bien une équivalence, ce qui démontre la première assertion de l'énoncé, et que $\tron(\Gamma)^{*}$
est un isomorphisme si et seulement si le 2-foncteur
$\Gamma^{*} : [\gcc[\gww^{-1}],D_{\delta}] \longrightarrow [\gcc,D_{\delta}]_{\gww}$ est un isomorphisme.

On suppose désormais que le pseudo-foncteur $\Gamma : \gcc \longrightarrow \gcc[\gww^{-1}]$ est
l'identité sur les objets. Comme le 2-foncteur $\Gamma^{*}$ est pleinement fidèle, on est ramené à vérifier
que $\Gamma^{*}$ est une bijection sur les objets. On note qu'une transformation pseudo-naturelle à valeur
dans $D_{\delta}$ est nécessairement strict. Soit un pseudo-foncteur $\Phi : \gcc \longrightarrow D_{\delta}$
envoyant $\gww$ dans les équivalences. Par la proposition \ref{p-str-ploc}, il existe un 2-foncteur
$\widetilde{\Phi} : \gcc[\gww^{-1}] \longrightarrow D_{\delta}$ et un isomorphisme naturel
$\alpha : \widetilde{\Phi} \circ \Gamma \stackrel{\sim}{\Longrightarrow} \Phi$, tels que, pour tout $c \in \gcc$,
$\widetilde{\Phi}(c) = \Phi(c)$ et $\alpha_{c} = Id_{\Phi(c)}$. Autrement dit, on a une égalité
$\widetilde{\Phi} \circ \Gamma = \Phi$, ce qui montre la surjectivité de $\Gamma^{*}$.
Soient deux pseudo-foncteurs $\Psi, \Psi' : \gcc[\gww^{-1}] \longrightarrow D_{\delta}$ tels que
$\Gamma^{*}(\Psi) = \Gamma^{*}(\Psi')$. Il existe une unique transformation naturelle
$\chi : \Psi \rightarrow \Psi'$ tel que $\Gamma^{*}(\chi) = Id_{\Gamma^{*}(\Psi)}$. Pour tout $c \in \gcc[\gww^{-1}]$,
on a $\chi_{c} = \chi_{\Gamma(c)} = \Gamma^{*}(\chi)_{c} =Id$, soit $\Psi = \Psi'$, ce qui montre l'injectivité
de $\Gamma^{*}$.
\med

\medskip

Soit $I = (0 \rightleftarrows 1)$ le groupoïde possédant deux objets et exactement deux morphismes non triviaux.
Un foncteur $I \rightarrow C \in \Cat$ correspond à la donnée d'un isomorphisme de $C$.

\begin{prop} \label{p-ptr-loceq}
Soit $\gcc$ une 2-catégorie. On suppose que tout objet de $\gcc$ possède une cotensorisation par $I$.
Alors la pseudo-troncation $\gcc_{\circ} \longrightarrow \ptron(\gcc)$ est la localisation relativement
aux équivalences de $\gcc_{\circ}$.
\end{prop}

\dem
On note $c^{I}$ la cotensorisation de $c$ par $I$. Chaque objet de $I$ détermine un 1-morphisme
$e_{i} : c^{I} \rightarrow c$, $i=0,1$, et l'équivalence de $I$ vers la catégorie ponctuelle détermine une
équivalence $\iota : c \rightarrow c^{I}$ tels que $e_{i} \circ \iota = Id$, pour $i=0,1$.

Soient $f, g \in \gcc_{\circ}(b,c)$ et $\tau : f \rightarrow g$ un 2-isomorphisme de $\gcc$.
Par définition de la cotensorisation, i.e. l'isomorphisme naturel $\Cat(I,\gcc(b,c)) \simeq \gcc(b,c^{I})$,
le 2-isomorphisme $\tau$ correspond à un 1-morphisme $h_{\tau} : b \rightarrow c^{I}$ tel que
$e_{0} \circ h_{\tau} = f$ et $e_{1} \circ h_{\tau} = g$.

Soient $D$ une catégorie et $F : \gcc_{\circ} \rightarrow D$ un foncteur envoyant les équivalences dans les
isomorphismes. Comme $F(\iota)$ est alors un isomorphisme d'inverse $F(e_{0})$ et $F(e_{1})$ , on a $F(e_{0}) = F(e_{1})$
et donc, pour tout 2-isomorphisme $\tau : f \rightarrow g$, $F(f) = F(g)$.
De la définition de $\gcc_{\circ} \longrightarrow \ptron(\gcc)$, il découle que l'on a vérifié la propriété
universelle de la localisation de $\gcc_{\circ}$ relativement aux équivalences.
\med

\section{Théories homotopiques de Quillen combinatoires}

On note $\modq$ la 2-catégorie des modèles de Quillen \cite{Ho1}.
Un 1-morphisme $F : \mm \rightarrow \nn$ est une adjonction de Quillen $F : \mm \rightleftarrows \nn : F^{\flat}$.
Un 2-morphisme $\tau : F \rightarrow G$ de $\modq$ est une transformation naturelle des adjoints à gauche
$F \rightarrow G$.
On note $\gqq$ la classe des 1-morphismes de $\modq$ constituée des équivalences de Quillen.

\subsection{Homotopies de Quillen et modèles chemins}

\subsubsection{Modèles chemins, modèles cylindres}

On note $[n]$ l'ensemble ordonné $\{0<...<n\}$ considéré comme une catégorie.
Soit $\mm$ une catégorie modèle. La factorisation de la codiagonale
$[0] \amalg [0] \stackrel{i}{\longrightarrow} [1] \stackrel{p}{\longrightarrow} [0]$
induit un diagramme~:
$$\xymatrix{
\mm  \ar@<.5ex>[r]^-{p^{*}} & \mm^{[1]} \ar@<.5ex>[l]^-{p_{*}} \ar@<.5ex>[r]^-{i^{*}} &
\mm \times \mm  \ar@<.5ex>[l]^-{i_{*}}        }$$
où $p^{*}$ est le foncteur ``diagramme constant'', $i^{*} = (e_{0},e_{1})$, $p_{*} = e_{0}$,
en notant $e_{i}$ l'évaluation en $i$, et $i_{*}(M_{0},M_{1}) = (M_{0} \times M_{1} \rightarrow M_{1})$.

En considérant la catégorie $[1]$ comme une catégorie indirecte, on peut munir $\mm^{[1]}$ d'une structure
de catégorie modèle de Reedy \cite[Chapter 15]{Hi}, qui n'est autre que la structure injective~: un morphisme
$f : X \rightarrow Y$ de $\mm^{[1]}$ est~:
\begin{itemize}
    \item une cofibration si $f_{i} : X_{i} \rightarrow Y_{i}$ est une cofibration pour $i = 0,1$,
    \item une équivalence faible si $f_{i} : X_{i} \rightarrow Y_{i}$ est une équivalence faible pour $i = 0,1$,
    \item une fibration si $f_{1} : X_{1} \rightarrow Y_{1}$ et $X_{0} \rightarrow X_{1} \times_{Y_{1}} Y_{0}$
    sont des fibrations.
\end{itemize}
Les adjonctions $(p^{*},p_{*})$ et $(i^{*},i_{*})$ sont ainsi des adjonctions de Quillen.

\bigskip

Dans le but de construire un ``objet chemin'' relativement aux équivalences de Quillen dans $\modq$, on considère
la localisation de Bousfield à droite \cite[3.3.1(2)]{Hi} relativement aux $e_{0}$-équivalences de la structure
injective de $\mm^{[1]}$.

\begin{prop}
La catégorie $\mm^{[1]}$ possède une structure de catégorie modèle telle qu'un morphisme $f : X \rightarrow Y$
est~:
\begin{itemize}
    \item une fibration si $f_{1} : X_{1} \rightarrow Y_{1}$ et $X_{0} \rightarrow X_{1} \times_{Y_{1}} Y_{0}$
    sont des fibrations,
    \item une équivalence faible si $f_{0} : X_{0} \rightarrow Y_{0}$ est une équivalence faible,
    \item une cofibration si $f_{i} : X_{i} \rightarrow Y_{i}$ est une cofibration pour $i = 0,1$, et
    $X_{1} \amalg_{X_{0}} Y_{0} \rightarrow  Y_{1}$ est une équivalence faible.
\end{itemize}
On note $\ch(\mm)$ cette catégorie modèle.
\end{prop}

\dem
Les fibrations et équivalences faibles de l'énoncé sont, par définition, celles de la structure localisée
de Bousfield à droite de la structure injective relativement aux $e_{0}$-équivalences. L'existence de la
structure de modèle localisé est équivalente à celle de la factorisation ``cofibration locale / fibration
triviale locale''. On va d'abord montrer que les cofibrations de l'énoncé sont des cofibrations locales,
puis produire la factorisation requise.

\smallskip

Soient deux morphismes $A \rightarrow B$ et $X \rightarrow Y$ de $\mm^{[1]}$, respectivement une cofibration
et une fibration triviale de l'énoncé. On fixe un morphisme du premier vers le second
$$\xymatrix{  A \ar[d] \ar[r] & X \ar[d] \\   B \ar[r] & Y }$$
et on note $f : A \rightarrow Y$ la composée. Il s'agit de construire un relèvement $B \rightarrow X$.

On se place dans la catégorie $\mm^{[1]}(f)$ des factorisations du morphisme $f : A \rightarrow Y$,
i.e. la catégorie $f/(\mm^{[1]}/Y)$ des objets au-dessous de $f$ dans celle des objets de $\mm^{[1]}$ au-dessus
de $Y$. Cette catégorie est munie de la structure modèle induite par la structure injective de $\mm^{[1]}$.
Comme il n'y aura pas d'ambiguité sur les morphismes structuraux, on dénotera les objets de $\mm^{[1]}(f)$ par
leur objet central.

On note $B'$ l'objet $(B_{0} \rightarrow B_{0} \amalg_{A_{0}} A_{1})$ de $\mm^{[1]}$ donné par l'injection
canonique.
Le fait que le morphisme $A \rightarrow B$ soit une cofibration de $\mm^{[1]}$ signifie, d'une part, que
les morphismes $A_{0} \rightarrow B_{0}$, $A_{1} \rightarrow B_{0} \amalg_{A_{0}} A_{1}$ et $A_{1} \rightarrow B_{1}$
sont des cofibrations, et, d'autre part, que le morphisme $B_{0} \amalg_{A_{0}} A_{1} \rightarrow B_{1}$ est
une équivalence faible. Autrement dit, que le morphisme $B' \rightarrow B$ de $\mm^{[1]}(f)$~:
$$(B_{0} \rightarrow B_{0} \amalg_{A_{0}} A_{1}) \longrightarrow (B_{0} \rightarrow B_{1})$$
est une équivalence faible entre cofibrants.
Le fait que le morphisme $X \rightarrow Y$ soit une fibration de $\mm^{[1]}$ signifie que $X$ est un objet fibrant
de $\mm^{[1]}(f)$.
Comme le morphisme $X \rightarrow Y$ est une fibration triviale de $\mm^{[1]}$, $X_{0} \rightarrow Y_{0}$ est une
fibration triviale, et comme $A_{0} \rightarrow B_{0}$ est une cofibration, on a un relèvement~:
$$\xymatrix{  A_{0} \ar[d] \ar[r] & X_{0} \ar[d] \\    B_{0} \ar[r] \ar[ru] & Y_{0} }$$
Le morphisme composé $B_{0} \rightarrow X_{0} \rightarrow X_{1}$ et le morphisme $A_{1} \rightarrow X_{1}$
déterminent un morphisme $B_{0} \amalg_{A_{0}} A_{1} \rightarrow  X_{1}$ dont on tire un morphisme
$B' \rightarrow X$ de $\mm^{[1]}$. Comme $B' \rightarrow B$ est une équivalence faible entre cofibrants et $X$
un objet fibrant, le morphisme $B' \rightarrow X$ s'étend, à homotopie près, en un morphisme $B \rightarrow X$
\cite[7.8.6]{Hi}, et on a ainsi obtenu le relèvement requis.

\smallskip

On va maintenant montrer que tout morphisme $f : X \rightarrow Y \in \mm^{[1]}$ se factorise en une
cofibration de l'énoncé suivie d'une fibration de l'énoncé (i.e. injective).
On factorise $f_{0}$ en une cofibration suivie d'une fibration triviale $X_{0} \rightarrow E_{0}' \rightarrow Y_{0}$.
Les morphismes $E_{0}' \rightarrow Y_{0} \rightarrow Y_{1}$ et $X_{1} \rightarrow Y_{1}$ induisent par propriété
universelle un morphisme $X_{1} \amalg_{X_{0}} E_{0}' \rightarrow Y_{1}$. On factorise ce dernier morphisme
en une cofibration triviale suivie d'une fibration $X_{1} \amalg_{X_{0}} E_{0}' \rightarrow E_{1} \rightarrow Y_{1}$.
Le morphisme composé $X_{1} \rightarrow X_{1} \amalg_{X_{0}} E_{0}' \rightarrow E_{1}$ est ainsi une cofibration.
En introduisant le produit fibré $E_{1} \times_{Y_{1}} Y_{0}$, on obtient donc le diagramme commutatif~:
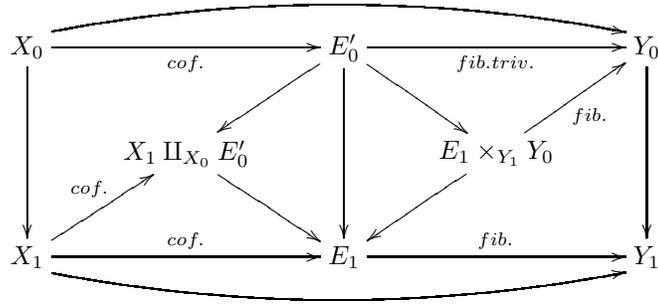
$$\xymatrix{
X_{0} \ar[dd] \ar[rr]_{cof.} \ar@/^1pc/@<1ex>[rrrr] &    & E_{0}' \ar[dd] \ar[rr]_{fib. triv.} \ar[dr] \ar[dl] &  & Y_{0} \ar[dd]  \\
 & X_{1} \amalg_{X_{0}} E_{0}' \ar[dr] &    & E_{1} \times_{Y_{1}} Y_{0} \ar[ld] \ar[ru]_{fib.} &                 \\
X_{1} \ar[rr]^{cof.} \ar[ru]^{cof.} \ar@/_1pc/@<-1ex>[rrrr] &    & E_{1} \ar[rr]^{fib.} &  & Y_{1}        }$$
Le morphisme $E_{0}' \rightarrow E_{1} \times_{Y_{1}} Y_{0}$ se factorise en une cofibration triviale suivie
d'une fibration $E_{0}' \rightarrow E_{0} \rightarrow E_{1} \times_{Y_{1}} Y_{0}$. Le morphisme $E_{0} \rightarrow Y_{0}$
ainsi obtenu est une fibration triviale, par l'axiome ``deux sur trois''.
Comme $E_{0}' \rightarrow E_{0}$ est une cofibration triviale, il en va de m\^{e}me de
$X_{1} \amalg_{X_{0}} E_{0}' \rightarrow X_{1} \amalg_{X_{0}} E_{0}$, et comme
$X_{1} \amalg_{X_{0}} E_{0}' \rightarrow E_{1}$ est également une cofibration triviale, le morphisme
$X_{1} \amalg_{X_{0}} E_{0} \rightarrow E_{1}$, obtenu par propriété universelle, est une équivalence faible.
On a donc la factorisation voulue~:
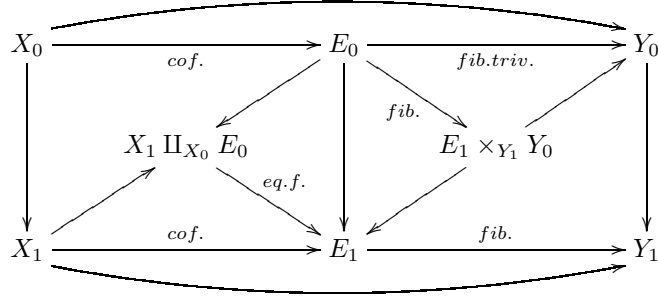
$$\xymatrix{
X_{0} \ar[dd] \ar[rr]_{cof.} \ar@/^1pc/@<1ex>[rrrr] &    & E_{0} \ar[dd] \ar[rr]_{fib. triv.} \ar[dr]_{fib.} \ar[dl] &  & Y_{0} \ar[dd]  \\
 & X_{1} \amalg_{X_{0}} E_{0} \ar[dr]^{eq.f.} &    & E_{1} \times_{Y_{1}} Y_{0} \ar[ld] \ar[ru] &                 \\
X_{1} \ar[rr]^{cof.} \ar[ru] \ar@/_1pc/@<-1ex>[rrrr] &    & E_{1} \ar[rr]^{fib.} &  & Y_{1}         }$$

\smallskip

Enfin, avec l'argument de retract habituel, la factorisation ci-dessus appliquée à une cofibration locale montre
que les cofibrations de l'énoncé coïncident avec les cofibrations locales.
\med

\smallskip

Il est clair que les foncteurs $p^{*} : \mm \rightarrow \ch(\mm)$ et $i^{*} : \ch(\mm) \rightarrow \mm \times \mm$
préservent cofibrations et cofibrations triviales, et que l'on a donc deux morphismes de $\modq$~:
$$\xymatrix{
\mm  \ar@<.5ex>[r]^-{p^{*}} & \ch(\mm) \ar@<.5ex>[l]^-{p_{*}} \ar@<.5ex>[r]^-{i^{*}} &
\mm \times \mm  \ar@<.5ex>[l]^-{i_{*}}        }$$
Qui plus est, comme l'unité de la pair $(p^{*},p_{*})$ est un isomorphisme et $e_{0}$ préserve et reflète
les équivalences faibles, cette pair est une équivalence de Quillen. La catégorie modèle $\ch(\mm)$
est appelée le \emph{modèle chemin} de $\mm$.

\smallskip

Dualement, on a un \emph{modèle cylindre} $\cyl(\mm)$ de $\mm$ défini par $\cyl(\mm) = \ch(\mm^{op})^{op}$,
et qui s'inscrit dans un diagramme~:
$$\xymatrix{
\mm \times \mm  \ar@<.5ex>[r]^-{i_{!}} & \cyl(\mm) \ar@<.5ex>[l]^-{i^{*}} \ar@<.5ex>[r]^-{p_{!}} &
\mm \ar@<.5ex>[l]^-{p^{*}}        }$$
où $p_{!} = e_{1}$. La catégorie modèle $\cyl(\mm)$ est la localisation de Bousfield à gauche relativement
aux $e_{1}$-équivalences de la structure projective de $\mm^{[1]}$.

\newpage

\subsubsection{Homotopies de Quillen}

\begin{defi} \label{def-htp-quil}
Soient $F_{0}, F_{1} : \mm \longrightarrow \nn$ deux morphismes de $\modq$.
\begin{enumerate}
    \item Une \emph{homotopie de Quillen} de $F_{0}$ vers $F_{1}$ est une transformation naturelle
          $\tau : F_{0} \rightarrow F_{1}$ telle que $\tau_{M} : F_{0}(M) \rightarrow F_{1}(M)$ soit
          une équivalence faible lorsque $M$ est cofibrant.
    \item Les morphismes $F_{0}$ et $F_{1}$ sont \emph{Quillen-homotopes} s'ils sont connectés par un zigzag
          d'homotopies de Quillen.
\end{enumerate}
On note $\qq_{\mm,\nn}$, ou plus simplement $\qq$, la classe des homotopies de Quillen de $\modq(\mm,\nn)$.
\end{defi}

\debrem
Suivant \cite[1.3.18]{Ho1}, la transformation naturelle de foncteurs dérivés totaux
$\Lb(\tau) : \Lb(F_{0}) \rightarrow \Lb(F_{1})$ induite par un 2-morphisme
$\tau : F_{0} \rightarrow F_{1}$ de $\modq$ est un isomorphisme si et seulement si
$\tau_{M} : F_{0}(M) \rightarrow F_{1}(M)$ est une équivalence faible pour tout $M$ cofibrant.
Il en découle qu'un morphisme de $\modq$ Quillen-homotope à une équivalence de Quillen est
également une équivalence de Quillen.
\finrem

\medskip

\begin{prop} \label{p-hq-hc}
Soient $F, G : \mm \longrightarrow \nn$ deux morphismes de $\modq$. La donnée d'une
homotopie de Quillen $\tau : F \rightarrow G$ est équivalente à la donnée d'un morphisme
$H : \mm \longrightarrow \ch(\nn)$ de $\modq$ faisant commuter le
diagramme
$$\xymatrix{
\mm \ar[r]^{H} \ar[dr]_{(F,G)} & \ch(\nn) \ar[d]^{i^{*}} \\
 & \nn \times \nn }$$
avec $H(M) = (\tau_{M} : F(M) \rightarrow G(M))$ pour tout $M \in \mm$.
\end{prop}

\dem
On note d'abord que la donnée d'un foncteur $H : \mm \longrightarrow \nn^{[1]}$ tel que $i^{*} \circ H = (F,G)$
est équivalente à celle d'une transformation naturelle $\tau^{H} : F \rightarrow G$ avec
$H(M) = (\tau^{H}_{M} : F(M) \rightarrow G(M))$ pour tout $M \in \mm$.

De plus, un tel foncteur $H$ possède un adjoint à droite si et seulement si $F$ et $G$ possèdent
un adjoint à droite. En effet, d'une part, $i^{*}$ possède un adjoint à droite, et, d'autre part, si on note
$F^{\flat}$ et $G^{\flat}$ des adjoints à droite de $F$ et $G$ respectivement, le foncteur
$H^{\flat} : \nn^{[1]} \longrightarrow \nn$, qui à $X = (X_{<} : X_{0} \rightarrow X_{1}) \in \nn^{[1]}$
associe le produit fibré de $F^{\flat}(X_{<}) : F^{\flat}(X_{0}) \longrightarrow F^{\flat}(X_{1})$ et de
$(\tau^{H})^{\flat}_{X_{1}} : G^{\flat}(X_{1}) \rightarrow F^{\flat}(X_{1})$ (où $(\tau^{H})^{\flat}$ désigne
la transformation conjuguée de $\tau^{H}$), est un adjoint à droite de $H$.

Un objet $X$ de $\ch(\nn)$ est cofibrant si et seulement si $X_{<} : X_{0} \rightarrow X_{1}$ est une
équivalence faible entre cofibrants. Donc un foncteur $H : \mm \longrightarrow \ch(\nn)$ tel que
$i^{*} \circ H = (F,G)$ préserve les cofibrants si et seulement si $\tau^{H}_{M} : F(M) \rightarrow G(M)$
est une équivalence faible pour tout $M$ cofibrant, autrement dit si $\tau^{H}$ est une homotopie de Quillen.
Il reste donc à vérifier qu'un foncteur $H : \mm \longrightarrow \ch(\nn)$ tel que $i^{*} \circ H = (F,G)$ est
un morphisme de $\modq$ si et seulement s'il préserve les cofibrants.

Lorsque $\nn^{[1]}$ est munie de la structure injective, $H : \mm \longrightarrow \nn^{[1]}$ est un morphisme
de $\modq$ si et seulement si $F$ et $G$ le sont également. Comme les cofibrations triviales injectives et les
cofibrations triviales locales coïncident, le foncteur $H : \mm \longrightarrow \ch(\nn)$ détermine un morphisme
de $\modq$ si et seulement s'il préserve les cofibrations entre cofibrants \cite[8.5.4]{Hi}. Cela nécessite évidemment
que $H$ préserve les cofibrants. Réciproquement, l'image par $H$ d'une cofibration entre cofibrants
$f :M \rightarrow M'$ de $\mm$ est représenté par le carré extérieur du diagramme commutatif
$$\xymatrix{
F(M) \ar[dd]_{\tau^{H}_{M}} \ar[rr]^{F(f)} & & F(M') \ar[dl] \ar[dd]^{\tau^{H}_{M'}} \\
 & F(M') \amalg_{F(M)} G(M) \ar[dr] &   \\
G(M) \ar[ru] \ar[rr]_{G(f)} &  & G(M') }$$
où $F(f)$ et $G(f)$ sont des cofibrations entre cofibrants et $\tau^{H}_{M}$ et $\tau^{H}_{M'}$ des équivalences
faibles entre cofibrants. Il s'en suit \cite[13.1.2]{Hi} que le morphisme
$F(M') \longrightarrow F(M') \amalg_{F(M)} G(M)$ et donc
le morphisme $F(M') \amalg_{F(M)} G(M) \longrightarrow G(M')$ sont des équivalences faibles, ce qui montre que si
$H$ préserve les cofibrants, alors il préserve les cofibrations entre cofibrants.
\med

\begin{corol}\label{c-eq-hq}
Soit $\modqd$ une sous-2-catégorie pleine de $\modq$.
Soient $\gcc$ une 2-catégorie et $\Phi : \modqd \longrightarrow \gcc$ un pseudo-foncteur envoyant
les équivalences de Quillen dans les équivalences de $\gcc$. Soient $\mm, \nn \in \modqd$ tels
que le modèle cylindre $\cyl(\mm)$ de $\mm$ ou le modèle chemin $\ch(\nn)$ de $\nn$ appartienne à $\modqd$.
Alors le foncteur $\Phi_{\mm,\nn} : \modqd(\mm,\nn) \longrightarrow \gcc(\Phi(\mm),\Phi(\nn))$ envoie les
homotopies de Quillen dans les isomorphismes.
\end{corol}

\dem
Soient $F, G : \mm \longrightarrow \nn$ deux morphismes de $\modqd$ et $\tau : F \rightarrow G$ une
homotopie de Quillen. Il faut montrer que $\Phi_{\mm,\nn}(\tau)$ est un isomorphisme.
On suppose que le modèle chemin $\ch(\nn)$ appartient à $\modqd$~: le cas où le modèle cylindre
$\cyl(\mm)$ appartient à $\modqd$ s'obtient par dualité.
Soit $H_{\tau} : \mm \longrightarrow \ch(\nn)$ le morphisme de $\modqd$ correspondant à $\tau$
(cf. proposition \ref{p-hq-hc}).
On note $\alpha : e_{0} \rightarrow e_{1}$ l'homotopie de Quillen tautologique, définie par
$\alpha_{f} = f : e_{0}(f) \rightarrow e_{1}(f)$ pour $f \in \ch(\nn)$, et correspondant au foncteur
identité $Id : \ch(\nn) \rightarrow \ch(\nn)$. On a ainsi $\tau = \alpha \circ H_{\tau}$.
En appliquant le pseudo-foncteur $\Phi$ au diagramme~:
$$\xymatrix{
\mm \ar[rr]^{H_{\tau}} \ar@/^2pc/@<1ex>[rrrr]^{F} \ar@/_2pc/@<-1ex>[rrrr]_{G} & &
\ch(\nn) \ar@/^1pc/[rr]^(.4){e_{0}} \ar@/_1pc/[rr]_(.4){e_{1}} & \Downarrow \text{\scriptsize{$\alpha$}} & \nn }$$
on constate que les 2-isomophismes structuraux du pseudo-foncteur $\Phi$ donnent un isomophisme entre $\Phi(\tau)$
et $\Phi(\alpha) \circ \Phi(H_{\tau})$, au sens où $\Phi(\tau)$ est composé de $\Phi(\alpha) \circ \Phi(H_{\tau})$
et de 2-isomophismes. Il suffit donc de vérifier que $\Phi(\alpha)$ est un isomorphisme.
En appliquant le pseudo-foncteur $\Phi$ au diagramme~:
$$\xymatrix{
\nn \ar[rr]^{p^{*}} \ar@/^2pc/@<1ex>[rrrr]^{Id} \ar@/_2pc/@<-1ex>[rrrr]_{Id} & &
\ch(\nn) \ar@/^1pc/[rr]^(.4){e_{0}} \ar@/_1pc/[rr]_(.4){e_{1}} & \Downarrow \text{\scriptsize{$\alpha$}} & \nn }$$
on obtient un isomorphisme entre $\Phi(\alpha) \circ \Phi(p^{*})$ et la transformation naturelle identité
du foncteur $Id : \Phi(\nn) \rightarrow \Phi(\nn)$, ce qui assure que $\Phi(\alpha) \circ \Phi(p^{*})$ est un
isomorphisme. Comme $p^{*}$ est une équivalence de Quillen, $\Phi(p^{*})$ est, par hypothèse, une équivalence,
ce qui implique que $\Phi(\alpha)$ est un isomorphisme.
\med

\begin{rem}\label{r-ch-gen}
Les résultats de cette section s'adaptent immédiatement aux 2-catégories $\modqpt$ et $\modqst$ des modèles de
Quillen pointés et des modèles de Quillen stables respectivement. En effet, si $\mm$ est un modèle pointé
(resp. stable), il en va de même du modèle chemin $\ch(\mm)$.

Ces résultats s'adaptent également à la 2-catégorie $\vv$-$\modq$ des $\vv$-catégories modèles, où $\vv$ est
une catégorie modèle monoïdale \cite{Ho1}. Il suffit de vérifier que pour tout $\vv$-modèle $\mm$, le modèle chemin
$\ch(\mm)$ est encore un $\vv$-modèle. Or, si $V \rightarrow W$ est une cofibration de $\vv$ et $X \rightarrow Y$
une fibration de $\ch(\mm)$, il est clair que le morphisme $\{W,X\} \rightarrow \{V,X\} \times_{\{V,Y\}} \{W,Y\}$,
où $\{-,-\}$ désigne la cotensorisation, est une fibration de $\ch(\mm)$, triviale si $V \rightarrow W$ est une
cofibration triviale. Et lorsque $X \rightarrow Y$ est aussi une $e_{0}$-équivalence, l'évaluation en $0$ de
$\{W,X\} \rightarrow \{V,X\} \times_{\{V,Y\}} \{W,Y\}$ n'est autre que
$\{W,X_{0}\} \rightarrow \{V,X_{0}\} \times_{\{V,Y_{0}\}} \{W,Y_{0}\}$, qui est bien une équivalence faible.
\end{rem}

\subsection{Modèles combinatoires}

Une catégorie modèle $\mm$ est dite combinatoire si la catégorie $\mm$ est localement présentable
et la catégorie modèle $\mm$ est engendrée par cofibrations \cite{D1}.
On note $\modqc$ la sous-2-catégorie pleine de $\modq$ constituée des modèles de Quillen combinatoires.

\subsubsection{Modèles présentables et résolutions}

Pour toute catégorie modèle combinatoire $\mm$ et toute petite catégorie $C$, la structure de catégorie modèle
projective sur $\mm^{C}$ existe et est combinatoire. Elle est également propre à gauche dès que $\mm$ est
propre à gauche.

Pour toute catégorie modèle combinatoire propre à gauche $\mm$ et tout ensemble de morphisme $S$ de $\mm$,
la localisation de Bousfield à gauche de $\mm$ relativement à $S$, notée $\mm/S$, existe et est combinatoire
propre à gauche \cite{D1,D2}.

Pour toute petite catégorie $C$, on note $\uu(C)$ la catégorie $\simp^{C^{op}}$ des préfaisceaux en ensembles
simplicaux sur $C$, que l'on munie de la structure de modèle projective. Il s'agit d'une catégorie modèle
combinatoire simpliciale propre.

\begin{defi} \label{def-pres} \cite{D1}
Une catégorie modèle de la forme $\uu(C)/S$, où $C$ est une petite catégorie et $S$ un ensemble de morphisme
de $\uu(C)$, est dite \emph{présentable}.
Une \emph{petite présentation} d'une catégorie modèle $\mm$ est la donnée d'une catégorie modèle présentable
$\uu(C)/S$ et d'une équivalence de Quillen $\uu(C)/S \longrightarrow \mm$.
\end{defi}

Un premier résultat important concernant les catégories modèles combinatoires et présentables
est le théorème suivant \cite[Theorem 1.1]{D2}.

\begin{theorem} \label{th-res-pres} \cite{D2}
Toute catégorie modèle combinatoire possède une petite présentation.  \carre
\end{theorem}

\medskip

De \cite{D2}, on tire également le résultat suivant, qui permet d'utiliser le corollaire \ref{c-eq-hq}.

\begin{prop} \label{p-cyl-comb}
Si $\mm$ est une catégorie modèle combinatoire propre à gauche, alors le modèle cylindre $\cyl(\mm)$ de $\mm$
est également une catégorie modèle combinatoire propre à gauche.
\end{prop}

\dem
La catégorie modèle $\cyl(\mm)$ est la localisation de Bousfield à gauche relativement aux $e_{1}$-équivalences
de la structure de catégorie modèle (combinatoire propre à gauche) projective de $\mm^{[1]}$. Les foncteurs de
l'adjonction de Quillen $e_{1} : \mm^{[1]} \rightleftarrows \mm : p^{*}$ préservent les équivalences faibles,
si bien que la co-unité de l'adjonction dérivée est un isomorphisme. Il découle alors de \cite[Proposition 3.2]{D2}
qu'il existe un ensemble $S$ de morphisme de $\mm^{[1]}$ tel que les $S$-équivalences coïncident avec les
$e_{1}$-équivalences. La catégorie modèle $\cyl(\mm) = \mm^{[1]}/S$ est donc combinatoire propre à gauche.
\med

\medskip

On termine cette section par une remarque générale qui sera utilisée dans la section suivante.

\begin{lem} \label{l-smod-loc}
Soient $\mm$ une catégorie modèle et $\nn$ une sous-catégorie pleine de $\mm$ munie de la structure induite.
On suppose que $\nn$ est stable relativement aux factorisations, i.e. que toutes les
factorisations ``cofibration trivale / fibration'' et ``cofibration / fibration trivale'' d'un
morphisme de $\nn$ sont également dans $\nn$.

Soit $\overline{\nn}$ l'image pleine de $\nn \hookrightarrow \mm \stackrel{\loc}{\rightarrow} \Ho(\mm)$.
Le foncteur $\loc : \nn \hookrightarrow \overline{\nn}$ est une localisation relativement aux équivalences
faibles de $\nn$.
\end{lem}

La stabilité relativement aux factorisations de l'énoncé est assurée dès que la sous-catégorie pleine $\nn$
est stable par équivalence faible dans $\mm$, i.e.~: dès que, pour toute équivalence faible $M \rightarrow M'$
de $\mm$, $M$ appartient à $\nn$ si et seulement si $M'$ appartient à $\nn$.

\medskip

\dem
Il s'agit de vérifier que le foncteur $\loc : \nn \hookrightarrow \overline{\nn}$ possède la propriété universelle
de la localisation. Pour cela, on reprend la construction classique de la catégorie homotopique d'une catégorie
modèle (e.g. \cite[8.3.5]{Hi}). Par définition de $\overline{\nn}$, on a
$$Ob(\overline{\nn}) = Ob(\nn) \quad \mbox{et} \quad \overline{\nn}(N,N') = \mm(N^{cf},N'^{cf})/\!\!\sim$$
où $i_{N} : N^{c} \rightarrow N$ et $j_{N^{c}} : N^{c} \rightarrow N^{cf}$ sont respectivement la
résolution cofibrante de $N$ et la résolution fibrante de $N^{c}$, obtenues par factorisation.
Le quotient est relatif à la relation d'homotopie : on a $f \sim f' \in \mm(N^{cf},N'^{cf})$ s'il existe
un objet cylindre dans $\mm$~:
$$N^{cf} \amalg N^{cf} \stackrel{i \amalg i'}{\longrightarrow} Cyl(N^{cf}) \rightarrow N^{cf}$$
constitué d'une cofibration suivie d'une fibration triviale, et un morphisme
$h : Cyl(N^{cf}) \rightarrow N'^{cf}$ tel que $h \circ i = f$ et $h \circ i' = f'$.
Comme $N^{cf}$ est cofibrant, on a en particulier une factorisation ``cofibration / fibration trivale''
$N^{cf} \stackrel{i}{\longrightarrow} Cyl(N^{cf}) \rightarrow N^{cf}$ de l'identité.
L'hypothèse de stabilité relativement aux factorisations montre alors que $N^{cf}$, $N'^{cf}$ et $Cyl(N^{cf})$
sont dans $\nn$. Ainsi, on a en fait $\overline{\nn}(N,N') = \nn(N^{cf},N'^{cf})/\!\!\sim$ , où la
relation d'homotopie est formée par des objets cylindres dans $\nn$.
Cela suffit pour vérifié la propriété universelle requise de la manière habituelle (e.g. \cite[8.3.5]{Hi}).
\med

\subsubsection{Invariance homotopique}

Soient $\mm$ une catégorie modèle et $C$ une petite catégorie. On note $\Delta$ la catégorie simpliciale et
$c\mm = \mm^{\Delta}$ la catégorie des objets cosimpliciaux dans $\mm$. On a une équivalence de catégories~:
$$\rr : \coc(\uu(C),\mm) \rightleftarrows c\mm^{C} : \qq$$ entre la catégorie des foncteurs
cocontinus de $\uu(C)$ dans $\mm$ et la catégorie des foncteurs de $C$ dans $c\mm$, avec,
pour $F \in \coc(\uu(C),\mm)$ et $\Gamma \in c\mm^{C}$~:
$$\rr(F) = F \circ h^{C \times \Delta}   \quad \text{et} \quad
\qq(\Gamma) : X \mapsto \int^{C \times \Delta} \Gamma^{n}(c) \cdot X_{n}(c)$$
où $h^{C \times \Delta} : C \times \Delta \hookrightarrow \uu(C)$ désigne le foncteur de Yoneda et
$\Gamma^{n}(c) \cdot X_{n}(c) = \amalg_{X_{n}(c)} \Gamma^{n}(c)$.

\smallskip

On munit $c\mm$ de la structure de catégorie modèle de Reedy.
Soit $E_{0} : c\mm \rightleftarrows \mm : E_{0}^{\flat}$ l'adjonction de Quillen constituée de l'évaluation
en $0$ et du foncteur ``objet cosimplicial constant''. On note
$E_{0*} : c\mm^{C} \rightleftarrows \mm^{C} : E_{0*}^{\flat}$ l'adjonction induite par postcomposition.
Il s'agit d'une adjonction de Quillen lorsque $c\mm^{C}$ et $\mm^{C}$ sont munies des structures de catégorie
modèle projectives, dont l'existence est assurée dès que le modèle $\mm$ est combinatoire.

\begin{defi} \label{def-res-cos}
On désigne par $\cres(C,\mm)$ la sous-catégorie pleine de $c\mm^{C}$ constituée des objets $X$
tels que l'unité d'adjonction $\eta_{X} : X \rightarrow (E_{0*}^{\flat} \circ E_{0*})(X)$ est
une équivalence faible argument par argument dans $c\mm^{C}$.
On désigne par $\cresc(C,\mm)$ la sous-catégorie pleine de $\cres(C,\mm)$ constiuée des objets
cofibrants injectifs dans $c\mm^{C}$, i.e. Reedy-cofibrants argument par argument.
\end{defi}

\rem La catégorie $\cresc(C,\mm)$ est la catégorie des résolutions cosimpliciales des foncteurs
de $C$ dans $\mm$ \cite{Hi}.

\medskip

\begin{prop}\label{p-fctq-resc}
Soit $\mm$ une catégorie modèle et $C$ une petite catégorie.
Les foncteurs décrits ci-dessus définissent une équivalence de catégories~:
$$\rr : \modq(\uu(C),\mm) \rightleftarrows \cresc(C,\mm) : \qq$$
qui fait correspondre aux homotopies de Quillen les équivalences faibles argument par argument.
\end{prop}

\dem
Il s'agit d'abord de vérifier que les foncteurs de l'équivalence de catégorie
$\rr : \coc(\uu(C),\mm) \rightleftarrows c\mm^{C} : \qq$ se restreignent aux sous-catégories de l'énoncé.

Soit $F \in \modq(\uu(C),\mm)$. L'unité d'adjonction
$\eta_{\rr(F)} : \rr(F) \rightarrow (E_{0*}^{\flat} \circ E_{0*})(\rr(F))$, évaluée en $c \in C$ et en
$n \in \Delta$, n'est autre que le morphisme $F(h^{C \times \Delta}_{c,n}) \rightarrow F(h^{C \times \Delta}_{c,0})$
induit par la projection $h^{C \times \Delta}_{c,n} = h^{C}_{c} \times \Delta^{n} \rightarrow h^{C}_{c}$.
Cette dernière est une équivalence faible entre cofibrants, donc $\eta_{\rr(F)}$ également, et on a ainsi
$\rr(F) \in \cres(C,\mm)$.
De plus, le foncteur induit $F^{\Delta} : \uu(C)^{\Delta} \rightarrow \mm^{\Delta}$ préserve les Reedy-cofibrants.
Comme, pour tout $c \in C$, $h^{C}_{c} \times \Delta^{\bullet}$ est Reedy-cofibrant dans $\uu(C)^{\Delta}$,
il en va de m\^{e}me de $\rr(F)(c) = F(h^{C}_{c} \times \Delta^{\bullet})$, si bien que $\rr(F) \in \cresc(C,\mm)$.

Soit $\Gamma \in \cresc(C,\mm)$. Comme, pour tout $c \in C$, $\Gamma(c)$ est Reedy-cofibrant,
\cite[16.5.4(2)]{Hi} assure que le foncteur de $\mm$ dans $\simp$ défini par $M \mapsto \mm(\Gamma(c),M)$
préserve fibrations et fibrations triviales, ce qui montre que $\qq(\Gamma)^{\flat}: \mm \rightarrow \uu(C)$
est un foncteur de Quillen à droite, soit $\qq(\Gamma) \in \modq(\uu(C),\mm)$.

Soit $\tau : F_{0} \rightarrow F_{1}$ une homotopie de Quillen dans $\modq(\uu(C),\mm)$. Comme, pour tout
$c \in C$ et $n \in \Delta$, $h^{C \times \Delta}_{c,n}$ est cofibrant dans $\uu(C)$, le morphisme
$\rr(\tau)(c)^{n} =
\tau_{h^{C \times \Delta}_{c,n}} : F_{0}(h^{C \times \Delta}_{c,n}) \longrightarrow F_{1}(h^{C \times \Delta}_{c,n})$
est une équivalence faible de $\mm$ et donc $\rr(\tau)$ est une équivalence faible de $\cresc(C,\mm)$.

Soit $f : \Gamma \rightarrow \Gamma'$ une équivalence faible de $\cresc(C,\mm)$. Par \cite[16.5.5(1)]{Hi},
l'équivalence faible entre résolutions cosimpliciales $f(c) : \Gamma(c) \rightarrow \Gamma'(c)$ induit,
pour tout $M \in \mm$ fibrant, une équivalence faible $\mm(\Gamma'(c),M) \rightarrow \mm(\Gamma(c),M)$.
Cela signifie que la transformation naturelle $\qq(\Gamma')^{\flat} \rightarrow \qq(\Gamma)^{\flat}$,
conjuguée de $\qq(f)$, est une équivalence faible sur tout fibrant de $\mm$, et donc que $\qq(f)$ est une
homotopie de Quillen.
\med

\medskip

Soient $C \in \cat$ et $S \subset \uu(C)$ un ensemble de morphisme.
Soit $T : \uu(C) \longrightarrow \uu(C)/S$ le morphisme de la localisaion de Bousfield.
Le foncteur $T^{*} : \modq(\uu(C)/S,\mm) \longrightarrow \modq(\uu(C),\mm)$, obtenu par précomposition par $T$,
est l'inclusion de la sous-catégorie pleine des morphismes $F : \uu(C) \rightarrow \mm$ tels que $\Lb(F)$ envoie
$S$ dans les isomorphismes.

On note encore $h^{C} \in \uu(C)^{C}$ le foncteur $C \hookrightarrow \ens^{C^{op}} \hookrightarrow \uu(C)$
associé au foncteur de Yoneda, et on considère le foncteur
$F^{C} : \uu(C)^{C} \longrightarrow \mm^{C}$ associé à $F \in \modq(\uu(C),\mm)$.

\begin{prop} \label{p-plgt-eq-loc}
Soient $\mm$ une catégorie modèle combinatoire, $C$ une petite catégorie et $S \subset \uu(C)$ un ensemble de morphisme.
Les foncteurs $T^{*} : \modq(\uu(C)/S,\mm) \hookrightarrow \modq(\uu(C),\mm)$ et $F \mapsto F^{C}(h^{C})$~:
$\modq(\uu(C),\mm) \rightarrow \mm^{C}$ induisent un plongement et une équivalence
$$\xymatrix{ \modq(\uu(C)/S,\mm)[\qq^{-1}] \ar@{^{(}->}[r] & \modq(\uu(C),\mm)[\qq^{-1}] \ar[r]^(.65){\sim} & \Ho(\mm^{C}) }$$
au niveau des catégories localisées relativement aux homotopies de Quillen et aux équivalences faibles
projectives respectivement (en particulier ces catégories sont localement petites).
\end{prop}

\dem
Le composé de l'énoncé est induit par le composé suivant~:
$$\modq(\uu(C)/S,\mm) \hookrightarrow \modq(\uu(C),\mm) \stackrel{\sim}{\rightarrow} \cresc(C,\mm) \hookrightarrow
\cres(C,\mm) \hookrightarrow c\mm^{C} \stackrel{E_{0*}}{\longrightarrow} \mm^{C}$$
dont le second foncteur est l'équivalence $\rr$ de la proposition \ref{p-fctq-resc}.

La sous-catégorie pleine $\cres(C,\mm)$ de $c\mm^{C}$ est constituée des objets $X \in c\mm^{C}$ tels que l'unité
$\eta_{X} : X \rightarrow (E_{0*}^{\flat} \circ E_{0*})(X)$ est une équivalence faible. Comme les foncteurs
$E_{0*}$ et $E_{0*}^{\flat}$ préservent les équivalences faibles, l'inclusion $\cres(C,\mm) \subset c\mm^{C}$
est stable par équivalence faible.

La sous-catégorie pleine $\cresc(C,\mm)$ de $\cres(C,\mm)$ est constituée des objets cofibrants injectifs
(i.e.~: Reedy-cofibrant argument par argument). Comme les cofibrations projectives sont également injectives,
toute cofibration projective de source un cofibrant injectif a pour but un cofibrant injectif, ce qui garantit
la stabilité relativement aux factorisations de $\cresc(C,\mm) \subset \cres(C,\mm)$.

Via l'équivalence $\rr : \modq(\uu(C),\mm) \rightarrow \cresc(C,\mm)$, on constate donc que le plongement
$\modq(\uu(C),\mm) \hookrightarrow c\mm^{C}$ vérifie les hypothèses de stabilité du lemme \ref{l-smod-loc}.
Comme la sous-catégorie pleine $\modq(\uu(C)/S,\mm)$ de $\modq(\uu(C),\mm)$ est stable par homotopies de Quillen,
le plongement $\modq(\uu(C)/S,\mm) \hookrightarrow c\mm^{C}$ vérifie également les hypothèses de stabilité du
lemme \ref{l-smod-loc}.

On voit donc que les foncteurs de la chaine ci-dessus, à l'exception du dernier, induisent des plongements
au niveau des catégories localisées, ce qui donne, d'une part, le plongement de l'énoncé et, d'autre part,
le plongement central du diagramme~:
$$\modq(\uu(C),\mm)[\qq^{-1}] \stackrel{\sim}{\rightarrow} \Ho(\cresc(C,\mm)) \hookrightarrow \Ho(\cres(C,\mm))
\stackrel{E_{0*}}{\longrightarrow} \Ho(\mm^{C})$$
dont la première flèche est l'équivalence induite par $\rr$. Ce plongement central est également essentiellement
surjectif, et donc une équivalence, du fait que tout objet de $\cresc(C,\mm)$ possède une résolution cofibrante
projective et donc injective. Enfin, par définition de $\cres(C,\mm)$, l'adjonction
$E_{0} : c\mm \rightleftarrows \mm : E_{0}^{\flat}$ se restreint en une adjonction
$E_{0} : \cres(C,\mm) \rightleftarrows \mm : E_{0}^{\flat}$ de foncteurs préservant les équivalences faibles,
dont la co-unité est un isomorphisme et l'unité une équivalence faible, et qui induit donc une équivalence de
catégories localisées.
\med

\medskip

Tout morphisme $G : \mm \rightarrow \mm'$ de $\modq$ induit, pour tout $\nn \in \modq$, des foncteurs
$$G_{*} : \modq(\nn,\mm) \longrightarrow \modq(\nn,\mm') \quad \text{et} \quad
G^{*} : \modq(\mm',\nn) \longrightarrow \modq(\mm,\nn)$$
par postcomposition $F \mapsto G \circ F$ et par précomposition $F \mapsto F \circ G$ respectivement.
Il est clair que chacun de ces foncteurs préserve les homotopies de Quillen. De plus, les isomorphismes
$\Lb(G_{*}(F)) \simeq \Lb(G) \circ \Lb(F)$ et $\Lb(G^{*}(F)) \simeq \Lb(F) \circ \Lb(G)$ montrent que si
$G$ est une équivalence de Quillen, les foncteurs $G_{*}$ et $G^{*}$ reflètent également les homotopies
de Quillen.

\begin{prop} \label{p-inv-htp}
Soit une équivalence de Quillen $G : \mm \rightarrow \nn$ dans $\modqc$. Pour tout modèle présentable $\pp$
de $\modqc$, le foncteur $\tilde{G}_{*} : \modq(\pp,\mm)[\qq^{-1}] \longrightarrow \modq(\pp,\nn)[\qq^{-1}]$
est une équivalence de catégorie.
\end{prop}

\dem
Soient une petite catégorie $C$ et un ensemble de morphisme $S \subset \uu(C)$. On considère le diagramme
comutatif~:
$$\xymatrix{
\modq(\uu(C)/S,\mm)[\qq^{-1}] \ar[d]_{\tilde{G}_{*}} \ar@{^{(}->}[r] & \modq(\uu(C),\mm)[\qq^{-1}]
\ar[d]_{\tilde{G}_{*}} \ar[r]^(.65){\sim} & \Ho(\mm^{C}) \ar@<-1ex>[d]_{\Lb(G^{C})}   \\
\modq(\uu(C)/S,\nn)[\qq^{-1}] \ar@{^{(}->}[r] & \modq(\uu(C),\nn)[\qq^{-1}] \ar[r]^(.65){\sim} & \Ho(\nn^{C})
\ar@<-1ex>[u]_{\Rb(G^{C\flat})} }$$
dans lequel, d'après la proposition \ref{p-plgt-eq-loc}, les flèches horizontales de gauche sont des
plongements et les flèches horizontales de droite sont des équivalences.
Comme $G : \mm \rightarrow \nn$ dans $\modq$ est une équivalence de Quillen, le foncteur $\Lb(G^{C})$ est
une équivalence, ainsi donc que
$\tilde{G}_{*} : \modq(\uu(C),\mm)[\qq^{-1}] \longrightarrow \modq(\uu(C),\nn)[\qq^{-1}]$. Cela implique, d'une part,
que le foncteur $\tilde{G}_{*} : \modq(\uu(C)/S,\mm)[\qq^{-1}] \longrightarrow \modq(\uu(C)/S,\nn)[\qq^{-1}]$
est pleinement fidèle, et, d'autre part, que, pour tout $F' \in \modq(\uu(C)/S,\nn)[\qq^{-1}]$, il existe
$F \in \modq(\uu(C),\mm)[\qq^{-1}]$ et un isomorphisme $G_{*}(F) \simeq F'$. Comme
$\Lb(F') \simeq \Lb(G_{*}(F)) \simeq \Lb(G) \circ \Lb(F)$ et $\Lb(G)$ est une équivalence, le foncteur $\Lb(F)$
envoie $S$ dans les isomorphismes, si bien que $F \in \modq(\uu(C)/S,\nn)$, et le foncteur
$\tilde{G}_{*} : \modq(\uu(C)/S,\mm)[\qq^{-1}] \longrightarrow \modq(\uu(C)/S,\nn)[\qq^{-1}]$
est ainsi essentiellement surjectif.
\med

\subsection{Pseudo-localisation de $\modqc$}
\label{sect-ps-loc-modqc}

L'objectif principal de cette section est de construire une pseudo-localisation relativement aux équivalences
de Quillende de $\modqc$ et, plus généralement, de certaines sous-2-catégories pleines de cette dernière.

\smallskip

Pour toute la section, on se donne une sous-2-catégorie pleine $\modqd$ de $\modqc$ satisfaisant les conditions
suivantes~:
\begin{itemize}
    \item tout modèle $\mm \in \modqd$ possède une petite présentation dans $\modqd$,
    \item pour tout modèle présentable $\mm \in \modqd$, le modèle cylindre $\cyl(\mm)$ appartient à $\modqd$.
\end{itemize}
Ces conditions sont évidemment vérifiées par $\modqc$.

\smallskip

Pour chaque $\mm \in \modqd$, on fixe une petite présentation $R_{\mm} : \mm^{p} \rightarrow \mm$,
en choisissant l'identité lorsque $\mm$ est présentable.
Pour tout $\nn \in \modqd$ présentable, par la proposition \ref{p-inv-htp}, le foncteur
$R_{\mm*} : \modq(\nn,\mm^{p}) \longrightarrow \modq(\nn,\mm)$ induit une équivalence de catégorie
$\tilde{R}_{\mm*} : \modqd(\nn,\mm^{p})[\qq^{-1}] \stackrel{\sim}{\longrightarrow} \modqd(\nn,\mm)[\qq^{-1}]$.
On fixe un quasi-inverse $\tilde{R}_{\mm*}^{\flat}$ de $\tilde{R}_{\mm*}$ et des isomorphismes naturels
$$\uni^{\mm} : Id \stackrel{\sim}{\longrightarrow} \tilde{R}_{\mm*}^{\flat} \circ \tilde{R}_{\mm*}
\qquad \text{et} \qquad
\cou^{\mm} : \tilde{R}_{\mm*} \circ \tilde{R}_{\mm*}^{\flat} \stackrel{\sim}{\longrightarrow} Id$$
tels que le quadruplet $(\tilde{R}_{\mm*}, \tilde{R}_{\mm*}^{\flat}, \uni^{\mm}, \cou^{\mm})$
constitue une adjonction. Lorsque $R_{\mm} = Id$, on choisit le foncteur $R_{\mm*}^{\flat} = Id$ et les
transformations naturelles $\uni^{\mm} = \cou^{\mm} = Id$.

\subsubsection{La 2-catégorie $\hmodqd$}

On définit d'abord une 2-catégorie $\hmodqd$ de la fa\c{c}on suivante.
On pose $Ob(\hmodqd) = Ob(\modqd)$ et, pour $\mm, \nn \in \hmodqd$~:
$$\hmodqd(\mm,\nn) = \modqd(\mm^{p},\nn^{p})[\qq^{-1}]$$
la localisation de la catégorie $\modqd(\mm^{p},\nn^{p})$ relativement aux homotopies de Quillen,
qui est une catégorie localement petite par la proposition \ref{p-plgt-eq-loc}.

La composition est donnée par la flèche diagonale ci-dessous.
Les deux flèche verticales supérieures sont les foncteurs localisations.
Comme les composées ``horizontales'' d'homotopies de Quillen sont des homotopies de Quillen, la propriété
universelle induit le foncteur horizontal inférieur.
$$\xymatrix{
\modqd(\mm_{1}^{p},\mm_{2}^{p}) \times \modqd(\mm_{2}^{p},\mm_{3}^{p}) \ar[d] \ar[r]^(.6){\circ} &
\modqd(\mm_{1}^{p},\mm_{3}^{p}) \ar[d]  \\
\left(\modqd(\mm_{1}^{p},\mm_{2}^{p}) \times \modqd(\mm_{2}^{p},\mm_{3}^{p})\right)[(\qq \times \qq)^{-1}]
\ar[d]^{\wr} \ar[r] &  \hmodqd(\mm_{1},\mm_{3})    \\
\hmodqd(\mm_{1},\mm_{2}) \times \hmodqd(\mm_{2},\mm_{3}) \ar[ur]_(.6){\circ} &      }$$
L'identité $Id_{\mm} \in \hmodqd(\mm,\mm)$ est définie par $Id_{\mm^{p}} \in \modqd(\mm^{p},\mm^{p})$.
Les axiomes d'associativité et d'unité découlent de la propriété universelle de localisation.

\subsubsection{Le pseudo-foncteur $\Gamma$}

On définit ensuite un pseudo-foncteur $\Gamma : \modqd \longrightarrow \hmodqd$ comme suit.
Pour les objets, on pose~: $\Gamma = Id : Ob(\modqd) \rightarrow Ob(\hmodqd)$.
Pour tout $\mm, \nn \in \modqd$, le foncteur
$\Gamma_{\mm,\nn} : \modqd(\mm,\nn) \longrightarrow \hmodqd(\mm,\nn)$
est défini par la composée~:
$$\modqd(\mm,\nn) \longrightarrow \modqd(\mm,\nn)[\qq^{-1}] \longrightarrow \modqd(\mm^{p},\nn)[\qq^{-1}]
\longrightarrow \modqd(\mm^{p},\nn^{p})[\qq^{-1}]$$
où la première flèche est le foncteur localisation $\loc$, la deuxième, $\tilde{R}_{\mm}^{*}$, est la précompostion
par $\loc R_{\mm}$, et la troisième est le foncteur $\tilde{R}_{\nn*}^{\flat}$. Pour $F \in \modqd(\mm,\nn)$,
on a donc~:
$$\Gamma_{\mm,\nn}(F) = (\tilde{R}_{\nn*}^{\flat} \circ \tilde{R}_{\mm}^{*})(\loc F)
= \tilde{R}_{\nn*}^{\flat}(\loc F \circ \loc R_{\mm})$$
En particulier, lorsque $\mm$ et $\nn$ sont présentables, on a $\Gamma_{\mm,\nn}(F) = \loc(F)$.
Par la suite, on s'abstiendra souvent d'indiquer la présence du foncteur $\loc$, pour des raisons de lisibilité.
On utilisera souvent le fait que~: $R_{\mm}^{*}(G) \circ F = G \circ R_{\mm} \circ F = G \circ R_{\mm*}(F)$.

\smallskip

Soient $R_{i} : \mm_{i}^{p} \rightarrow \mm_{i}$, $i = 1,2,3,4$, les petites présentations fixées et
$(\tilde{R}_{i*}, \tilde{R}_{i*}^{\flat}, \uni^{i}, \cou^{i})$ le quadruplets fixés associés aux l'équivalences~:
$$\tilde{R}_{i*} : \modqd(\mm_{i-1}^{p},\mm_{i}^{p})[\qq^{-1}] \stackrel{\sim}{\rightleftarrows}
\modqd(\mm_{i-1}^{p},\mm_{i})[\qq^{-1}] : \tilde{R}_{i*}^{\flat}$$
Le 2-isomorphisme
$\un^{\Gamma}_{\mm_{1}} : Id_{\Gamma(\mm_{1})} \stackrel{\sim}{\longrightarrow} \Gamma(Id_{\mm_{1}})$
est défini par $$\un^{\Gamma}_{\mm_{1}} = \uni^{1}_{Id_{\mm_{1}}} :
Id_{\mm_{1}} \stackrel{\sim}{\longrightarrow} (\tilde{R}_{1*}^{\flat} \circ \tilde{R}_{1*})(Id_{\mm_{1}})$$
et le 2-isomorphisme
$\com^{\Gamma}_{G,F} : \Gamma(G) \circ \Gamma(F) \stackrel{\sim}{\longrightarrow} \Gamma(G \circ F)$,
naturel en $F \in \modqd(\mm_{1},\mm_{2})$ et en $G \in \modqd(\mm_{2},\mm_{3})$, est le composé~:
$$\xymatrix{
\tilde{R}_{3*}^{\flat}(G \circ R_{2}) \circ \tilde{R}_{2*}^{\flat}(F \circ R_{1})
\ar[d]^{\uni^{3}_{\left[\tilde{R}_{3*}^{\flat}(G \circ R_{2}) \circ \tilde{R}_{2*}^{\flat}(F \circ R_{1})\right]}}  \\
(\tilde{R}_{3*}^{\flat} \circ \tilde{R}_{3*})
\left(\tilde{R}_{3*}^{\flat}(G \circ R_{2}) \circ \tilde{R}_{2*}^{\flat}(F \circ R_{1})\right)
\ar[d]^{\tilde{R}_{3*}^{\flat}\left(\cou^{3}_{(G \circ R_{2})} \circ \tilde{R}_{2*}^{\flat}(F \circ R_{1})\right)}  \\
\tilde{R}_{3*}^{\flat}\left(G \circ (\tilde{R}_{2*} \circ \tilde{R}_{2*}^{\flat})(F \circ R_{1})\right)
\ar[d]^{\tilde{R}_{3*}^{\flat}\left(G \circ \cou^{2}_{(F \circ R_{1})}\right)}  \\
\tilde{R}_{3*}^{\flat}(G \circ F \circ R_{1})    }$$

Les données qui précèdent constituent un pseudo-foncteur $\Gamma : \modqd \longrightarrow \hmodqd$.
La vérification des axiomes est située à la section \ref{app-gam}.

\begin{prop}\label{prop-gamma-eqQ}
Le pseudo-foncteur $\Gamma : \modqd \longrightarrow \hmodqd$
envoie les équivalences de Quillen dans les équivalences.
\end{prop}

\dem
Soit $F \in \modqd(\mm_{1},\mm_{2})$ une équivalence de Quillen. Le 2-isomorphisme
$$\cou^{2}_{(F \circ R_{1})} : \tilde{R}_{2*}(\Gamma(F)) =
(\tilde{R}_{2*} \circ \tilde{R}_{2*}^{\flat})(F \circ R_{1}) \stackrel{\sim}{\longrightarrow} F \circ R_{1}$$
induit, pour tout $\nn \in \hmodqd$, un 2-isomorphisme
$$\xymatrix{
\hmodqd(\nn,\mm_{1}) = \modqd(\nn^{p},\mm_{1}^{p})[\qq^{-1}] \ar[r]^(.6){\tilde{R}_{1*}} \ar@<-12ex>[d]_{\Gamma(F)_{*}}
& \modqd(\nn^{p},\mm_{1})[\qq^{-1}] \ar[d]^{\tilde{F}_{*}}   \\
\hmodqd(\nn,\mm_{2}) = \modqd(\nn^{p},\mm_{2}^{p})[\qq^{-1}] \ar[r]_(.6){\tilde{R}_{2*}} \ar@{=>}[ur]^{\sim} &
\modqd(\nn^{p},\mm_{2})[\qq^{-1}]       }$$
D'après la proposition \ref{p-inv-htp}, les foncteurs $\tilde{R}_{1*}$, $\tilde{R}_{2*}$ et $\tilde{F}_{*}$
sont des équivalences de catégories. Il en est donc de même de $\Gamma(F)_{*}$, si bien que $\Gamma(F)$ est
une équivalence de $\hmodqd$.
\med

\subsubsection{Propriété universelle}

Cette section est consacrée à la vérification du résultat suivant.

\begin{theorem}\label{th-pu-hmodqd}
Le pseudo-foncteur $\Gamma : \modqd \longrightarrow \hmodqd$ est une pseudo-localisation de $\modqd$
relativement aux équivalences de Quillen.
\end{theorem}

On commence par vérifier une partie de l'énoncé qui précède.

\begin{prop}\label{p-pu-exist}
Soit un pseudo-foncteur $\Phi : \modqd \longrightarrow \gcc$ envoyant les équivalences de Quillen dans les équivalences.
Il existe un pseudo-foncteur $\widetilde{\Phi} : \hmodqd \longrightarrow \gcc$ et un isomorphisme pseudo-naturel
$\alpha : \widetilde{\Phi} \circ \Gamma \stackrel{\sim}{\longrightarrow} \Phi$ tels que, pour tout objet
$\mm \in Ob(\hmodqd) = Ob(\modqd)$, on a $\widetilde{\Phi}(\mm) = \Phi(\mm)$ et $\; \alpha_{\mm} = Id_{\Phi(\mm)}$.
\end{prop}

\dem
Pour tout $\mm \in \modqd$, le 1-morphisme $\Phi(R_{\mm}) : \Phi(\mm^{p}) \rightarrow \Phi(\mm)$ est
une équivalence de $\gcc$. On choisit un quasi-inverse, noté $\Phi(R_{\mm})^{\sharp}$, et des 2-isomorphismes
$$\uni^{\Phi(R_{\mm})} : Id \stackrel{\sim}{\longrightarrow} \Phi(R_{\mm}) \circ \Phi(R_{\mm})^{\sharp} \quad
\text{et} \quad \cou^{\Phi(R_{\mm})} : \Phi(R_{\mm})^{\sharp} \circ \Phi(R_{\mm}) \stackrel{\sim}{\longrightarrow} Id$$
tels que le quadruplet $(\Phi(R_{\mm})^{\sharp}, \Phi(R_{\mm}), \uni^{\Phi(R_{\mm})}, \cou^{\Phi(R_{\mm})})$
constitue une adjonction.

\smallskip

D'après le corollaire \ref{c-eq-hq} et la proposition \ref{p-cyl-comb}, pour tout $\mm, \nn \in \modqd$,
où $\mm$ est présentable, le foncteur $\Phi_{\mm,\nn}$ transforme les homotopies de Quillen en isomorphismes.
On note $\bar{\Phi}_{\mm,\nn}$ l'unique foncteur tel que $\bar{\Phi}_{\mm,\nn} \circ \loc = \Phi_{\mm,\nn}$.
Pour construire $\widetilde{\Phi}$ et $\alpha$, on considère le diagramme (non commutatif) suivant~:
$$\xymatrix{
\modqd(\mm,\nn) \ar[d]_{{R_{\mm}^{*}}} \ar[rr]_{\loc} \ar@/^2pc/[rrrr]^{\Phi_{\mm,\nn}} & & \modqd(\mm,\nn)[\qq^{-1}]
\ar[d]_{{\tilde{R}_{\mm}^{*}}} & &  \gcc(\Phi(\mm),\Phi(\nn))  \\
\modqd(\mm^{p},\nn) \ar[rr]_{\loc} & & \modqd(\mm^{p},\nn)[\qq^{-1}] \ar[d]_{{\tilde{R}_{\nn*}^{\flat}}}
\ar[rr]_(.55){\bar{\Phi}_{\mm^{p},\nn}} & &  \gcc(\Phi(\mm^{p}),\Phi(\nn)) \ar[u]_{(\Phi(R_{\mm})^{\sharp})^{*}}  \\
\modqd(\mm^{p},\nn^{p}) \ar[u]^{R_{\nn*}^{\flat}} \ar[rr]_{\loc} & & \modqd(\mm^{p},\nn^{p})[\qq^{-1}]
\ar[rr]_(.55){\bar{\Phi}_{\mm^{p},\nn^{p}}} & & \gcc(\Phi(\mm^{p}),\Phi(\nn^{p})) \ar[u]_{\Phi(R_{\nn})_{*}}      }$$

On construit le pseudo-foncteur $\widetilde{\Phi} : \hmodqd \longrightarrow \gcc$ de la façon suivante.
Pour les objets, la fonction $\widetilde{\Phi} : Ob(\hmodqd) \longrightarrow Ob(\gcc)$ est égale à
$\Phi : Ob(\modqd) \longrightarrow Ob(\gcc)$. Pour $\mm, \nn \in \modqd$, le foncteur
$\widetilde{\Phi}_{\mm,\nn} : \hmodqd(\mm,\nn) \longrightarrow \gcc(\Phi(\mm),\Phi(\nn))$
est défini par la composée~:
$$\modqd(\mm^{p},\nn^{p})[\qq^{-1}]  \longrightarrow  \gcc(\Phi(\mm^{p}),\Phi(\nn^{p}))  \longrightarrow
\gcc(\Phi(\mm^{p}),\Phi(\nn))  \longrightarrow  \gcc(\Phi(\mm),\Phi(\nn))$$
où la première flèche est le foncteur $\bar{\Phi}_{\mm^{p},\nn^{p}}$,
la deuxième le foncteur $\Phi(R_{\nn})_{*}$, postcomposition par $\Phi(R_{\nn})$,
et la troisième le foncteur $(\Phi(R_{\mm})^{\sharp})^{*}$, précomposition par $\Phi(R_{\mm})^{\sharp}$.
Ainsi, $\widetilde{\Phi}_{\mm,\nn}$ est l'unique foncteur tel que
$$\widetilde{\Phi}_{\pmm,\pnn} \circ \loc =
(\Phi(R_{\pmm})^{-1})^{*} \circ \Phi(R_{\pnn})_{*} \circ \Phi_{\mm^{p},\nn^{p}}$$
et pour tout $F \in \hmodqd(\mm^{p},\nn^{p})$, on a donc~:
$\widetilde{\Phi}(F) = \Phi(R_{\nn}) \circ \bar{\Phi}_{\mm,\nn}(F) \circ \Phi(R_{\mm})^{\sharp}$.

\smallskip

Le 2-isomorphisme
$\un^{\widetilde{\Phi}}_{\mm} : Id_{\widetilde{\Phi}(\mm)} \stackrel{\sim}{\longrightarrow} \widetilde{\Phi}(Id_{\mm})$
est le composé~:
$$Id_{\Phi(\mm)} \stackrel{\sim}{\longrightarrow} \Phi(R_{\mm}) \circ \Phi(R_{\mm})^{\sharp} \stackrel{\sim}{\longrightarrow}
\Phi(R_{\mm}) \circ \Phi_{\mm^{p},\mm^{p}}(Id_{\mm^{p}}) \circ \Phi(R_{\mm})^{\sharp}$$
des 2-isomorphismes $\uni^{\Phi(R_{\mm})}$ et $\Phi(R_{\mm}) \circ \un^{\Phi}_{\mm^{p}} \circ \Phi(R_{\mm})^{\sharp}$.
$$\xymatrix{
\hmodqd(\mm_{1},\mm_{2}) \times \hmodqd(\mm_{2},\mm_{3}) \ar[rr]
\ar@/_2pc/@<-20ex>[dddd]_{\widetilde{\Phi} \times \widetilde{\Phi}} & &
\hmodqd(\mm_{1},\mm_{3}) \ar@/^2pc/@<9ex>[dddd]^{\widetilde{\Phi}} \\
\modqd(\mm_{1}^{p},\mm_{2}^{p}) \times \modqd(\mm_{2}^{p},\mm_{3}^{p})
\ar[d]_{\Phi \times \Phi} \ar[rr] \ar[u]^{\loc \times \loc} & &
\modqd(\mm_{1}^{p},\mm_{3}^{p}) \ar[d]^{\Phi} \ar[u]_{\loc}  \\
\gcc(\Phi(\mm_{1}^{p}),\Phi(\mm_{2}^{p})) \times \gcc(\Phi(\mm_{2}^{p}),\Phi(\mm_{3}^{p}))
\ar[rr] \ar[d]_{(\Phi(R_{2}))_{*} \times (\Phi(R_{3}))_{*}} \ar@{=>}[urr]_(.65){\com^{\Phi}} & &
\gcc(\Phi(\mm_{1}^{p}),\Phi(\mm_{3}^{p})) \ar[d]^{(\Phi(R_{3}))_{*}}  \\
\gcc(\Phi(\mm_{1}^{p}),\Phi(\mm_{2})) \times \gcc(\Phi(\mm_{2}^{p}),\Phi(\mm_{3}))
\ar[d]_{(\Phi(R_{1})^{\sharp})^{*} \times (\Phi(R_{2})^{\sharp})^{*}}  & &
\gcc(\Phi(\mm_{1}^{p}),\Phi(\mm_{3})) \ar[d]^{(\Phi(R_{1})^{\sharp})^{*}}  \\
\gcc(\Phi(\mm_{1}),\Phi(\mm_{2})) \times \gcc(\Phi(\mm_{2}),\Phi(\mm_{3})) \ar[rr]
\ar@{=>}@<-2ex>[uurr]^(.65){\tilde{\cou}} & &  \gcc(\Phi(\mm_{1}),\Phi(\mm_{3}))  \\    }$$
La transformation naturelle $\com^{\widetilde{\Phi}}$ est définie par
$$\com^{\widetilde{\Phi}} = \left[ [(\Phi(R_{1})^{\sharp})^{*} \circ (\Phi(R_{3}))_{*} \circ \com^{\bar{\Phi}} \right]
\ast [\tilde{\cou} \circ (\Phi \times \Phi)]$$
où $\tilde{\cou}$ est la transformation naturelle en $F \in \gcc(\Phi(\mm_{1}^{p}),\Phi(\mm_{2}^{p}))$ et en
$G \in \gcc(\Phi(\mm_{2}^{p}),\Phi(\mm_{3}^{p}))$~:
$$\tilde{\cou}_{G,F} = \Phi(R_{3}) \circ G \circ \cou^{\Phi(R_{2})} \circ F \circ \Phi(R_{1})^{\sharp}$$
et où $\com^{\bar{\Phi}}$ est l'unique transformation naturelle telle que
$\com^{\bar{\Phi}} \circ (\loc \times \loc) = \loc \circ \com^{\Phi}$.

Pour tout $F \in \modqd(\mm_{1}^{p},\mm_{2}^{p})$, $G \in \modqd(\mm_{2}^{p},\mm_{3}^{p})$
et $H \in \modqd(\mm_{3}^{p},\mm_{4})$, l'axiome de composition de $\Phi$ induit, par la propriété universelle,
le diagramme commutatif
$$\xymatrix{
\bar{\Phi}(H) \circ \bar{\Phi}(G) \circ \bar{\Phi}(F)  \ar[d]_{\com^{\bar{\Phi}}_{H,G} \circ \bar{\Phi}(F)}
\ar[rrr]^{\bar{\Phi}(H) \circ \com^{\bar{\Phi}}_{G,F}} & & & \bar{\Phi}(H) \circ \bar{\Phi}(G \circ F)
\ar[d]^{\com^{\bar{\Phi}}_{H,G \circ F}} \\
\bar{\Phi}(H \circ G) \circ \bar{\Phi}(F)  \ar[rrr]_{\com^{\bar{\Phi}}_{H \circ G,F}} & & &
\bar{\Phi}(H \circ G \circ F)   }$$
qui participe à la vérification de l'axiome d'associativité de $\widetilde{\Phi}$.

Le 2-isomorphisme
$\com^{\widetilde{\Phi}}_{G,F} : \widetilde{\Phi}(G) \circ \widetilde{\Phi}(F)
\stackrel{\sim}{\longrightarrow} \widetilde{\Phi}(G \circ F)$,
naturel en $F \in \hmodqd(\mm_{1},\mm_{2})$ et en $G \in \hmodqd(\mm_{2},\mm_{3})$,
est donc le composé~:
$$\xymatrix{
\Phi(R_{3}) \circ \bar{\Phi}(G) \circ \Phi(R_{2})^{\sharp} \circ \Phi(R_{2}) \circ \bar{\Phi}(F) \circ \Phi(R_{1})^{\sharp}
\ar[d]^{\Phi(R_{3}) \circ \bar{\Phi}(G) \circ \cou^{\Phi(R_{2})} \circ \bar{\Phi}(F) \circ \Phi(R_{1})^{\sharp}}  \\
\Phi(R_{3}) \circ \bar{\Phi}(G)  \circ \bar{\Phi}(F) \circ \Phi(R_{1})^{\sharp}
\ar[d]^{\Phi(R_{3}) \circ \com^{\bar{\Phi}}_{G,F}  \circ \Phi(R_{1})^{\sharp}}  \\
\Phi(R_{3}) \circ \bar{\Phi}(G \circ F) \circ \Phi(R_{1})^{\sharp}   }$$

Les données qui précèdent constituent un pseudo-foncteur $\widetilde{\Phi} : \hmodqd \longrightarrow \gcc$.
La vérification des axiomes est située à la section \ref{app-phi}.

\smallskip

On construit maintenant un isomorphisme pseudo-naturelle
$\alpha : \widetilde{\Phi} \circ \Gamma \stackrel{\sim}{\longrightarrow} \Phi$
de la façon suivante.
Pour tout $\mm \in \modqd$, on pose $\alpha_{\mm} = Id_{\Phi(\mm)}$ et pour tout $F \in \modqd(\mm_{1},\mm_{2})$,
le 2-isomorphisme naturel
$\tpn^{\alpha}_{F} : (\widetilde{\Phi} \circ \Gamma)(F) \stackrel{\sim}{\longrightarrow} \Phi(F)$
est défini par la composée~:
$$\xymatrix{
\Phi(R_{2}) \circ \bar{\Phi}\left(\tilde{R}_{2*}^{\flat}(F \circ R_{1})\right) \circ \Phi(R_{1})^{\sharp}
\ar[d]^{\com^{\bar{\Phi}}_{\gamma \! R_{2},R_{2*}^{\flat}(F \circ R_{1})} \circ \Phi(R_{1})^{\sharp}}  \\
\bar{\Phi}\left((\tilde{R}_{2*} \circ \tilde{R}_{2*}^{\flat})(F \circ R_{1})\right) \circ \Phi(R_{1})^{\sharp}
\ar[d]^{\bar{\Phi}\left(\cou^{2}_{(F \circ R_{1})}\right) \circ \Phi(R_{1})^{\sharp}}  \\
\Phi(F \circ R_{1}) \circ \Phi(R_{1})^{\sharp}
\ar[d]^{(\com^{\Phi}_{F,R_{1}})^{-1} \circ \Phi(R_{1})^{\sharp}}  \\
\Phi(F) \circ \Phi(R_{1}) \circ \Phi(R_{1})^{\sharp}
\ar[d]^{\Phi(F) \circ \left(\uni^{\Phi(R_{1})}\right)^{-1}}  \\
\Phi(F)     }$$

Les données qui précèdent constituent une transformation pseudo-naturelle
$\alpha : \widetilde{\Phi} \circ \Gamma \longrightarrow \Phi$.
La vérification des axiomes est située à la section \ref{app-alph}.
\med

\medskip

On en arrive maintenant à la démonstration du principal résultat de cette section.

\medskip

\demr{\ref{th-pu-hmodqd}}
La proposition \ref{prop-gamma-eqQ} indique que le pseudo-foncteur $\Gamma : \modqd \longrightarrow \hmodqd$
envoie les équivalences de Quillen dans les équivalences.
La proposition \ref{p-pu-exist} montre que, pour toute 2-catégorie $\gcc$, le 2-foncteur
$\Gamma^{*} : [\hmodqd,\gcc] \longrightarrow [\modqd,\gcc]_{\qq}$ est essentiellement surjectif.
Il reste donc à vérifier que, pour tout $\Psi, \Psi' \in [\hmodqd,\gcc]$, le foncteur
$\Gamma^{*}_{\Psi, \Psi'} : [\hmodqd,\gcc](\Psi, \Psi') \longrightarrow [\modqd,\gcc](\Gamma^{*}(\Psi),\Gamma^{*}(\Psi'))$
est un isomorphisme.

Soit un 1-morphisme $\chi : \Psi \rightarrow \Psi'$, c'est-à-dire une transformation pseudo-naturelle.
Il résulte du fait que $\Gamma$ est l'identité sur les objets que, pour tout $\mm \in \modqd$, on a
$\Gamma^{*}(\chi)_{\mm} = \chi_{\mm}$, donc que les 1-morphismes structuraux de $\chi$ et $\Gamma^{*}(\chi)$
se déterminent l'un l'autre. D'autre part, pour tout $\mm_{1}, \mm_{2} \in \modqd$, on a
$\tpn^{\Gamma^{*}(\chi)}_{\mm_{1},\mm_{2}} = \tpn^{\chi}_{\mm_{1},\mm_{2}} \circ \Gamma$.

On étudie d'abord le cas où $\mm_{1} = \mm_{1}^{p}$. Par construction, on a un foncteur localisation
$\loc : \modqd(\mm_{1}^{p},\mm_{2}^{p}) \rightarrow \hmodqd(\mm_{1}^{p},\mm_{2})$
et pour tout $G \in \modqd(\mm_{1}^{p},\mm_{2}^{p})$, un isomorphisme
$\uni^{2}_{\gamma G} : \loc G \stackrel{\sim}{\longrightarrow} \tilde{R}_{2*}^{\flat} \tilde{R}_{2*}( \loc G)
= \Gamma(R_{2} G)$.
La naturalité de $\tpn^{\chi}_{\mm_{1}^{p},\mm_{2}}$ relativement à l'isomorphisme $\uni^{2}_{\gamma G}$
$$\xymatrix{
\chi_{\mm_{2}} \circ \Psi(\loc G)
\ar[d]_{\tpn^{\chi}_{\gamma G}} \ar[rr]^(.45){\chi_{\mm_{2}} \circ \Psi(\uni^{2}_{\gamma G})}_(.4){\sim}  & &
\chi_{\mm_{2}} \circ \Psi(\tilde{R}_{2*}^{\flat} \tilde{R}_{2*}(\loc G))
\ar[d]^{\tpn^{\chi}_{\Gamma(R_{2} G)}} \ar@{=}[r] &
\chi_{\mm_{2}} \circ \Psi(\Gamma(R_{2}G))  \ar[d]^{\tpn^{\Gamma^{*}(\chi)}_{R_{2} G}}   \\
\Psi'(\loc G) \circ \chi_{\mm_{1}^{p}} \ar[rr]_(.45){\Psi'(\uni^{2}_{\gamma G}) \circ \chi_{\mm_{1}^{p}}}^(.4){\sim}
& & \Psi'(\tilde{R}_{2*}^{\flat} \tilde{R}_{2*}(\loc G)) \circ \chi_{\mm_{1}^{p}} \ar@{=}[r] &
\Psi'(\Gamma(R_{2}G)) \circ \chi_{\mm_{1}^{p}}    }$$
montre que $\tpn^{\chi}_{\mm_{1}^{p},\mm_{2}}$ est l'unique 2-morphisme tel que
$$\tpn^{\chi}_{\mm_{1}^{p},\mm_{2}} \circ \loc = [\Psi'(\uni^{2} \loc)^{-1} \circ \chi_{\mm_{1}^{p}}] \circ
\left[\tpn^{\Gamma^{*}(\chi)}_{\mm_{1}^{p},\mm_{2}} \circ R_{2*}\right] \circ [\chi_{\mm_{2}} \circ \Psi(\uni^{2} \loc)]$$
et donc que les 2-morphismes structuraux de $\chi$ sont déterminés par ceux de $\Gamma^{*}(\chi)$ lorsque
$\mm_{1} = \mm_{1}^{p}$.
On en déduit, dans le cas général, que, pour tout $F \in \hmodqd(\mm_{1},\mm_{2})$, l'axiome de composition
pour la composée $F \circ \Gamma(R_{1}) : \mm_{1}^{p} \longrightarrow \mm_{1} \longrightarrow \mm_{2}$
$$\xymatrix{
\chi_{\mm_{2}} \circ  \Psi(F) \circ \Psi(\Gamma(R_{1}))  \ar[d]_{Id \ast \com^{\Psi}_{F,\Gamma(R_{1})}}
\ar[rr]^{\tpn^{\chi}_{F} \ast Id} &
& \Psi'(F) \circ \chi_{\mm_{1}} \circ \Psi(\Gamma(R_{1}))  \ar[rr]^{Id \ast \tpn^{\chi}_{\Gamma(R_{1})}} & &
\Psi'(F) \circ \Psi'(\Gamma(R_{1})) \circ \chi_{\mm_{1}^{p}} \ar[d]^{\com^{\Psi'}_{F,\Gamma(R_{1})} \ast Id}   \\
\chi_{\mm_{2}} \circ \Psi(F \circ \Gamma(R_{1})) \ar[rrrr]_{\tpn^{\chi}_{F \circ \Gamma(R_{1})}}  & & & &
\Psi'(F \circ \Gamma(R_{1})) \circ \chi_{\mm_{1}^{p}}   }$$
détermine le 2-morphisme $\tpn^{\chi}_{F} \circ \Psi(\Gamma(R_{1}))$. Comme $\Psi(\Gamma(R_{1}))$ est
une équivalence, on a un isomorphisme
$Id \stackrel{\sim}{\longrightarrow} \Psi(\Gamma(R_{1})) \circ \Psi(\Gamma(R_{1}))^{-1}$ et un diagramme commutatif
$$\xymatrix{
\chi_{\mm_{2}} \circ  \Psi(F)  \ar[d]_{\tpn^{\chi}_{F}} \ar[r]^(.28){\sim}  &
\chi_{\mm_{2}} \circ  \Psi(F) \circ \Psi(\Gamma(R_{1})) \circ \Psi(\Gamma(R_{1}))^{-1}
\ar[d]^{\tpn^{\chi}_{F} \circ \Psi(\Gamma(R_{1})) \circ \Psi(\Gamma(R_{1}))^{-1}}  \\
\Psi'(F) \circ \chi_{\mm_{1}} \ar[r]^(.28){\sim}  &
\Psi'(F) \circ \chi_{\mm_{1}} \circ \Psi(\Gamma(R_{1})) \circ \Psi(\Gamma(R_{1}))^{-1}    }$$
qui montrent que $\tpn^{\chi}_{F}$ est uniquement déterminé par les 2-morphismes structuraux de $\Gamma^{*}(\chi)$.
On a ainsi vérifié que le foncteur $\Gamma^{*}_{\Psi, \Psi'}$ est injectif sur les objets.

Soit un 1-morphisme $\chi': \Gamma^{*}(\Psi) \rightarrow \Gamma^{*}(\Psi')$. On construit un 1-morphisme $\chi$
tel que $\Gamma^{*}(\chi) = \chi'$ de la façon suivante. Pour tout $\mm \in \hmodqd$, on pose
$\chi_{\mm} = \chi'_{\mm}$. Pour définir $\tpn^{\chi}_{\mm_{1},\mm_{2}}$, on procède en deux temps.
Dans le cas où $\mm_{1} = \mm_{1}^{p}$, on définit $\tpn^{\chi}_{\mm_{1}^{p},\mm_{2}}$ comme étant l'unique
2-morphisme tel que
$$\tpn^{\chi}_{\mm_{1}^{p},\mm_{2}} \circ \loc = [\Psi'(\uni^{2} \loc)^{-1} \circ \chi_{\mm_{1}^{p}}] \circ
[\tpn^{\chi'}_{\mm_{1}^{p},\mm_{2}} \circ R_{2*}] \circ [\chi_{\mm_{2}} \circ \Psi(\uni^{2} \loc)]$$
Comme $\Gamma(R_{1})$ est une équivalence, on peut se fixer un quasi-inverse $\Gamma(R_{1})^{\sharp}$
et un isomorphisme
$\uni^{\Gamma(R_{1})} : Id \stackrel{\sim}{\longrightarrow} \Gamma(R_{1}) \circ \Gamma(R_{1})^{\sharp}$.
On en tire un isomorphisme composé $\uni^{\Psi\Gamma(R_{1})}$~:
$$\xymatrix{Id \ar[r]^(.4){\sim}_(.4){\un^{\Psi}} &
\Psi(Id) \ar[rr]^(.35){\sim}_(.35){\Psi(\uni^{\Gamma(R_{1})})} & & \Psi(\Gamma(R_{1}) \circ \Gamma(R_{1})^{\sharp})
\ar[rrr]^{\sim}_{\left(\com^{\Psi}_{\Gamma(R_{1}),\Gamma(R_{1})^{\sharp}}\right)^{-1}} & & &
\Psi(\Gamma(R_{1})) \circ \Psi(\Gamma(R_{1})^{\sharp})        }$$
On définit alors, pour tout $F \in \hmodqd(\mm_{1},\mm_{2})$, le 2-isomorphisme $\tpn^{\chi}_{F}$, comme le
composé des 2-isomorphismes naturels en $F$~:
$$\xymatrix{
\chi_{\mm_{2}} \circ \Psi(F)
\ar[d]^{\chi_{\mm_{2}} \circ \Psi(F) \circ \uni^{\Psi(\Gamma(R_{1}))}}   \\
\chi_{\mm_{2}} \circ \Psi(F) \circ \Psi(\Gamma(R_{1})) \circ \Psi(\Gamma(R_{1})^{\sharp})
\ar[d]^{\chi_{\mm_{2}} \circ \com^{\Psi}_{F,\Gamma(R_{1})} \circ \Psi(\Gamma(R_{1})^{\sharp})}   \\
\chi_{\mm_{2}} \circ \Psi(F \circ \Gamma(R_{1})) \circ \Psi(\Gamma(R_{1})^{\sharp})
\ar[d]^{\tpn^{\chi}_{F \circ \Gamma(R_{1})} \circ \Psi(\Gamma(R_{1})^{\sharp}) }   \\
\Psi'(F \circ \Gamma(R_{1})) \circ \chi_{\mm_{1}^{p}} \circ \Psi(\Gamma(R_{1})^{\sharp})
\ar[d]^{\left(\com^{\Psi'}_{F,\Gamma(R_{1})}\right)^{-1} \circ \chi_{\mm_{1}^{p}} \circ \Psi(\Gamma(R_{1})^{\sharp})}  \\
\Psi'(F) \circ  \Psi(\Gamma(R_{1})) \circ \chi_{\mm_{1}^{p}} \circ \Psi(\Gamma(R_{1})^{\sharp})
\ar[d]^{\Psi'(F) \circ \left(\tpn^{\chi}_{\Gamma(R_{1})}\right)^{-1} \circ \Psi(\Gamma(R_{1})^{\sharp})}  \\
\Psi'(F) \circ \chi_{\mm_{1}} \circ \Psi(\Gamma(R_{1})) \circ \Psi(\Gamma(R_{1})^{\sharp})
\ar[d]^{\Psi'(F) \circ \chi_{\mm_{1}} \circ \left(\uni^{\Psi(\Gamma(R_{1}))}\right)^{-1}}  \\
\Psi'(F) \circ \chi_{\mm_{1}}    }$$

Les données qui précèdent constituent une transformation pseudo-naturelle $\chi : \Psi \rightarrow \Psi'$~:
la vérification des axiomes est située à la section \ref{app-chi}.
Il est clair que $\Gamma^{*}(\chi)_{\mm} = \chi'_{\mm}$ pour tout $\mm \in \modqd$. Il faut vérifier que,
pour tout $F \in \modqd(\mm_{1},\mm_{2})$, on a $\tpn^{\chi}_{\Gamma(F)} = \tpn^{\chi'}_{F}$. Via un isomorphisme
$Id \stackrel{\sim}{\longrightarrow} \Psi(\Gamma(R_{1})) \circ \Psi(\Gamma(R_{1}))^{-1}$ et les axiomes d'associativité
de $\chi$ et $\chi'$, on se ramène au cas où $\mm_{1} = \mm_{1}^{p}$. Soit $F \in \modqd(\mm_{1}^{p},\mm_{2})$,
et soit $G \in \modqd(\mm_{1}^{p},\mm_{2}^{p})$ tel que $\loc G = \Gamma(F) = \tilde{R}_{2*}^{\flat}(\loc F)$.
Dans le diagramme
$$\xymatrix{
\chi_{\mm_{2}} \circ \Psi(\Gamma(F))
\ar[d]_{\tpn^{\chi}_{\Gamma(F)}} \ar[rr]^(.5){\chi_{\mm_{2}} \circ \Psi(\uni^{2}_{\gamma G})}_(.4){\sim}  & &
\chi_{\mm_{2}} \circ \Psi(\Gamma(R_{2} G))  \ar[d]^{\tpn^{\chi'}_{R_{2} G}}
\ar[rrr]^(.55){\chi_{\mm_{2}} \circ \Psi(\tilde{R}_{2*}^{\flat}(\cou^{2}_{\gamma F}))}_(.6){\sim} & & &
\chi_{\mm_{2}} \circ \Psi(\Gamma(F))  \ar[d]^{\tpn^{\chi}_{\Gamma(F)}}   \\
\Psi'(\Gamma(F)) \circ \chi_{\mm_{1}^{p}}
\ar[rr]_(.5){\Psi'(\uni^{2}_{\gamma G}) \circ \chi_{\mm_{1}^{p}}}^(.4){\sim} & &
\Psi'(\Gamma(R_{2} G)) \circ \chi_{\mm_{1}^{p}}
\ar[rrr]_(.55){\Psi'(\tilde{R}_{2*}^{\flat}(\cou^{2}_{\gamma F}))\circ \chi_{\mm_{1}^{p}}}^(.6){\sim} & & &
\Psi'(\Gamma(F)) \circ \chi_{\mm_{1}^{p}}    }$$
le carré gauche est commutatif par définition de $\tpn^{\chi}$. Les composés horizontales sont les flèches
identités, par les identités triangulaires, si bien que le carré droit est également commutatif. Mais le
diagramme
$$\xymatrix{
\chi_{\mm_{2}} \circ \Psi(\Gamma(R_{2} G))  \ar[d]_{\tpn^{\chi'}_{R_{2} G}}
\ar[rrr]^(.5){\chi_{\mm_{2}} \circ \Psi(\tilde{R}_{2*}^{\flat}(\cou^{2}_{\gamma F}))}_{\sim} & & &
\chi_{\mm_{2}} \circ \Psi(\Gamma(F))  \ar[d]^{\tpn^{\chi'}_{F}}   \\
\Psi'(\Gamma(R_{2} G)) \circ \chi_{\mm_{1}^{p}}
\ar[rrr]_(.55){\Psi'(\tilde{R}_{2*}^{\flat}(\cou^{2}_{\gamma F}))\circ \chi_{\mm_{1}^{p}}}^{\sim} & & &
\Psi'(\Gamma(F)) \circ \chi_{\mm_{1}^{p}}    }$$
est aussi commutatif. Cela résulte de la naturalité de $\tpn^{\chi'}_{F}$ et de l'existence d'un zigzag
d'homotopies de Quillen entre $F$ et $R_{2}G$ dans $\modqd(\mm_{1}^{p},\mm_{2})$ réalisant le 2-isomorphisme
$\cou^{2}_{\gamma F}$ via $\loc$. Ceci achève de montrer que $\Gamma^{*}(\chi) = \chi'$.

Ainsi, on a vérifié que le foncteur $\Gamma^{*}_{\Psi, \Psi'}$ est bijectif sur les objets. Il reste à
montrer qu'il est pleinement fidèle.

Soit un 2-morphisme $\theta : \chi \rightarrow \chi'$, c'est-à-dire une modification, de $[\hmodqd,\gcc](\Psi, \Psi')$.
Comme $\Gamma^{*}(\theta)_{\mm} = \theta_{\mm}$ pour tout $\mm \in \modqd$, il est clair que $\Gamma^{*}_{\Psi, \Psi'}$
est fidèle.
Soient $\chi$ et $\chi'$ deux morphismes de $[\hmodqd,\gcc](\Psi, \Psi')$ et un 2-morphisme
$\theta' : \Gamma^{*}(\chi) \rightarrow \Gamma^{*}(\chi')$. On définit un 2-morphisme $\theta : \chi \rightarrow \chi'$
tel que $\Gamma^{*}(\theta) = \theta'$ en posant, pour tout $\mm \in \hmodqd$, $\theta_{\mm} = \theta'_{\mm}$
et en vérifiant la commutativité du diagramme
$$\xymatrix{
\chi_{\mm_{2}} \circ \Psi(F)  \ar[r]^{\tpn^{\chi}_{F}} \ar[d]_{\theta_{\mm_{2}} \ast Id}  &
\Psi'(F) \circ \chi_{\mm_{1}}  \ar[d]^{Id \ast \theta_{\mm_{1}}}   \\
\chi'_{\mm_{2}} \circ \Psi(F)  \ar[r]_{\tpn^{\chi'}_{F}}  &  \Psi'(F) \circ \chi'_{\mm_{1}}      }$$
pour tout $F \in \hmodqd(\mm_{1},\mm_{2})$. A nouveau, via un isomorphisme
$Id \stackrel{\sim}{\longrightarrow} \Psi(\Gamma(R_{1})) \circ \Psi(\Gamma(R_{1}))^{-1}$ et l'axiome d'associativité
de $\chi$ et $\chi'$, on se ramène au cas où $\mm_{1} = \mm_{1}^{p}$. En choisissant
$G \in \modqd(\mm_{1}^{p},\mm_{2}^{p})$ tel que $\loc G = F$ et en utilisant l'isomorphisme
$\uni^{2}_{\gamma G} : \loc G \stackrel{\sim}{\longrightarrow} \tilde{R}_{2*}^{\flat} \tilde{R}_{2*}( \loc G)
= \Gamma(R_{2} G)$, la commutativité découle du fait que $\Gamma^{*}(\theta) = \theta'$ est une modification.
Le foncteur $\Gamma^{*}_{\Psi, \Psi'}$ est ainsi pleinement fidèle, et donc une équivalence.
Ceci achève la démonstration.
\rmed

\subsubsection{Autres propriétés}

On note $\Catad$ la sous-2-catégorie de $\Cat$ ayant les adjonctions pour 1-morphismes.
On a un pseudo-foncteur $\Ho : \modqd \longrightarrow \Catad$ qui associe à une catégorie modèle $\mm \in \modqd$
sa catégorie homotopique $\Ho(\mm)$, et à un morphisme $F : \mm \rightarrow \nn$ de $\modqd$ l'adjonction de
foncteurs dérivés totaux $\Lb(F) : \Ho(\mm) \rightleftarrows \Ho(\nn) : \Rb(F^{\flat})$.

Le pseudo-foncteur $\Ho$ envoie les équivalences de Quillen dans les équivalences, donc, d'après la proposition
\ref{p-pu-exist}, il existe un pseudo-foncteur $\widetilde{\Ho} : \hmodqd \longrightarrow \Catad$ et
un isomorphisme pseudo-naturel $\alpha : \widetilde{\Ho} \circ \Gamma \stackrel{\sim}{\longrightarrow} \Ho$.

\begin{prop}\label{prop-gamma-eqQ-rec}
Un morphisme $F : \mm \rightarrow \nn$ de $\modqd$ est une équivalence de Quillen si et seulement si
$\Gamma(F)$ est une équivalence de $\hmodqd$.
\end{prop}

\dem
La proposition \ref{prop-gamma-eqQ} indique que si $F$ est une équivalence de Quillen, alors
$\Gamma(F)$ est une équivalence.
Réciproquement, si $\Gamma(F)$ est une équivalence, alors $\Ho(F) \simeq \widetilde{\Ho}(\Gamma(F))$
est aussi une équivalence et $F$ est donc une équivalence de Quillen.
\med

\medskip

Deux catégories modèles $\mm$ et $\nn$ sont dites Quillen-équivalentes si elles sont connectées par
un zigzag d'équivalences de Quillen.

\begin{prop}\label{prop-gamma-eq-zigzag-Qeq}
Deux catégories modèles $\mm$ et $\nn$ de $\modqd$ sont Quillen-équivalentes si et
seulement si $\Gamma(\mm)$ et $\Gamma(\nn)$ sont équivalentes dans $\hmodqd$.
\end{prop}

\dem
Si $\Gamma(\mm)$ et $\Gamma(\nn)$ sont équivalentes, il en va de même  $\Gamma(\mm^{p})$ et $\Gamma(\nn^{p})$.
Comme le foncteur
$\Gamma_{\mm^{p},\nn^{p}} = \loc : \modqd(\mm^{p},\nn^{p}) \longrightarrow \hmodqd(\Gamma(\mm^{p}),\Gamma(\nn^{p}))$
est surjectif, il existe un morphisme $F : \mm^{p} \rightarrow \nn^{p}$ qui induit l'équivalence
de $\Gamma(\mm^{p})$ vers $\Gamma(\nn^{p})$, et qui doit donc être une équivalence de Quillen d'après
la proposition \ref{prop-gamma-eqQ-rec}. On obtient ainsi un zigzag d'équivalences de Quillen
$$\xymatrix{ \mm & \mm^{p} \ar[l]_{R_{\mm}} \ar[r]^{F} & \nn^{p} \ar[r]^{R_{\nn}} & \nn   }$$
La réciproque découle immédiatemet de la proposition \ref{prop-gamma-eqQ-rec}.
\med

\bigskip

On précise maintenant la nature des troncations de la pseudo-localisation de $\modqd$ relativement aux équivalences
de Quillen.

\begin{prop} \label{p-tr-modq}
La 1-troncation $\tron(\Gamma) : \tron(\modqd) \longrightarrow \tron(\hmodqd)$ de la pseudo-localisation de
$\modqd$ relativement aux équivalences de Quillen est une localisation de $\tron(\modqd)$ relativement
aux équivalences de Quillen.
\end{prop}

\dem
Cela résulte du théorème \ref{th-pu-hmodqd} et de la proposition \ref{p-tr-ploc}, en notant que le pseudo-foncteur
$\Gamma : \modqd \longrightarrow \hmodqd$ est l'identité sur les objets.
\med

\begin{prop} \label{p-ptr-modq}
La pseudo-troncation $\ptron(\Gamma) : \ptron(\modqd) \longrightarrow \ptron(\hmodqd)$ de la pseudo-localisation
de $\modqd$ relativement aux équivalences de Quillen est une localisation de $\ptron(\modqd)$ relativement
aux équivalences de Quillen.
\end{prop}

\dem
Soient une catégorie $C$ et un foncteur $\Phi : \ptron(\modqd) \longrightarrow C$ envoyant les équivalences de
Quillen dans les isomorphismes. On remarque d'abord que pour tout $\mm, \nn \in \modqd$, on a
$\ptron(\hmodqd)(\mm,\nn) = \ptron(\modqd(\mm^{p},\nn^{p})[\qq^{-1}]) = \modqd_{\circ}(\mm^{p},\nn^{p})/\!\!\!\sim$,
où $F \sim G \in \modqd_{\circ}(\mm^{p},\nn^{p})$ si et seulement si $F$ et $G$ sont Quillen-homotopes.
Par la version cylindrique de la proposition \ref{p-hq-hc}, l'argument habituel montre qu'alors $\Phi(F) = \Phi(G)$,
si bien qu'il existe une unique application
$\bar{\Phi}_{\mm,\nn} : \ptron(\hmodqd)(\mm,\nn) \longrightarrow C(\Phi(\mm^{p}),\Phi(\nn^{p}))$ tel que
$\bar{\Phi}_{\mm,\nn} \circ \ptron(\Gamma)_{\mm^{p},\nn^{p}} = \Phi_{\mm^{p},\nn^{p}}$.

On définit un foncteur $\widetilde{\Phi} : \ptron(\hmodqd) \longrightarrow C$ en posant
$\widetilde{\Phi} = \Phi : Ob(\ptron(\hmodqd)) \longrightarrow Ob(C)$,
et pour tout $F \in \ptron(\hmodqd)(\mm,\nn)$~:
$\widetilde{\Phi}(F) = \Phi(R_{\nn}) \circ \bar{\Phi}(F) \circ \Phi(R_{\mm})^{\sharp}$.

Pour $F \in \modqd(\mm,\nn)$, les morphismes $R_{\nn} \circ \tilde{R}_{\nn*}^{\flat}(F \circ R_{\mm})$ et
$F \circ R_{\mm}$ sont Quillen-homotopes, donc
$\Phi(R_{\nn}) \circ \Phi(\tilde{R}_{\nn*}^{\flat}(F \circ R_{\mm})) = \Phi(F) \circ \Phi(R_{\mm})$,
soit $\Phi(F) =  \Phi(R_{\nn}) \circ \bar{\Phi}(\ptron(\Gamma)(F)) \circ \Phi(R_{\mm})^{\sharp}$.
Comme $\ptron(\Gamma)$ est l'identité sur les objets, on a bien $\widetilde{\Phi} \circ \ptron(\Gamma) = \Phi$.

Enfin, soit un foncteur $\Psi : \ptron(\hmodqd) \longrightarrow C$ vérifiant
$\Psi \circ \ptron(\Gamma) = \Phi$. D'une part $\widetilde{\Phi}$ et $\Psi$ coïncident sur les
objets du fait que $\ptron(\Gamma)$ est l'identité sur les objets. D'autre part, pour $F \in \ptron(\hmodqd)(\mm,\nn)$,
il existe $\bar{F} \in \modqd_{\circ}(\mm^{p},\nn^{p})$ tel que
$F = \ptron(\Gamma)(R_{\nn} \circ \bar{F}) \circ \ptron(\Gamma)(R_{\mm})^{\sharp}$, ce qui implique
que $\Psi(F) = \widetilde{\Phi}(F)$. On a ainsi montré $\Psi = \widetilde{\Phi}$.
\med

\begin{corol}
On suppose que la sous-2-catégorie $\modqd \subset \modqc$ est stable par équivalence.
Le foncteur composé $\modqd_{\circ} \longrightarrow \ptron(\modqd) \longrightarrow \ptron(\hmodqd)$ est une
localisation de $\modqd_{\circ}$ relativement aux équivalences de Quillen.
\end{corol}

\dem
Comme les équivalences de $\modqd$ sont des équivalences de Quillen, il suffit, compte tenu de la proposition
\ref{p-ptr-modq}, de vérifier que $\modqd_{\circ} \longrightarrow \ptron(\modqd)$ est la localisation de
$\modqd_{\circ}$ relativement aux équivalences. Par la proposition \ref{p-ptr-loceq}, il suffit de vérifier que
tout objet de $\modqd_{\circ}$ possède une cotensorisation par $I = (0 \rightleftarrows 1)$.
Soient $\mm, \nn \in \modqd$. Soit $\mm^{I}$ la catégorie des foncteurs de $I$ dans $\mm$. Via l'équivalence
$e_{0} : \mm^{I} \rightarrow \mm$, on constate que $\mm^{I}$ est munie d'une structure de catégorie modèle dont les
fibrations, cofibrations et équivalences faibles sont définies argument par argument.
Soient $F, G \in \modqd_{\circ}(\nn,\mm)$. Il est clair que la donnée d'un 2-isomorphisme $\tau : F \rightarrow G$
de $\modqd$ est équivalente à la donnée d'un foncteur $H_{\tau} : \nn \rightarrow \mm^{I}$ tel que
$e_{0} \circ H_{\tau} = F$ et $e_{1} \circ H_{\tau} = G$.
Or, un tel foncteur détermine toujours un morphisme de $\modqd_{\circ}$. En effet, d'une part il
possède un adjoint à droite qui associe à $f : M_{0} \stackrel{\sim}{\rightarrow} M_{1} \in \mm^{I}$ le produit fibré
de $F^{\flat}(f) : F^{\flat}(M_{0}) \longrightarrow F^{\flat}(M_{1})$ et de
$\bar{\tau} : G^{\flat}(M_{1}) \longrightarrow F^{\flat}(M_{1})$ (où $\bar{\tau}$ désigne le conjugué de $\tau$),
et d'autre part l'identité $e_{0} \circ H_{\tau} = F$ montre qu'il est alors un foncteur de Quillen à gauche.
On a ainsi un isomorphisme $\Cat(I,\modqd(\nn,\mm)) \simeq \modqd(\nn,\mm^{I})$ qui montre que le modèle $\mm^{I}$
est la cotensorisation de $\mm$ par $I$ dans $\modqd$.
\med

\subsection{Variantes}

Dans cette section, on indique brièvement comment les constructions des sections qui précèdent peuvent être
adaptées, d'une part aux modèles pointés et aux modèles stables, et d'autre part aux modèles simpliciaux et
aux modèles spectraux.

\subsubsection{Modèles pointés et modèles stables}

On note $\modqpt$ et $\modqst$ les sous-2-catégories de $\modq$ constituées des modèles de Quillen pointés et
des modèles de Quillen stables respectivement.

L'inclusion $\modqpt \hookrightarrow \modq$ possède un 2-adjoint à gauche $\Phi_{*} : \modq \rightarrow \modqpt$
tel que, pour tout $\mm \in \modq$, $\Phi_{*}(\mm) = \ast/\mm$, la catégorie des objets en-dessous de l'objet
terminal $\ast$ de $\mm$. L'unité de la 2-adjonction $\Phi_{*} : \modq \rightleftarrows \modqpt$ est l'adjonction
de Quillen $\mm \rightleftarrows \Phi_{*}(\mm)$, dont l'adjoint à gauche est défini, pour tout $M \in \mm$, par
$M \mapsto M_{+} = M \amalg \ast$.

\medskip

Pour toute petite catégorie $C$, on note $\uu_{\ast}(C)$ la catégorie $\simp_{\ast}^{C^{op}}$ des préfaisceaux en
ensembles simplicaux pointés sur $C$, que l'on munie de la structure de modèle projective.
Pour $C$ une petite catégorie et $S$ un ensemble de morphisme de $\uu(C)$, on a un isomorphisme de modèle
$\Phi_{*}(\uu (C)/S) \simeq \uu_{\ast}(C)/S_{+}$.

\begin{defi} \label{def-pres-pt}
Une catégorie modèle de la forme $\uu_{\ast}(C)/S_{+}$, où $C$ est une petite catégorie et $S$ un ensemble de
morphisme de $\uu(C)$, est dite \emph{présentable pointée}.
Une \emph{petite présentation pointée} d'une catégorie modèle pointée $\mm$ est la donnée d'une catégorie modèle
présentable pointée $\uu_{\ast}(C)/S_{+}$ et d'une équivalence de Quillen $\uu_{\ast}(C)/S_{+} \longrightarrow \mm$.
\end{defi}

Par la proposition \cite[Prop. 4.7(a)]{D3}, on a~:

\begin{prop} \label{th-res-pres-pt}
Tout modèle pointé combinatoire possède une petite présentation pointée. \carre
\end{prop}

Pour tout $\mm \in \modqc$ et $\nn \in \modqcpt$, l'isomorphisme $\modqpt(\Phi_{*}(\mm),\nn) \simeq \modq(\mm,\nn)$
préserve les homotopies de Quillen, si bien que l'on a un isomorphisme des localisations relativement aux homotopies
de Quillen $\modqpt(\Phi_{*}(\mm),\nn)[\qq^{-1}] \simeq \modq(\mm,\nn)[\qq^{-1}]$.

\begin{prop} \label{p-inv-htp-pt}
Soit une équivalence de Quillen $G : \mm \rightarrow \nn$ dans $\modqcpt$. Pour tout modèle présentable pointé $\pp$
de $\modqcpt$, le foncteur $\tilde{G}_{*} : \modqpt(\pp,\mm)[\qq^{-1}] \longrightarrow \modqpt(\pp,\nn)[\qq^{-1}]$
est une équivalence de catégorie.
\end{prop}

\dem
On a $\pp = \uu_{\ast}(C)/S_{+} \simeq \Phi_{*}(\uu (C)/S)$ pour une petite catégorie $C$ et un ensemble $S$ de
morphisme de $\uu(C)$, donc la proposition découle de la proposition \ref{p-inv-htp} par l'isomorphisme induit
par la 2-adjonction $\Phi_{*} : \modq \rightleftarrows \modqpt$ décrit ci-dessus.
\med

Soit maintenant $\modqdpt$ une sous-2-catégorie pleine de $\modqcpt$ satisfaisant les conditions suivantes~:
\begin{itemize}
    \item tout modèle $\mm \in \modqdpt$ possède une petite présentation pointée dans $\modqdpt$,
    \item pour tout modèle présentable $\mm \in \modqdpt$, le modèle cylindre $\cyl(\mm)$ appartient à $\modqdpt$.
\end{itemize}
Compte tenu des énoncés ci-dessus et de la remarque \ref{r-ch-gen}, il est possible de reproduire pour $\modqdpt$
les résultats de la section \ref{sect-ps-loc-modqc}, en particulier la construction d'une pseudo-localisation
$\Gamma : \modqdpt \longrightarrow \hmodqdpt$ de $\modqdpt$ relativement aux équivalences de Quillen.

Les conditions sont vérifiées par la 2-catégorie $\modqcpt$, ainsi que par sa sous-2-catégorie pleine $\modqcst$.

\subsubsection{Modèles simpliciaux et modèles spectraux}

Soient $\vv$ une catégorie modèle monoïdale combinatoire \cite{Ho1}, et $\vv$-$\modqc$ la 2-catégorie des
$\vv$-catégories modèles combinatoires.

Pour toute petite catégorie $C$, on note $\uu_{\vv}(C)$ la catégorie $\vv^{C^{op}}$ des foncteurs contravariants
de $C$ dans $\vv$, que l'on munie de la structure de modèle projective.

\begin{defi} \label{def-pres-enr}
Une $\vv$-catégorie modèle de la forme $\uu_{\vv}(C)/S$, où $C$ est une petite catégorie et $S$ un ensemble de
morphisme $\uu_{\vv}(C)$, est dite \emph{$\vv$-présentable}.
Une \emph{petite $\vv$-présentation} d'une $\vv$-catégorie modèle $\mm$ est la donnée d'une $\vv$-catégorie modèle
$\vv$-présentable $\uu_{\vv}(C)/S$ et d'une équivalence de Quillen $\uu_{\vv}(C)/S \longrightarrow \mm$ dans
$\vv$-$\modqc$.
\end{defi}

Soient $\mm$ une $\vv$-catégorie modèle et $C$ une petite catégorie. On a une équivalence de catégories
$\coc_{\vv}(\uu_{\vv}(C),\mm) \rightleftarrows \mm^{C}$, entre la catégorie des $\vv$-foncteurs cocontinus
de $\uu(C)$ dans $\mm$ et la catégorie des foncteurs de $C$ dans $\mm$, qui se restreint en une équivalence de
catégories $\mbox{$\vv$-$\modq$}(\uu_{\vv}(C),\mm) \rightleftarrows \mm_{c}^{C}$ dont le but est la
sous-catégorie pleine de $\mm^{C}$ constituée des cofibrants injectifs, et qui fait correspondre aux
homotopies de Quillen les équivalences faibles argument par argument.

Une démonstration similaire (mais plus simple) à celle de la proposition \ref{p-plgt-eq-loc} prouve alors
la proposition suivante.

\begin{prop} \label{p-plgt-eq-loc-enr}
Soient $\mm$ un $\vv$-modèle combinatoire, $C$ une petite catégorie et $S \subset \uu_{\vv}(C)$
un ensemble de morphisme.
Les foncteurs $T^{*} :\mbox{$\vv$-$\modq$}(\uu_{\vv}(C)/S,\mm) \hookrightarrow \mbox{$\vv$-$\modq$}(\uu_{\vv}(C),\mm)$
et $F \mapsto F^{C}(h_{\vv}^{C})$~: $\mbox{$\vv$-$\modq$}(\uu(C),\mm) \rightarrow \mm^{C}$ induisent un plongement
et une équivalence
$$\xymatrix{ \mbox{$\vv$-$\modq$}(\uu_{\vv}(C)/S,\mm)[\qq^{-1}] \ar@{^{(}->}[r] &
\mbox{$\vv$-$\modq$}(\uu_{\vv}(C),\mm)[\qq^{-1}] \ar[r]^(.65){\sim} & \Ho(\mm^{C}) }$$
au niveau des catégories localisées relativement aux homotopies de Quillen et aux équivalences faibles
projectives respectivement.   \carre
\end{prop}

De cette dernière proposition, on déduit (cf. proposition \ref{p-inv-htp})~:

\begin{prop} \label{p-inv-htp-enr}
Soit une équivalence de Quillen $G : \mm \rightarrow \nn$ dans $\vv$-$\modqc$. Pour tout $\vv$-modèle $\vv$-présentable
$\pp$, le foncteur
$\tilde{G}_{*} : \mbox{$\vv$-$\modq$}(\pp,\mm)[\qq^{-1}] \longrightarrow \mbox{$\vv$-$\modq$}(\pp,\nn)[\qq^{-1}]$
est une équivalence de catégorie.  \carre
\end{prop}

Soit alors $\vv$-$\modqd$ une sous-2-catégorie pleine de $\vv$-$\modqc$ satisfaisant les conditions suivantes~:
\begin{itemize}
    \item tout modèle $\mm \in \mbox{$\vv$-$\modqd$}$ possède une petite $\vv$-présentation dans $\vv$-$\modqd$,
    \item pour tout modèle $\vv$-présentable $\mm \in \mbox{$\vv$-$\modqd$}$, le modèle cylindre $\cyl(\mm)$ appartient
    à $\vv$-$\modqd$.
\end{itemize}
Compte tenu des énoncés ci-dessus et de la remarque \ref{r-ch-gen}, il est possible de reproduire pour $\vv$-$\modqd$
les résultats de la section \ref{sect-ps-loc-modqc}, en particulier la construction d'une pseudo-localisation
$\Gamma : \mbox{$\vv$-$\modqd$} \longrightarrow \mbox{$\vv$-$\modqd$}$ de $\vv$-$\modqd$ relativement aux équivalences
de Quillen.

\medskip

On se restreint maintenant aux cas où $\vv$ est la catégorie modèle $\simp$ des ensembles simpliciaux ou la catégorie
modèle $\Sps = \Sps(\simp;S^{1})$ des spectres symétriques.

\begin{prop} \label{th-res-pres-simp-spec}
Tout modèle simplicial combinatoire possède une petite $\simp$-présentation dans $\simp$-$\modqc$.
Tout modèle spectral combinatoire possède une petite $\Sps$-présentation dans $\Sps$-$\modqc$.
\end{prop}

\dem
La première assertion résulte de la proposition \cite[Prop. 2.3]{D1} (cf. également la proposition \cite[Prop. 5.4]{D3}).
La deuxième assertion est prouvée dans la démonstration de \cite[Prop. 6.4]{D3}).
\med

Les conditions ci-dessus sont donc vérifiées par les 2-catégories $\simp$-$\modqc$ et $\Sps$-$\modqc$, ainsi que par la
sous-2-catégorie pleine $\simp$-$\modqcst$ de $\simp$-$\modqc$ constituée des modèles simpliciaux stables combinatoires.

\begin{prop} \label{p-comp-htp}
Les pseudo-foncteurs $\mbox{$\simp$-$\hmodqc$} \longrightarrow \hmodqc$ et
$\mbox{$\Sps$-$\hmodqc$} \longrightarrow \mbox{$\simp$-$\hmodqcst$}$, induits par les 2-foncteurs de changement de base,
sont des biéquivalences.
\end{prop}

De plus, les constructions de \cite{Ho2} permettent d'obtenir un pseudo-foncteur
$\widetilde{\Sps} : \mbox{$\simp$-$\hmodqc$} \longrightarrow \mbox{$\simp$-$\hmodqcst$}$ bi-adjoint à gauche
de l'inclusion $\mbox{$\simp$-$\hmodqcst$} \hookrightarrow \mbox{$\simp$-$\hmodqc$}$.

\bigskip

\dem
Du théorème \cite[Theorem 1.1]{D2} (cf. théorème \ref{th-res-pres}) et de \cite[\S 6.1]{D3}, il découle que les
pseudo-foncteurs $\mbox{$\simp$-$\hmodqc$} \longrightarrow \hmodqc$ et
$\mbox{$\Sps$-$\hmodqc$} \longrightarrow \mbox{$\simp$-$\hmodqcst$}$ sont 2-essentiellement surjectifs.
Pour montrer que le premier est également une équivalence locale, on est amené à vérifier que, pour tout
$\simp$-présentable $\mm \in \mbox{$\simp$-$\modqc$}$, pour toute petite catégorie $C$ et tout ensemble de morphisme
$S \subset \uu_{\simp}(C) = \uu(C)$, le foncteur
$\mbox{$\simp$-$\modq$}(\uu(C)/S,\mm)[\qq^{-1}] \longrightarrow \modq(\uu(C)/S,\mm)[\qq^{-1}]$ est une équivalence
de catégorie.
Des propositions \ref{p-plgt-eq-loc} et \ref{p-plgt-eq-loc-enr}, on tire le diagramme commutatif suivant~:
$$\xymatrix{
\mbox{$\simp$-$\modq$}(\uu(C)/S,\mm)[\qq^{-1}] \ar[d]_{} \ar@{^{(}->}[r] &
\mbox{$\simp$-$\modq$}(\uu(C),\mm)[\qq^{-1}] \ar[d]_{} \ar[r]^(.67){\sim} & \Ho(\mm^{C}) \ar@{=}[d]_{} \\
\modq(\uu(C)/S,\mm)[\qq^{-1}] \ar@{^{(}->}[r] & \modq(\uu(C),\mm)[\qq^{-1}] \ar[r]^(.65){\sim} & \Ho(\mm^{C})
}$$
dont l'équivalence de catégorie requise se déduit de la même façon que dans la démonstration de la proposition
\ref{p-inv-htp}. On vérifie de même que le pseudo-foncteur
$\mbox{$\Sps$-$\hmodqc$} \longrightarrow \mbox{$\simp$-$\hmodqcst$}$ est une équivalence locale.
\med

\section{Dérivateurs}

\subsection{Prédérivateurs, dérivateurs}

On rappelle dans cette section des définitions et résultats concernant la notion de dérivateur,
détaillés dans \cite{M} et \cite{C1,C2}.

\begin{defi} \label{def-pder}
Un \emph{prédérivateur} (de domaine $\cat$) est un 2-foncteur $\cat^{\circ} \longrightarrow \Cat$,
où $\cat^{\circ}$ est la 2-catégorie déduite de $\cat$ en inversant le sens des 1-morphismes
et des 2-morphismes.

On note $\pder$ la 2-catégorie dont les objets sont les prédérivateurs, les 1-morphismes les
transformations pseudo-naturelles, et les 2-morphismes les modifications.
\end{defi}

Ainsi, un prédérivateur $\D$ associe à une petite catégorie $A$ une catégorie $\D(A)$,
à un foncteur $u : A \rightarrow B$ un foncteur $u^{*} : \D(B) \rightarrow \D(A)$,
et à une transformation naturelle $\alpha : u \Rightarrow v$ une transformation naturelle
$\alpha^{*} : v^{*} \Rightarrow u^{*}$.

Tout prédérivateur $\D$ détermine un prédérivateur \emph{opposé} $\D^{\circ}$ tel que
$\D^{\circ}(A) = \D(A^{op})^{op}$ pour tout $A \in \cat$.

\medskip

On utilise dans la suite la construction que voici. Soit un 2-morphisme de $\Cat$
$$\xymatrix{
C \ar[d]_{f} \ar[r]^{k} & C' \ar[d]^{f'} \\
D \ar[r]_{l} \ar@{=>}[ur]_{\theta} & D' }$$
tel que $f$ et $f'$ possèdent chacun un adjoint à gauche, notés $g$ et $g'$ respectivement, avec
$$\eta : Id \longrightarrow f \circ g  \qquad  \varepsilon : g \circ f \longrightarrow Id
\qquad \eta' : Id \longrightarrow f' \circ g'  \qquad  \varepsilon' : g' \circ f' \longrightarrow Id$$
pour morphismes d'adjonction. On tire de cette situation un nouveau 2-morphisme
$$\xymatrix{
C \ar[r]^{k} & C' \\
D \ar[u]^{g} \ar[r]_{l} & D' \ar[u]_{g'} \ar@{=>}[ul]_{\bar{\theta}}  }$$
où $\bar{\theta} = (g' \circ l \circ \eta) \circ (g' \circ \theta \circ g) \circ (\varepsilon' \circ k \circ g)~:
g' \circ l \longrightarrow g' \circ l \circ f \circ g \longrightarrow
g' \circ f' \circ k \circ g \longrightarrow k \circ g$.

\bigskip

On commence par utiliser cette construction dans la situation suivante.
Pour tout $A \in \cat$ et $a \in A$, on note $i_{A,a} : e \rightarrow A$ le foncteur de source la catégorie
ponctuelle et de valeur $a$.
Soient $u : A \rightarrow B$ dans $\cat$ et $b \in B$. On note $b/A$ la catégorie dont les objets sont les
couples $(a,f)$ tels que $a \in A$ et $f : b \rightarrow u(a) \in B$, et $j_{u,b} : b/A \rightarrow A$ le
foncteur envoyant le couple $(a,f)$ sur $a$. Les morphismes structuraux des objets de $b/A$ définissent
un 2-morphisme $\alpha_{u,b} : u \circ j_{u,b} \rightarrow i_{B,b} \circ p_{b/A}$
$$\xymatrix{
b/A \ar[d]_{p_{b/A}} \ar[r]^{j_{u,b}} & A \ar[d]^{u}  \\
e \ar[r]_{i_{B,b}} \ar@{=>}[ur]_(.4){\alpha_{u,b}} & B }$$
Soit $\D$ un prédérivateur.
Si les foncteurs $u^{*}$ et $(p_{b/A})^{*}$ possèdent chacun un adjoint à gauche, noté respectivement $u_{!}$ et
$(p_{b/A})_{!}$, la construction ci-dessus appliquée au 2-morphisme $(\alpha_{u,b})^{*}$ définit le 2-morphisme
de \emph{changement de base} $c_{u,b} : (p_{b/A})_{!} \circ j_{u,b}^{*} \longrightarrow i_{B,b}^{*} \circ u_{!}$
$$\xymatrix{
\D(b/A) \ar[d]_{(p_{b/A})_{!}} \ar@{=>}[drr]_{c_{u,b}} & & \D(A) \ar[ll]_{j_{u,b}^{*}} \ar[d]^{u_{!}}  \\
\D(e) & & \D(B) \ar[ll]^{i_{B,b}^{*}}  }$$

\begin{defi} \label{def-der-dr}
Un \emph{dérivateur faible à droite} est un prédérivateur $\D$ satisfaisant aux axiomes suivants~:

\noindent\textbf{Der~1} (a) $\D(\emptyset) = e$, la catégorie finale;

\noindent (b) Pour tout $A,B \in \cat$, le foncteur induit par les foncteurs canoniques $A \rightarrow A \amalg B$
et $B \rightarrow A \amalg B$~:
$$\D(A \amalg B) \longrightarrow \D(A) \times \D(B)$$
est une équivalence de catégories.

\noindent\textbf{Der~2} La famille des foncteurs $i_{A,a}^{*} : \D(A) \longrightarrow \D(e)$, $a \in A$, est
conservative (i.e. reflète les isomorphismes).

\noindent\textbf{Der~3d} Pour tout $u : A \rightarrow B$ dans $\cat$, le foncteur $u^{*} : \D(B) \rightarrow \D(A)$
possède un adjoint à gauche $u_{!} : \D(A) \rightarrow \D(B)$.

\noindent\textbf{Der~4d} Pour tout $u : A \rightarrow B$ dans $\cat$ et $b \in B$, le morphisme de changement
de base $c_{u,b} : (p_{b/A})_{!} \circ j_{u,b}^{*} \longrightarrow i_{B,b}^{*} \circ u_{!}$ est un isomorphisme.
\end{defi}

Soient $F : \D \rightarrow \D'$ un morphisme de $\pder$ entre dérivateurs faibles à droite et
$u : A \rightarrow B$ dans $\cat$. La construction développée plus haut appliquée au diagramme
commutatif
$$\xymatrix{ \D(A) \ar[r]^{F} & \D'(A)  \\   \D(B) \ar[r]_{F} \ar[u]^{u^{*}} & \D'(B) \ar[u]_{u^{*}}  }$$
définit un 2-morphisme $u_{!} \circ F \rightarrow F \circ u_{!}$.

\begin{defi} \label{def-der-coc}
Un morphisme $F : \D \rightarrow \D'$ entre dérivateur faible à droite est \emph{cocontinu} si, pour tout
$u : A \rightarrow B$ dans $\cat$, le 2-morphisme $u_{!} \circ F \rightarrow F \circ u_{!}$ ci-dessus
est un isomorphisme.
On note $\derd$ la sous-2-catégorie de $\pder$ constituée des dérivateurs faible à droite
et des morphismes cocontinus.
\end{defi}

On définit dualement la notion de dérivateur faible à gauche, vérifiant, outre les axiomes \textbf{Der~1} et
\textbf{Der~2}, les axiomes \textbf{Der~3g} et \textbf{Der~4g} duaux des axiomes \textbf{Der~3d} et \textbf{Der~4d}.
Un prédérivateur $\D$ est un dérivateur faible à gauche si et seulement si son opposé $\D^{\circ}$ est un dérivateur
faible à droite.

\begin{defi} \label{def-der}
Un \emph{dérivateur} est un prédérivateur qui est un dérivateur faible à droite et à gauche.
On note $\der$ la sous-2-catégorie pleine de $\pder$ constituée des dérivateurs.
\end{defi}

On désigne par $\derad$ la 2-catégorie dont les objets sont les dérivateurs et dont les 1-morphismes
sont les adjonctions entre ceux-ci. On constate facilement qu'un adjoint à gauche est cocontinu,
si bien que l'on a un 2-foncteur ``oubli de l'adjoint à droite" $\derad \rightarrow \derd$.

\subsection{Dérivateur d'un modèle de Quillen}

On décrit dans cette section un pseudo-foncteur de la 2-catégorie des théories homotopiques de Quillen
combinatoires dans celle des dérivateurs. On rappelle pour cela la construction de \cite{C1} d'un pseudo-foncteur
$\Phi : \modq \longrightarrow \derad$ de la 2-catégorie des modèles de Quillen dans celle des dérivateurs.

\medskip

Pour toute catégorie modèle $\mm$, on définit un 2-foncteur $\Phi(\mm) : \cat^{\circ} \longrightarrow \Cat$ de la
façon suivante. Pour toute petite catégorie $C$, on pose $\Phi(\mm)(C) = \Ho(\mm^{C^{op}})$, la catégorie localisée
de $\mm^{C^{op}}$ relativement aux équivalences faibles argument par argument. Pour tout foncteur
$u : A \rightarrow B \in \cat$, le foncteur $\Phi(\mm)(u) : \Ho(\mm^{B^{op}}) \longrightarrow \Ho(\mm^{A^{op}})$
est induit par le foncteur $(u^{op})^{*} : \mm^{B^{op}} \longrightarrow \mm^{A^{op}}$, qui préserve les équivalences
faibles. Une transformation naturelle $\alpha : u \Rightarrow v$ induit alors une transformation naturelle
$\Phi(\mm)(\alpha) : \Phi(\mm)(v) \Rightarrow \Phi(\mm)(u)$ de manière évidente.

Il est démontré dans \cite{C1} que le prédérivateur ainsi défini est un dérivateur.
A titre indicatif, lorsque la catégorie modèle $\mm$ est combinatoire, la catégorie $\mm^{C^{op}}$ peut être munie
de la structure de modèle projective, ce qui assure l'existence de $\Ho(\mm^{C^{op}})$ et que le foncteur
$(u^{op})^{*} : \mm^{B^{op}} \longrightarrow \mm^{A^{op}}$ est un foncteur de Quillen à droite. On en déduit
facilement les axiomes \textbf{Der~1}, \textbf{Der~2} et \textbf{Der~3d}.

Pour tout morphisme $F : \mm \rightarrow \nn$ de $\modq$ et toute petite catégorie $C$, on a une adjonction
de foncteurs dérivés totaux
$\Lb(F^{C^{op}}) : \Ho(\mm^{C^{op}}) \rightleftarrows \Ho(\nn^{C^{op}}) : \Rb(F^{\flat C^{op}})$
qui détermine un morphisme $\Phi(F) : \Phi(\mm) \longrightarrow \Phi(\nn)$ de $\derad$.
Tout 2-morphisme $\tau : F \rightarrow G$ de $\modq$ détermine également ainsi un 2-morphisme
$\Phi(\tau) : \Phi(F) \longrightarrow \Phi(G)$ de $\derad$.
Cette procédure définit un pseudo-foncteur $\Phi : \modq \longrightarrow \derad$ qui envoie les équivalences
de Quillen dans les équivalences.

\medskip

Du théorème \ref{th-pu-hmodqd}, on déduit l'existence d'un pseudo-foncteur
$\widetilde{\Phi} : \hmodqc \longrightarrow \derad$ et d'un isomorphisme pseudo-naturel
$\alpha : \widetilde{\Phi} \circ \Gamma \stackrel{\sim}{\longrightarrow} \Phi$,
déterminés à isomorphisme pseudo-naturel unique près. Par la proposition \ref{p-pu-exist}, on peut choisir
$\widetilde{\Phi}$ et $\alpha$ tels que $\widetilde{\Phi}(\mm) = \Phi(\mm)$ et $\; \alpha_{\mm} = Id_{\Phi(\mm)}$.

\medskip

Soient $C$ une petite catégorie et $S$ un ensemble de morphisme de $\uu(C)$. On note $\hot_{C}$ pour
$\Phi(\uu(C))$ et $\hot_{C,S}$ pour $\Phi(\uu(C)/S)$. Le morphisme $T : \uu(C) \longrightarrow \uu(C)/S$
de $\modqc$ induit un morphisme $\Phi(T) : \hot_{C} \longrightarrow \hot_{C,S}$ de $\derad$.
On remarque pour la suite que la co-unité de cette dernière adjonction de $\der$ est un isomorphisme,
et que la classe image inverse des isomorphismes par le foncteur $\Phi(T)(e) : \hot_{C}(e) \longrightarrow \hot_{C,S}(e)$
n'est autre que la classe des $S$-équivalences dans $\hot_{C}(e)$.

\subsection{Une équivalence locale}

Soient $\D$ un dérivateur faible à droite et $C \in \cat$.
On note à nouveau $h^{C}$ l'objet de $\hot_{C}(C^{op}) = \Ho(\uu(C)^{C})$ associé au foncteur de Yoneda.
Un morphisme $F : \hot_{C} \rightarrow \D$ détermine un foncteur
$F(C^{op}) : \hot_{C}(C^{op}) \longrightarrow \D(C^{op})$. On en tire un foncteur
$$\derd(\hot_{C},\D) \longrightarrow \D(C^{op}) \quad \quad F \mapsto F(C^{op})(h^{C})$$
On a le théorème de représentation suivant, qui est une spécialisation de \cite[Corollaire 3.26]{C2}.

\begin{theorem}\cite{C2} \label{th-der-rep}
Pour tout $\D \in \der$ et tout $C \in \cat$, le foncteur $\derd(\hot_{C},\D) \longrightarrow \D(C^{op})$
est une équivalence de catégorie.   \carre
\end{theorem}

Ce théorème est le principal outil dans la vérification du résultat suivant.

\begin{theorem}\label{th-eq-loc}
Le pseudo-foncteur $\widetilde{\Phi} : \hmodqc \longrightarrow \derad$ est une équivalence locale.
\end{theorem}

\dem
On doit vérifier que, pour tout $\nn, \mm \in \hmodqc$, le foncteur
$$\widetilde{\Phi}_{\nn,\mm} : \hmodqc(\nn,\mm) \longrightarrow \derad(\widetilde{\Phi}(\nn),\widetilde{\Phi}(\mm))$$
est une équivalence de catégorie. Comme l'équivalence de Quillen $R_{\nn} : \nn^{p} \rightarrow \nn$ induit une
équivalence dans $\hmodqc$, via $\Gamma$, et dans $\derad$, via $\widetilde{\Phi}$, on peut supposer que $\nn$ est
présentable. La même remarque s'applique à $\mm$.
Il s'agit donc de montrer que pour tout $C \in \cat$, tout ensemble de morphisme $S \subset \uu(C)$ et tout
$\mm \in \hmodqc$ présentable, le foncteur
$$\widetilde{\Phi}_{\uu(C)/S,\mm} : \hmodqc(\uu(C)/S,\mm) \longrightarrow \derad(\hot_{C,S},\widetilde{\Phi}(\mm))$$
est une équivalence de catégorie.

\smallskip

On considère d'abord le cas particulier $S = \varnothing$. Dans le diagramme suivant, le carré de gauche est
commutatif et le carré de droite est un 2-isomorphisme induit par une résolution cofibrante de $h^{C}$ dans
$\uu(C)^{C}$~:
$$\xymatrix{
\modqc(\uu(C),\mm) \ar[d]_{} \ar[r]^{\loc} & \hmodqc(\uu(C),\mm) \ar[d]^{\wr} \ar[r]^{\widetilde{\Phi}} &
\derd(\hot_{C},\widetilde{\Phi}(\mm)) \ar[d]^{\wr} \ar@{=>}[dl]_{\sim}   \\
\mm^{C} \ar[r]_(.45){\loc} & \Ho\left(\mm^{C}\right) \ar@{=}[r] & \widetilde{\Phi}(\mm)(C^{op})  }$$
L'équivalence centrale est celle de la proposition \ref{p-plgt-eq-loc} et l'équivalence de droite
celle du théorème \ref{th-der-rep} ci-dessus. Il s'en suit que la composée
$$\xymatrix{
\hmodqc(\uu(C),\mm) \ar[r]^{\widetilde{\Phi}} & \derad(\hot_{C},\widetilde{\Phi}(\mm)) \ar@{^{(}->}[r] &
\derd(\hot_{C},\widetilde{\Phi}(\mm))      }$$
dont le second foncteur est le plongement ``oubli de l'adjoint à droite", est une équivalence, si bien que
$\widetilde{\Phi}_{\uu(C),\mm} : \hmodqc(\uu(C),\mm) \longrightarrow \derad(\hot_{C},\widetilde{\Phi}(\mm))$
est également une équivalence.

\smallskip

On revient maintenant au cas général.
Du morphisme $T : \uu(C) \longrightarrow \uu(C)/S \in \modqc$, on tire le diagramme suivant~:
$$\xymatrix{
\modqc(\uu(C)/S,\mm) \ar[d]_{T^{*}} \ar[r]^{\loc} & \hmodqc(\uu(C)/S,\mm) \ar[d]_{\Gamma(T)^{*}} \ar[r]^{\widetilde{\Phi}} &
\derad(\hot_{C,S},\widetilde{\Phi}(\mm)) \ar[d]^{\Phi(T)^{*}} \ar@{^{(}->}[r] \ar@{=>}[dl]_{\sim} &
\pder(\hot_{C,S},\widetilde{\Phi}(\mm)) \ar[d]^{\Phi(T)^{*}}  \\
\modqc(\uu(C),\mm)  \ar[r]_{\loc} & \hmodqc(\uu(C),\mm)  \ar[r]_{\widetilde{\Phi}} &
\derad(\hot_{C},\widetilde{\Phi}(\mm)) \ar@{^{(}->}[r] & \pder(\hot_{C},\widetilde{\Phi}(\mm))    }$$
Le carré central est composé des 2-isomorphismes $\Phi(T)^{*} \simeq \widetilde{\Phi}(\Gamma(T))^{*}$ et
$\com^{\widetilde{\Phi}}_{T,-} : \widetilde{\Phi}(\Gamma(T))^{*} \circ  \widetilde{\Phi}
\stackrel{\sim}{\longrightarrow} \widetilde{\Phi} \circ \Gamma(T)^{*}$.
Les carrés latéraux sont commutatifs, les plongements de droite étant les foncteurs ``oubli de l'adjoint à droite".

Comme la co-unité $\varepsilon : \Phi(T) \circ \Phi(T)^{\flat} \longrightarrow Id$ de l'adjonction
$\Phi(T) : \hot_{C} \rightleftarrows \hot_{C,S} : \Phi(T)^{\flat}$ de $\pder$ est un isomorphisme,
l'unité $\widetilde{\Phi}(\varepsilon)^{*}$ de l'adjonction induite $(\Phi(T)^{\flat*},\Phi(T)^{*})$
est un isomorphisme, ce qui fait du foncteur
$\Phi(T)^{*} : \pder(\hot_{C,S},\widetilde{\Phi}(\mm)) \longrightarrow \pder(\hot_{C},\widetilde{\Phi}(\mm))$,
et donc du foncteur
$\Phi(T)^{*} : \derad(\hot_{C,S},\widetilde{\Phi}(\mm)) \longrightarrow \derad(\hot_{C},\widetilde{\Phi}(\mm))$,
un plongement.
Par ailleurs, la proposition \ref{p-plgt-eq-loc} indique que le foncteur
$\Gamma(T)^{*} = \widetilde{T^{*}} : \hmodqc(\uu(C)/S,\mm) \longrightarrow \hmodqc(\uu(C),\mm)$ est un plongement.

Le foncteur $\Phi(T)^{*} \circ \widetilde{\Phi}$, isomorphe au plongement $\widetilde{\Phi} \circ \Gamma(T)^{*}$,
est donc lui-même un plongement, et comme $\Phi(T)^{*}$ est également un plongement, on a montré que
$\widetilde{\Phi}_{\uu(C)/S,\mm}$ est pleinement fidèle.

\smallskip

Il reste à vérifier que $\widetilde{\Phi}_{\uu(C)/S,\mm}$ est essentiellement surjectif.
Soit $G \in \derad(\hot_{C,S},\widetilde{\Phi}(\mm))$. Comme $\widetilde{\Phi}_{\uu(C),\mm}$ est une équivalence,
il existe $F \in \modqc(\uu(C),\mm)$ et un 2-isomorphisme $\Phi(F) \simeq \Phi(T)^{*}(G)$. Le 2-isomorphisme
$$\Lb(F) = \Phi(F)(e) \simeq G(e) \circ \Phi(T)(e) : \Ho(\uu(C)) \longrightarrow \Ho(\uu(C)/S) \longrightarrow \Ho(\mm)$$
montre que $\Lb(F)$ envoie l'ensemble $S$ dans les isomorphismes, ce qui assure l'existence d'un morphisme
$\bar{F} \in \modqc(\uu(C)/S,\mm)$ tel que $T^{*}(\bar{F}) = F$ et donc
$\Phi(T)^{*}(\widetilde{\Phi}(\loc\bar{F})) \simeq \widetilde{\Phi}(F) \simeq \Phi(T)^{*}(G)$.
Comme $\Phi(T)^{*}$ est un plongement, on a obtenu un morphisme $\loc\bar{F} \in \hmodqc(\uu(C)/S,\mm)$ et un
2-isomorphisme $\widetilde{\Phi}(\loc\bar{F}) \simeq G$, démontrant ainsi la surjectivité essentielle de
$\widetilde{\Phi}$.
\med

\subsection{Dérivateurs de petite présentation}

Dans cette section, on introduit les dérivateurs de petite présentation dans le but de déduire du
théorème \ref{th-eq-loc} une biéquivalence.

\begin{defi} \label{def-loc-der}
Soit $\D$ un prédérivateur. Une \emph{localisation} de $\D$ est la donnée d'une adjonction
$F : \D \rightleftarrows \D' : F^{\flat}$ de $\pder$ dont la co-unité
$\varepsilon^{F} : F \circ F^{\flat} \rightarrow Id$ est un isomorphisme.
\end{defi}

Il s'avère que la localisation préserve un certain nombre de propriété.

\begin{lem} \label{l-loc-der} \cite[Lemme 4.2]{C2}
Une localisation d'un dérivateur est également un dérivateur. \carre
\end{lem}

\begin{defi} \label{def-der-ppres}
Soit $\D$ un prédérivateur.
Une \emph{petite génération} de $\D$ est la donnée d'une petite catégorie $C \in \cat$ et d'une localisation
$\hot_{C} \rightleftarrows \D$.
Une \emph{petite présentation} de $\D$ est la donnée d'une petite génération $\hot_{C} \rightleftarrows \D$
et d'un ensemble $S$ de morphisme de $\hot_{C}(e)$ tels que les $S$-équivalences coïncident avec l'image
inverse des isomorphismes par le foncteur $\hot_{C}(e) \rightarrow \D(e)$.
Le prédérivateur $\D$ est dit de petite présentation s'il possède une petite présentation.
\end{defi}

Il découle du lemme \ref{l-loc-der} qu'un prédérivateur de petite présentation est nécessairement un dérivateur.
L'exemple fondamental de petite présentation est la localisation
$\Phi(T) : \hot_{C} \rightleftarrows \hot_{C,S} : \Phi(T)^{\flat}$ associé à une petite catégorie $C$
et à un ensemble $S$ de morphisme de $\hot_{C}(e)$.

On note $\deradpp$ la sous-2-catégorie pleine de $\derad$ constituée des dérivateurs de petite présentation.
Le théorème \ref{th-res-pres} montre que le pseudo-foncteur $\Phi : \modqc \longrightarrow \derad$
est en fait à valeur dans $\deradpp$.

\begin{theorem}\label{th-bieq}
Le pseudo-foncteur $\widetilde{\Phi} : \hmodqc \longrightarrow \deradpp$ est une biéquivalence.
\end{theorem}

\dem
Vu le théorème \ref{th-eq-loc}, il reste à vérifier que le pseudo-foncteur $\widetilde{\Phi}$ est 2-essentiellement
surjectif.
Soit $\D \in \deradpp$. Soient $C \in \cat$, un ensemble $S$ de morphisme de $\hot_{C}(e)$ et une localisation
$F : \hot_{C} \rightleftarrows \D$ constituant une petite présentation de $\D$. De la localisation
$\Phi(T) : \hot_{C} \rightleftarrows \hot_{C,S} : \Phi(T)^{\flat}$, on tire deux morphismes
$\widetilde{F} : \hot_{C,S} \rightleftarrows \D : \widetilde{F}^{\flat}$ de $\der$ définis par
$\widetilde{F} = F \circ \Phi(T)^{\flat}$ et $\widetilde{F}^{\flat} = \Phi(T) \circ F^{\flat}$.
Le 2-morphisme composé
$$\xymatrix{
Id \ar[rr]_{(\varepsilon^{F})^{-1}}^{\sim} & & FF^{\flat} \ar[rr]_(.35){F\eta^{\Phi(T)}F^{\flat}} & &
F\Phi(T)^{\flat}\Phi(T)F^{\flat} = \widetilde{F}\widetilde{F}^{\flat}      }$$
est un 2-isomorphisme. En effet, par l'axiome \textbf{Der~2}, il suffit de vérifier que $F\eta^{\Phi(T)}$
est un isomorphisme dans $\D(e)$, ce qui découle de ce que $\eta^{\Phi(T)}$ est une $S$-équivalence, car
$\Phi(T)(\eta^{\Phi(T)})$ est un isomorphisme, via une identité triangulaire. De même, le 2-morphisme composé
$$\xymatrix{
Id \ar[rr]_{(\varepsilon^{\Phi(T)})^{-1}}^{\sim} & & \Phi(T)\Phi(T)^{\flat} \ar[rr]_(.4){\Phi(T)\eta^{F}\Phi(T)^{\flat}}
& & \Phi(T)F^{\flat}F\Phi(T)^{\flat} = \widetilde{F}^{\flat}\widetilde{F}      }$$
est un 2-isomorphisme, et on a ainsi obtenu une équivalence entre $\D$ et $\hot_{C,S} = \widetilde{\Phi}(\uu(C)/S)$.
\med

\newpage

\appendix

\section{Quelques diagrammes}

Dans cet appendice sont indiqués quelques diagrammes commutatifs auxquels font référence certaines démonstrations
des pages qui précèdent.

\subsection{Le pseudo-foncteur $\Gamma$}
\label{app-gam}

Axiomes d'unité
$$\xymatrix{
Id \, \tilde{R}_{2*}^{\flat}(F R_{1}) \ar[rr]^{\uni^{2}_{Id}  \, \tilde{R}_{2*}^{\flat}(F R_{1})}
\ar[dr]_(.4){\uni^{2}_{\left[\tilde{R}_{2*}^{\flat}(F R_{1})\right]}} \ar@{=}[ddd]  & &
\tilde{R}_{2*}^{\flat}(Id \, R_{2}) \tilde{R}_{2*}^{\flat}(F R_{1})
\ar[d]^{\uni^{2}_{\left[\tilde{R}_{2*}^{\flat}(Id \, R_{2}) \tilde{R}_{2*}^{\flat}(F R_{1})\right]}}
\\
 & (\tilde{R}_{2*}^{\flat} \tilde{R}_{2*})\left(\tilde{R}_{2*}^{\flat}(FR_{1})\right)
\ar[r] \ar@{=}[dr] &
(\tilde{R}_{2*}^{\flat} \tilde{R}_{2*})\left(\tilde{R}_{2*}^{\flat}(Id \, R_{2}) \tilde{R}_{2*}^{\flat}(FR_{1})\right)
\ar[d]^{\tilde{R}_{2*}^{\flat}\left(\cou^{2}_{(Id \, R_{2})} \tilde{R}_{2*}^{\flat}(F R_{1})\right)}
\\
 & &
\tilde{R}_{2*}^{\flat}\left(Id \, (\tilde{R}_{2*} \tilde{R}_{2*}^{\flat})(F R_{1})\right)
\ar[d]^{\tilde{R}_{2*}^{\flat}\left(Id \cou^{2}_{(F R_{1})}\right)}
\\
\tilde{R}_{2*}^{\flat}(F R_{1}) \ar@{=}[rr] & & \tilde{R}_{2*}^{\flat}(Id \, F R_{1})    }$$

\bigskip

$$\xymatrix{
\tilde{R}_{2*}^{\flat}(F R_{1}) Id \ar[rr]^{\tilde{R}_{2*}^{\flat}(F R_{1}) \, \uni^{1}_{Id}}
\ar[dr]_(.4){\uni^{2}_{\left[\tilde{R}_{2*}^{\flat}(F R_{1})\right]}} \ar@{=}[ddd]  & &
\tilde{R}_{2*}^{\flat}(F R_{1}) \tilde{R}_{1*}^{\flat}(Id \, R_{1})
\ar[d]^{\uni^{2}_{\left[\tilde{R}_{2*}^{\flat}(F R_{1}) \tilde{R}_{1*}^{\flat}(Id \, R_{1})\right]}}
\\
 & (\tilde{R}_{2*}^{\flat} \tilde{R}_{2*})\left(\tilde{R}_{2*}^{\flat}(FR_{1})\right)  \ar[r]  \ar@{=}[dr] &
(\tilde{R}_{2*}^{\flat} \tilde{R}_{2*})\left(\tilde{R}_{2*}^{\flat}(F R_{1}) \tilde{R}_{1*}^{\flat}(Id \, R_{1})\right)
\ar[d]^{\tilde{R}_{2*}^{\flat}\left(\cou^{2}_{(F R_{1})} \tilde{R}_{1*}^{\flat}(Id \, R_{1})\right)}
\\
 & & \tilde{R}_{2*}^{\flat}\left(F (\tilde{R}_{2*} \tilde{R}_{2*}^{\flat})(Id \, R_{1})\right)
\ar[d]^{\tilde{R}_{2*}^{\flat}\left(F \cou^{2}_{(Id \, R_{1})}\right)}
\\
\tilde{R}_{2*}^{\flat}(F R_{1}) \ar@{=}[rr] & & \tilde{R}_{2*}^{\flat}(F Id \, R_{1})    }$$

\newpage

\oddsidemargin 0cm
\setlength{\textheight}{22.2cm}
\setlength{\textwidth}{20cm}
\begin{landscape}
$$\xymatrix{
\tilde{R}_{4*}^{\flat}(H R_{3}) \tilde{R}_{3*}^{\flat}(G R_{2}) \tilde{R}_{2*}^{\flat}(F R_{1})
\ar[dddd]_{\uni^{3}_{\left[\tilde{R}_{3*}^{\flat}(G R_{2}) \tilde{R}_{2*}^{\flat}(F R_{1})\right]}}
\ar[ddr]^{\uni^{4}_{}}
\ar[rr]^{\uni^{4}_{\left[\tilde{R}_{4*}^{\flat}(H R_{3}) \tilde{R}_{3*}^{\flat}(G R_{2})\right]}} & &
(\tilde{R}_{4*}^{\flat}\tilde{R}_{4*})
\left[\tilde{R}_{4*}^{\flat}(H R_{3}) \tilde{R}_{3*}^{\flat}(G R_{2})\right] \tilde{R}_{2*}^{\flat}(F R_{1})
\ar[r]^(.53){\cou^{4}_{HR_{3}}} \ar[dd]^{\uni^{4}_{}} &
\tilde{R}_{4*}^{\flat}
\left[H (\tilde{R}_{3*} \tilde{R}_{3*}^{\flat})(G R_{2})\right] \tilde{R}_{2*}^{\flat}(F R_{1})
\ar[r]^(.55){\cou^{3}_{GR_{2}}} \ar[dd]^{\uni^{4}_{}} &
\tilde{R}_{4*}^{\flat}(H G R_{2}) \tilde{R}_{2*}^{\flat}(F R_{1})
\ar[dd]^{\uni^{4}_{}}
 \\
 & & & &  \\
 &
(\tilde{R}_{4*}^{\flat}\tilde{R}_{4*})
\left[\tilde{R}_{4*}^{\flat}(H R_{3}) \tilde{R}_{3*}^{\flat}(G R_{2}) \tilde{R}_{2*}^{\flat}(F R_{1})\right]
\ar[r]^(.45){\uni^{4}_{}}
\ar[dd]_{\uni^{3}_{\left[\tilde{R}_{3*}^{\flat}(G R_{2}) \tilde{R}_{2*}^{\flat}(F R_{1})\right]}}
\ar@{=}[ddr] &
(\tilde{R}_{4*}^{\flat}\tilde{R}_{4*})\left[(\tilde{R}_{4*}^{\flat}\tilde{R}_{4*})
\left[\tilde{R}_{4*}^{\flat}(H R_{3}) \tilde{R}_{3*}^{\flat}(G R_{2})\right] \tilde{R}_{2*}^{\flat}(F R_{1})\right]
\ar[r]^(.53){\cou^{4}_{HR_{3}}} \ar[dd]^{\cou^{4}_{}} &
(\tilde{R}_{4*}^{\flat}\tilde{R}_{4*})\left[\tilde{R}_{4*}^{\flat}
\left[H (\tilde{R}_{3*} \tilde{R}_{3*}^{\flat})(G R_{2})\right] \tilde{R}_{2*}^{\flat}(F R_{1})\right]
\ar[r]^(.55){\cou^{3}_{GR_{2}}} \ar[dd]^{\cou^{4}_{}} &
(\tilde{R}_{4*}^{\flat}\tilde{R}_{4*})
\left[\tilde{R}_{4*}^{\flat}(H G R_{2}) \tilde{R}_{2*}^{\flat}(F R_{1})\right]
\ar[dddd]^{\cou^{4}_{HGR_{2}}}
 \\
 & & & &  \\
\tilde{R}_{4*}^{\flat}(H R_{3}) (\tilde{R}_{3*}^{\flat}\tilde{R}_{3*})
\left[\tilde{R}_{3*}^{\flat}(G R_{2}) \tilde{R}_{2*}^{\flat}(F R_{1})\right]
\ar[r]^(.45){\uni^{4}_{}} \ar[dddd]_{\cou^{3}_{G R_{2}}}  &
(\tilde{R}_{4*}^{\flat}\tilde{R}_{4*})
\left[\tilde{R}_{4*}^{\flat}(H R_{3}) (\tilde{R}_{3*}^{\flat}\tilde{R}_{3*})
\left(\tilde{R}_{3*}^{\flat}(G R_{2}) \tilde{R}_{2*}^{\flat}(F R_{1})\right)\right]
\ar[dddd]^{\cou^{3}_{GR_{2}}} \ar[ddr]_{\cou^{4}_{HR_{3}}} &
\tilde{R}_{4*}^{\flat} \left[(\tilde{R}_{4*}\tilde{R}_{4*}^{\flat})(H R_{3})
\tilde{R}_{3*}^{\flat}(G R_{2}) \tilde{R}_{2*}^{\flat}(F R_{1})\right]
\ar[l]^(.45){\uni^{3}} \ar[r]^(.53){\cou^{4}_{HR_{3}}} &
\tilde{R}_{4*}^{\flat} \left(H (\tilde{R}_{3*}\tilde{R}_{3*}^{\flat})(G R_{2}) \tilde{R}_{2*}^{\flat}(F R_{1})\right)
\ar[ddl]^{\uni^{3}} \ar[ddr]^(.55){\cou^{3}_{GR_{2}}} \ar@{=}[dd] &
 \\
 & & & &  \\
  &  &
\tilde{R}_{4*}^{\flat} \left[H (\tilde{R}_{3*}\tilde{R}_{3*}^{\flat})
\left((\tilde{R}_{3*}\tilde{R}_{3*}^{\flat})(G R_{2}) \tilde{R}_{2*}^{\flat}(F R_{1})\right)\right]
\ar[ddr]_(.4){\cou^{3}_{GR_{2}}}
\ar[r]_(.7){\cou^{3}_{\left[(\tilde{R}_{3*}\tilde{R}_{3*}^{\flat})(G R_{2}) \tilde{R}_{2*}^{\flat}(F R_{1})\right]}} &
\tilde{R}_{4*}^{\flat} \left(H (\tilde{R}_{3*}\tilde{R}_{3*}^{\flat})(G R_{2}) \tilde{R}_{2*}^{\flat}(F R_{1})\right)
\ar[r]_(.55){\cou^{3}_{GR_{2}}} &
\tilde{R}_{4*}^{\flat} \left(H G (\tilde{R}_{2*}\tilde{R}_{2*}^{\flat})(F R_{1})\right)
\ar[dddd]^{\cou^{2}_{FR_{1}}}
 \\
 & & & &  \\
\tilde{R}_{4*}^{\flat}(H R_{3}) \tilde{R}_{3*}^{\flat}\left(G (\tilde{R}_{2*}\tilde{R}_{2*}^{\flat})(F R_{1})\right)
\ar[r]_(.45){\uni^{4}_{}} \ar[dd]_{\cou^{2}_{FR_{1}}}  &
(\tilde{R}_{4*}^{\flat}\tilde{R}_{4*})\left[\tilde{R}_{4*}^{\flat}(H R_{3})
\tilde{R}_{3*}^{\flat}\left(G (\tilde{R}_{2*}\tilde{R}_{2*}^{\flat})(F R_{1})\right)\right]
\ar[rr]_{\cou^{4}_{HR_{3}}} \ar[dd]^{\cou^{2}_{FR_{1}}} & &
\tilde{R}_{4*}^{\flat}\left[H
(\tilde{R}_{3*}\tilde{R}_{3*}^{\flat})\left(G (\tilde{R}_{2*}\tilde{R}_{2*}^{\flat})(F R_{1})\right)\right]
\ar[uur]_{\cou^{3}_{}} \ar[dd]^{\cou^{2}_{FR_{1}}} &
 \\
 & & & &  \\
\tilde{R}_{4*}^{\flat}(H R_{3}) \tilde{R}_{3*}^{\flat}(G F R_{1})
\ar[r]_(.45){\uni^{4}_{\left[\tilde{R}_{4*}^{\flat}(H R_{3}) \tilde{R}_{3*}^{\flat}(G F R_{1})\right]}} &
(\tilde{R}_{4*}^{\flat}\tilde{R}_{4*})
\left[\tilde{R}_{4*}^{\flat}(H R_{3}) \tilde{R}_{3*}^{\flat}(G F R_{1})\right]
\ar[rr]_{\cou^{4}_{HR_{3}}} & &
\tilde{R}_{4*}^{\flat}
\left(H (\tilde{R}_{3*}\tilde{R}_{3*}^{\flat})(G F R_{1})\right)
\ar[r]_{\cou^{3}_{G F R_{1}}} &
\tilde{R}_{4*}^{\flat} (H G F R_{1})
}$$
\end{landscape}
\oddsidemargin .3cm
\setlength{\textheight}{64cm}
\begin{landscape}
$$\xymatrix{
\tilde{R}_{4*}^{\flat}(H R_{3}) \tilde{R}_{3*}^{\flat}(G R_{2}) \tilde{R}_{2*}^{\flat}(F R_{1})
\ar[dddd]_{\uni^{3}_{\left[\tilde{R}_{3*}^{\flat}(G R_{2}) \tilde{R}_{2*}^{\flat}(F R_{1})\right]}}
\ar[ddr]^{\uni^{4}_{}}
\ar[rr]^{\uni^{4}_{\left[\tilde{R}_{4*}^{\flat}(H R_{3}) \tilde{R}_{3*}^{\flat}(G R_{2})\right]}} & &
(\tilde{R}_{4*}^{\flat}\tilde{R}_{4*})
\left[\tilde{R}_{4*}^{\flat}(H R_{3}) \tilde{R}_{3*}^{\flat}(G R_{2})\right] \tilde{R}_{2*}^{\flat}(F R_{1})
\ar[r]^(.53){\cou^{4}_{HR_{3}}} \ar[dd]^{\uni^{4}_{}} &
\tilde{R}_{4*}^{\flat}
\left[H (\tilde{R}_{3*} \tilde{R}_{3*}^{\flat})(G R_{2})\right] \tilde{R}_{2*}^{\flat}(F R_{1})
\ar[r]^(.55){\cou^{3}_{GR_{2}}} \ar[dd]^{\uni^{4}_{}} &
\tilde{R}_{4*}^{\flat}(H G R_{2}) \tilde{R}_{2*}^{\flat}(F R_{1})
\ar[dd]^{\uni^{4}_{}}
 \\
 & & & &  \\
 &
(\tilde{R}_{4*}^{\flat}\tilde{R}_{4*})
\left[\tilde{R}_{4*}^{\flat}(H R_{3}) \tilde{R}_{3*}^{\flat}(G R_{2}) \tilde{R}_{2*}^{\flat}(F R_{1})\right]
\ar[r]^(.45){\uni^{4}_{}}
\ar[dd]_{\uni^{3}_{\left[\tilde{R}_{3*}^{\flat}(G R_{2}) \tilde{R}_{2*}^{\flat}(F R_{1})\right]}}
\ar@{=}[ddr] &
(\tilde{R}_{4*}^{\flat}\tilde{R}_{4*})\left[(\tilde{R}_{4*}^{\flat}\tilde{R}_{4*})
\left[\tilde{R}_{4*}^{\flat}(H R_{3}) \tilde{R}_{3*}^{\flat}(G R_{2})\right] \tilde{R}_{2*}^{\flat}(F R_{1})\right]
\ar[r]^(.53){\cou^{4}_{HR_{3}}} \ar[dd]^{\cou^{4}_{}} &
(\tilde{R}_{4*}^{\flat}\tilde{R}_{4*})\left[\tilde{R}_{4*}^{\flat}
\left[H (\tilde{R}_{3*} \tilde{R}_{3*}^{\flat})(G R_{2})\right] \tilde{R}_{2*}^{\flat}(F R_{1})\right]
\ar[r]^(.55){\cou^{3}_{GR_{2}}} \ar[dd]^{\cou^{4}_{}} &
(\tilde{R}_{4*}^{\flat}\tilde{R}_{4*})
\left[\tilde{R}_{4*}^{\flat}(H G R_{2}) \tilde{R}_{2*}^{\flat}(F R_{1})\right]
\ar[dddd]^{\cou^{4}_{HGR_{2}}}
 \\
 & & & &  \\
\tilde{R}_{4*}^{\flat}(H R_{3}) (\tilde{R}_{3*}^{\flat}\tilde{R}_{3*})
\left[\tilde{R}_{3*}^{\flat}(G R_{2}) \tilde{R}_{2*}^{\flat}(F R_{1})\right]
\ar[r]^(.45){\uni^{4}_{}} \ar[dddd]_{\cou^{3}_{G R_{2}}}  &
(\tilde{R}_{4*}^{\flat}\tilde{R}_{4*})
\left[\tilde{R}_{4*}^{\flat}(H R_{3}) (\tilde{R}_{3*}^{\flat}\tilde{R}_{3*})
\left(\tilde{R}_{3*}^{\flat}(G R_{2}) \tilde{R}_{2*}^{\flat}(F R_{1})\right)\right]
\ar[dddd]^{\cou^{3}_{GR_{2}}} \ar[ddr]_{\cou^{4}_{HR_{3}}} &
\tilde{R}_{4*}^{\flat} \left[(\tilde{R}_{4*}\tilde{R}_{4*}^{\flat})(H R_{3})
\tilde{R}_{3*}^{\flat}(G R_{2}) \tilde{R}_{2*}^{\flat}(F R_{1})\right]
\ar[l]^(.45){\uni^{3}} \ar[r]^(.53){\cou^{4}_{HR_{3}}} &
\tilde{R}_{4*}^{\flat} \left(H (\tilde{R}_{3*}\tilde{R}_{3*}^{\flat})(G R_{2}) \tilde{R}_{2*}^{\flat}(F R_{1})\right)
\ar[ddl]^{\uni^{3}} \ar[ddr]^(.55){\cou^{3}_{GR_{2}}} \ar@{=}[dd] &
 \\
 & & & &  \\
  &  &
\tilde{R}_{4*}^{\flat} \left[H (\tilde{R}_{3*}\tilde{R}_{3*}^{\flat})
\left((\tilde{R}_{3*}\tilde{R}_{3*}^{\flat})(G R_{2}) \tilde{R}_{2*}^{\flat}(F R_{1})\right)\right]
\ar[ddr]_(.45){\cou^{3}_{GR_{2}}}
\ar[r]_(.7){\cou^{3}_{\left[(\tilde{R}_{3*}\tilde{R}_{3*}^{\flat})(G R_{2}) \tilde{R}_{2*}^{\flat}(F R_{1})\right]}} &
\tilde{R}_{4*}^{\flat} \left(H (\tilde{R}_{3*}\tilde{R}_{3*}^{\flat})(G R_{2}) \tilde{R}_{2*}^{\flat}(F R_{1})\right)
\ar[r]_(.55){\cou^{3}_{GR_{2}}} &
\tilde{R}_{4*}^{\flat} \left(H G (\tilde{R}_{2*}\tilde{R}_{2*}^{\flat})(F R_{1})\right)
\ar[dddd]^{\cou^{2}_{FR_{1}}}
 \\
 & & & &  \\
\tilde{R}_{4*}^{\flat}(H R_{3}) \tilde{R}_{3*}^{\flat}\left(G (\tilde{R}_{2*}\tilde{R}_{2*}^{\flat})(F R_{1})\right)
\ar[r]_(.45){\uni^{4}_{}} \ar[dd]_{\cou^{2}_{FR_{1}}}  &
(\tilde{R}_{4*}^{\flat}\tilde{R}_{4*})\left[\tilde{R}_{4*}^{\flat}(H R_{3})
\tilde{R}_{3*}^{\flat}\left(G (\tilde{R}_{2*}\tilde{R}_{2*}^{\flat})(F R_{1})\right)\right]
\ar[rr]_{\cou^{4}_{HR_{3}}} \ar[dd]^{\cou^{2}_{FR_{1}}} & &
\tilde{R}_{4*}^{\flat}\left[H
(\tilde{R}_{3*}\tilde{R}_{3*}^{\flat})\left(G (\tilde{R}_{2*}\tilde{R}_{2*}^{\flat})(F R_{1})\right)\right]
\ar[uur]_{\cou^{3}_{}} \ar[dd]^{\cou^{2}_{FR_{1}}} &
 \\
 & & & &  \\
\tilde{R}_{4*}^{\flat}(H R_{3}) \tilde{R}_{3*}^{\flat}(G F R_{1})
\ar[r]_(.45){\uni^{4}_{\left[\tilde{R}_{4*}^{\flat}(H R_{3}) \tilde{R}_{3*}^{\flat}(G F R_{1})\right]}} &
(\tilde{R}_{4*}^{\flat}\tilde{R}_{4*})
\left[\tilde{R}_{4*}^{\flat}(H R_{3}) \tilde{R}_{3*}^{\flat}(G F R_{1})\right]
\ar[rr]_{\cou^{4}_{HR_{3}}} & &
\tilde{R}_{4*}^{\flat}
\left(H (\tilde{R}_{3*}\tilde{R}_{3*}^{\flat})(G F R_{1})\right)
\ar[r]_{\cou^{3}_{G F R_{1}}} &
\tilde{R}_{4*}^{\flat} (H G F R_{1})
}$$
\end{landscape}
\setlength{\textheight}{23cm}

\subsection{Le pseudo-foncteur $\widetilde{\Phi}$}
\label{app-phi}

Axiomes d'unité

\bigskip

$$\xymatrix{
\Phi(R_{2}) \bar{\Phi}(F) \Phi(R_{1})^{\sharp} \ar[r]^(.4){\uni^{\Phi(R_{2})}} \ar@{=}[ddr] &
\Phi(R_{2}) \Phi(R_{2})^{\sharp} \Phi(R_{2}) \bar{\Phi}(F) \Phi(R_{1})^{\sharp} \ar[r]^(.46){\un^{\Phi}_{}} \ar[dd]^{\cou^{\Phi(R_{2})}} &
\Phi(R_{2}) \bar{\Phi}(Id) \Phi(R_{2})^{\sharp} \Phi(R_{2}) \bar{\Phi}(F) \Phi(R_{1})^{\sharp} \ar[dd]^{\cou^{\Phi(R_{2})}}  \\
 & &  \\
 & \Phi(R_{2}) \bar{\Phi}(F) \Phi(R_{1})^{\sharp} \ar[r]^(.46){\un^{\Phi}_{}} \ar@{=}[ddr] &
\Phi(R_{2}) \bar{\Phi}(Id) \bar{\Phi}(F) \Phi(R_{1})^{\sharp} \ar[dd]^{\com^{\bar{\Phi}}_{Id,F}}  \\
 & &  \\
 & & \Phi(R_{2}) \bar{\Phi}(Id F) \Phi(R_{1})^{\sharp}     }$$

\bigskip

$$\xymatrix{
\Phi(R_{2}) \bar{\Phi}(F) \Phi(R_{1})^{\sharp} \ar[r]^(.4){\uni^{\Phi(R_{1})}} \ar@{=}[ddr] &
\Phi(R_{2}) \bar{\Phi}(F) \Phi(R_{1})^{\sharp} \Phi(R_{1}) \Phi(R_{1})^{\sharp} \ar[r]^(.46){\un^{\Phi}_{}} \ar[dd]^{\cou^{\Phi(R_{1})}} &
\Phi(R_{2}) \bar{\Phi}(F) \Phi(R_{1})^{\sharp} \Phi(R_{1}) \bar{\Phi}(Id) \Phi(R_{1})^{\sharp} \ar[dd]^{\cou^{\Phi(R_{1})}}  \\
 & &  \\
 & \Phi(R_{2}) \bar{\Phi}(F) \Phi(R_{1})^{\sharp} \ar[r]^(.46){\un^{\Phi}_{}} \ar@{=}[ddr] &
\Phi(R_{2}) \bar{\Phi}(F) \bar{\Phi}(Id) \Phi(R_{1})^{\sharp} \ar[dd]^{\com^{\bar{\Phi}}_{F,Id}}  \\
 & &  \\
 & & \Phi(R_{2}) \bar{\Phi}(F Id) \Phi(R_{1})^{\sharp}     }$$

Axiome d'associativité

{\tiny
$$\xymatrix{
\Phi(R_{4}) \bar{\Phi}(H) \Phi(R_{3})^{\sharp} \Phi(R_{3}) \bar{\Phi}(G) \Phi(R_{2})^{\sharp} \Phi(R_{2}) \bar{\Phi}(F) \Phi(R_{1})^{\sharp}
\ar[r]^(.57){\cou^{\Phi(R_{2})}}  \ar[dd]_{\cou^{\Phi(R_{3})}} &
\Phi(R_{4}) \bar{\Phi}(H) \Phi(R_{3})^{\sharp} \Phi(R_{3}) \bar{\Phi}(G) \bar{\Phi}(F) \Phi(R_{1})^{\sharp}
\ar[r]^(.52){\com^{\bar{\Phi}}_{G,F}}  \ar[dd]^{\cou^{\Phi(R_{3})}} &
\Phi(R_{4}) \bar{\Phi}(H) \Phi(R_{3})^{\sharp} \Phi(R_{3}) \bar{\Phi}(G F) \Phi(R_{1})^{\sharp}  \ar[dd]^{\cou^{\Phi(R_{3})}}  \\
 & &  \\
\Phi(R_{4}) \bar{\Phi}(H) \bar{\Phi}(G) \Phi(R_{2})^{\sharp} \Phi(R_{2}) \bar{\Phi}(F) \Phi(R_{1})^{\sharp}
\ar[r]_(.57){\cou^{\Phi(R_{2})}}  \ar[dd]_{\com^{\bar{\Phi}}_{H,G}} &
\Phi(R_{4}) \bar{\Phi}(H) \bar{\Phi}(G) \bar{\Phi}(F) \Phi(R_{1})^{\sharp}
\ar[r]_(.52){\com^{\bar{\Phi}}_{G,F}}  \ar[dd]^{\com^{\bar{\Phi}}_{H,G}} &
\Phi(R_{4}) \bar{\Phi}(H) \bar{\Phi}(G F) \Phi(R_{1})^{\sharp}  \ar[dd]^{\com^{\bar{\Phi}}_{H,GF}}  \\
 & &  \\
\Phi(R_{4}) \bar{\Phi}(H G) \Phi(R_{2})^{\sharp} \Phi(R_{2}) \bar{\Phi}(F) \Phi(R_{1})^{\sharp}
\ar[r]_(.57){\cou^{\Phi(R_{2})}}  &
\Phi(R_{4}) \bar{\Phi}(H G) \bar{\Phi}(F) \Phi(R_{1})^{\sharp}
\ar[r]_(.52){\com^{\bar{\Phi}}_{HG,F}}  &
\Phi(R_{4}) \bar{\Phi}(H G F) \Phi(R_{1})^{\sharp}    }$$
}

\newpage

\subsection{La transformation pseudo-naturelle $\alpha$}
\label{app-alph}

Le 2-isomorphisme
$\un^{\widetilde{\Phi} \circ \Gamma}_{\mm_{1}} : Id_{\left(\widetilde{\Phi} \circ \Gamma\right)(\mm_{1})}
\stackrel{\sim}{\longrightarrow} \left(\widetilde{\Phi} \circ \Gamma\right)(Id_{\mm_{1}})$
est la composée~:
$$\xymatrix{
Id_{(\widetilde{\Phi} \circ \Gamma)(\mm_{1})} \ar[rr]^{\sim}_{\un^{\widetilde{\Phi}}_{\Gamma(\mm_{1})}} & &
\widetilde{\Phi}(Id_{\Gamma(\mm_{1})}) \ar[rr]^{\sim}_{\widetilde{\Phi}(\un^{\Gamma}_{\mm_{1}})} & &
(\widetilde{\Phi} \circ \Gamma)(Id_{\mm_{1}})    }$$
soit, en détail~:
$$\xymatrix{
Id_{\Phi(\mm_{1})}
\ar[d]^{\uni^{\Phi(R_{1})}}  \\
\Phi(R_{1}) \circ \Phi(R_{1})^{\sharp}
\ar[d]^{\Phi(R_{1}) \circ \un^{\Phi}_{\mm_{1}} \circ \Phi(R_{1})^{\sharp}}  \\
\Phi(R_{1}) \circ \Phi(Id_{\mm_{1}}) \circ \Phi(R_{1})^{\sharp}
\ar[d]^{\Phi(R_{1}) \circ \bar{\Phi}\left(\uni^{1}_{Id_{\mm_{1}}}\right) \circ \Phi(R_{1})^{\sharp}}  \\
\Phi(R_{1}) \circ \bar{\Phi}
\left((\tilde{R}_{1*}^{\flat} \circ \tilde{R}_{1*})(Id_{\mm_{1}})\right) \circ \Phi(R_{1})^{\sharp}    }$$

\bigskip

Le 2-isomorphisme
$\com^{\widetilde{\Phi} \circ \Gamma}_{G,F} :
\left(\widetilde{\Phi} \circ \Gamma\right)(G) \circ \left(\widetilde{\Phi} \circ \Gamma\right)(F)
\stackrel{\sim}{\longrightarrow} \left(\widetilde{\Phi} \circ \Gamma\right)(G \circ F)$
est la composée~:
$$\xymatrix{
\left(\widetilde{\Phi} \circ \Gamma\right)(G) \circ \left(\widetilde{\Phi} \circ \Gamma\right)(F)
\ar[rr]^(.6){\sim}_(.6){\com^{\widetilde{\Phi}}_{\Gamma(G),\Gamma(F)}} & &
\widetilde{\Phi} (\Gamma(G) \circ \Gamma(F))
\ar[rr]^{\sim}_{\widetilde{\Phi}(\com^{\Gamma}_{G,F})} & &
\left(\widetilde{\Phi} \circ \Gamma\right)(G \circ F)     }$$
soit, en détail~:
$$\xymatrix{
\Phi(R_{3}) \circ \bar{\Phi}\left(\tilde{R}_{3*}^{\flat}(G \circ R_{2})\right) \circ \Phi(R_{2})^{\sharp}
\circ \Phi(R_{2}) \circ \bar{\Phi}\left(\tilde{R}_{2*}^{\flat}(F \circ R_{1})\right) \circ \Phi(R_{1})^{\sharp}
\ar[d]^{\Phi(R_{3}) \circ \bar{\Phi}\left(\tilde{R}_{3*}^{\flat}(G \circ R_{2})\right) \circ \cou^{\Phi(R_{2})}
\circ \bar{\Phi}\left(\tilde{R}_{2*}^{\flat}(F \circ R_{1})\right) \circ \Phi(R_{1})^{\sharp}}
\\
\Phi(R_{3}) \circ \bar{\Phi}\left(\tilde{R}_{3*}^{\flat}(G \circ R_{2})\right)
\circ \bar{\Phi}\left(\tilde{R}_{2*}^{\flat}(F \circ R_{1})\right) \circ \Phi(R_{1})^{\sharp}
\ar[d]^{\Phi(R_{3}) \circ
\com^{\bar{\Phi}}_{\left(\tilde{R}_{3*}^{\flat}(G \circ R_{2}),\tilde{R}_{2*}^{\flat}(F \circ R_{1})\right)}
\circ \Phi(R_{1})^{\sharp}}
\\
\Phi(R_{3}) \circ \bar{\Phi}\left(\tilde{R}_{3*}^{\flat}(G \circ R_{2}) \circ \tilde{R}_{2*}^{\flat}(F \circ R_{1})\right)
\circ \Phi(R_{1})^{\sharp}
\ar[d]^{\Phi(R_{3}) \circ
\bar{\Phi}\left(\uni^{3}_{\left[\tilde{R}_{3*}^{\flat}(G \circ R_{2}) \circ \tilde{R}_{2*}^{\flat}(F \circ R_{1})\right]}\right)
\circ \Phi(R_{1})^{\sharp}}
\\
\Phi(R_{3}) \circ \bar{\Phi}\left((\tilde{R}_{3*}^{\flat} \circ \tilde{R}_{3*})
\left[\tilde{R}_{3*}^{\flat}(G \circ R_{2}) \circ \tilde{R}_{2*}^{\flat}(F \circ R_{1})\right]\right)
\circ \Phi(R_{1})^{\sharp}
\ar[d]^{\Phi(R_{3}) \circ
\bar{\Phi}\left(\tilde{R}_{3*}^{\flat}\left(\cou^{3}_{(G \circ R_{2})} \circ \tilde{R}_{2*}^{\flat}(F \circ R_{1})\right)\right)
\circ \Phi(R_{1})^{\sharp}}
\\
\Phi(R_{3}) \circ \bar{\Phi}\left(
\tilde{R}_{3*}^{\flat}\left[G \circ (\tilde{R}_{2*} \circ \tilde{R}_{2*}^{\flat})(F \circ R_{1})\right]\right)
\circ \Phi(R_{1})^{\sharp}
\ar[d]^{\Phi(R_{3}) \circ
\bar{\Phi}\left(\tilde{R}_{3*}^{\flat}\left(G \circ \cou^{2}_{(F \circ R_{1})}\right)\right)
\circ \Phi(R_{1})^{\sharp}}
\\
\Phi(R_{3}) \circ \bar{\Phi}\left(\tilde{R}_{3*}^{\flat}(G \circ F \circ R_{1})\right) \circ \Phi(R_{1})^{\sharp}  }$$

\oddsidemargin 0cm
\setlength{\textheight}{20.6cm}
\begin{landscape}
$$\xymatrix{
\Phi(R_{3}) \bar{\Phi}\left(\tilde{R}_{3*}^{\flat}(G R_{2})\right) \Phi(R_{2})^{\sharp}
\Phi(R_{2}) \bar{\Phi}\left(\tilde{R}_{2*}^{\flat}(F R_{1})\right) \Phi(R_{1})^{\sharp}
\ar[d]_{\com^{\bar{\Phi}}_{R_{3},\tilde{R}_{3*}^{\flat}(G R_{2})}} \ar[r]^(.57){\cou^{\Phi(R_{2})}} &
\Phi(R_{3}) \bar{\Phi}\left(\tilde{R}_{3*}^{\flat}(G R_{2})\right)
\bar{\Phi}\left(\tilde{R}_{2*}^{\flat}(F R_{1})\right) \Phi(R_{1})^{\sharp}
\ar[d]_{\com^{\bar{\Phi}}} \ar[r]^(.52){\com^{\bar{\Phi}}} &
\Phi(R_{3}) \bar{\Phi}\left(\tilde{R}_{3*}^{\flat}(G R_{2}) (\tilde{R}_{2*}^{\flat}(F R_{1})\right)
\Phi(R_{1})^{\sharp}
\ar[d]_{\com^{\bar{\Phi}}}  \ar[r]^(.45){\bar{\Phi}(\uni^{3})} &
\Phi(R_{3}) \bar{\Phi}\left((\tilde{R}_{3*}^{\flat} \tilde{R}_{3*})
\left[\tilde{R}_{3*}^{\flat}(G R_{2}) (\tilde{R}_{2*}^{\flat}(F R_{1})\right]\right)
\Phi(R_{1})^{\sharp}
\ar[d]_{\com^{\bar{\Phi}}} \ar[r]^(.55){\cou^{3}_{G R_{2}}} &
\Phi(R_{3}) \bar{\Phi}\left(
\tilde{R}_{3*}^{\flat}\left[G (\tilde{R}_{2*} \tilde{R}_{2*}^{\flat})(F R_{1})\right]\right)
\Phi(R_{1})^{\sharp}
\ar[d]_{\com^{\bar{\Phi}}} \ar[r]^(.57){\cou^{2}_{F R_{1}}} &
\Phi(R_{3}) \bar{\Phi}\left(\tilde{R}_{3*}^{\flat}(G F R_{1})\right) \Phi(R_{1})^{\sharp}
\ar[d]^{\com^{\bar{\Phi}}}
\\
\bar{\Phi}\left((\tilde{R}_{3*}\tilde{R}_{3*}^{\flat})(G R_{2})\right) \Phi(R_{2})^{\sharp}
\Phi(R_{2}) \bar{\Phi}\left(\tilde{R}_{2*}^{\flat}(F R_{1})\right) \Phi(R_{1})^{\sharp}
\ar[d]_{\bar{\Phi}(\cou^{3}_{G R_{2}})} \ar[r]^(.57){\cou^{\Phi(R_{2})}} &
\bar{\Phi}\left((\tilde{R}_{3*}\tilde{R}_{3*}^{\flat})(G R_{2})\right)
\bar{\Phi}\left(\tilde{R}_{2*}^{\flat}(F R_{1})\right) \Phi(R_{1})^{\sharp}
\ar[d]_{\bar{\Phi}(\cou^{3}_{G R_{2}})} \ar[r]^(.52){\com^{\bar{\Phi}}} &
\bar{\Phi}\left((\tilde{R}_{3*}\tilde{R}_{3*}^{\flat})(G R_{2})\tilde{R}_{2*}^{\flat}(F R_{1})\right)
\Phi(R_{1})^{\sharp}
\ar@{=}[dr] \ar[r]^(.55){\bar{\Phi}(R_{3} \uni^{3})} &
\bar{\Phi}\left((\tilde{R}_{3*}\tilde{R}_{3*}^{\flat}\tilde{R}_{3*})
\left[\tilde{R}_{3*}^{\flat}(G R_{2}) (\tilde{R}_{2*}^{\flat}(F R_{1})\right]\right)
\Phi(R_{1})^{\sharp}
\ar[d]^{\bar{\Phi}\left(\cou^{3}_{(\tilde{R}_{3*}\tilde{R}_{3*}^{\flat})(G R_{2}) \tilde{R}_{2*}^{\flat}(F R_{1})}\right)}
\ar[r]^(.53){\cou^{3}_{G R_{2}}} &
\bar{\Phi}\left((\tilde{R}_{3*}\tilde{R}_{3*}^{\flat})
\left[G (\tilde{R}_{2*}\tilde{R}_{2*}^{\flat})(F R_{1})\right]\right)
\Phi(R_{1})^{\sharp}
\ar[ddddl]^{\bar{\Phi}\left(\cou^{3}_{G (\tilde{R}_{2*}\tilde{R}_{2*}^{\flat})(F R_{1})}\right)}
\ar[r]^(.57){\cou^{2}_{F R_{1}}} &
\bar{\Phi}\left((\tilde{R}_{3*}\tilde{R}_{3*}^{\flat})(G F R_{1})\right) \Phi(R_{1})^{\sharp}
\ar[ddddd]^{\bar{\Phi}(\cou^{3}_{G F R_{1}})}
\\
\Phi(G R_{2}) \Phi(R_{2})^{\sharp} \Phi(R_{2}) \bar{\Phi}\left(\tilde{R}_{2*}^{\flat}(F R_{1})\right)
\Phi(R_{1})^{\sharp}
\ar[d]_{(\com^{\Phi}_{G,R_{2}})^{-1}} \ar[r]^(.57){\cou^{\Phi(R_{2})}} &
\Phi(G R_{2}) \bar{\Phi}\left(\tilde{R}_{2*}^{\flat}(F R_{1})\right) \Phi(R_{1})^{\sharp}
\ar[d]_{(\com^{\Phi}_{G,R_{2}})^{-1}} \ar@{=}[dr] &
&
\bar{\Phi}\left((\tilde{R}_{3*}\tilde{R}_{3*}^{\flat})(G R_{2}) \tilde{R}_{2*}^{\flat}(F R_{1})\right)
\Phi(R_{1})^{\sharp}
\ar[ddd]_{\bar{\Phi}\left(\cou^{3}_{G R_{2}}\tilde{R}_{2*}^{\flat}(F R_{1})\right)} & &
\\
\Phi(G) \Phi(R_{2}) \Phi(R_{2})^{\sharp} \Phi(R_{2}) \bar{\Phi}\left(\tilde{R}_{2*}^{\flat}(F R_{1})\right)
\Phi(R_{1})^{\sharp}
\ar[d]_{(\uni^{\Phi(R_{2})})^{-1}} \ar[r]^(.57){\cou^{\Phi(R_{2})}} &
\Phi(G) \Phi(R_{2}) \bar{\Phi}\left(\tilde{R}_{2*}^{\flat}(F R_{1})\right)
\Phi(R_{1})^{\sharp}
\ar[ddl]^{\com^{\Phi}_{R_{2},\tilde{R}_{2*}^{\flat}(F R_{1})}} \ar[r]_{\com^{\Phi}_{G,R_{2}}}  &
\Phi(G R_{2}) \bar{\Phi}\left(\tilde{R}_{2*}^{\flat}(F R_{1})\right)
\Phi(R_{1})^{\sharp}
\ar[ddr]^{\com^{\Phi}_{G R_{2},\tilde{R}_{2*}^{\flat}(F R_{1})}} &
& &
\\
\Phi(G) \Phi(R_{2}) \bar{\Phi}\left(\tilde{R}_{2*}^{\flat}(F R_{1})\right) \Phi(R_{1})^{\sharp}
\ar[d]_{\com^{\Phi}_{R_{2},\tilde{R}_{2*}^{\flat}(F R_{1})}} \ar@{=}[ur] &
& & & &
\\
\Phi(G) \bar{\Phi}\left((\tilde{R}_{2*}\tilde{R}_{2*}^{\flat})(F R_{1})\right) \Phi(R_{1})^{\sharp}
\ar[d]_{\bar{\Phi}(\cou^{2}_{F R_{1}})}  & & &
\bar{\Phi}\left(G (\tilde{R}_{2*}\tilde{R}_{2*}^{\flat})(F R_{1})\right) \Phi(R_{1})^{\sharp}
\ar[lll]^{(\com^{\Phi}_{G,(\tilde{R}_{2*}\tilde{R}_{2*}^{\flat})(F R_{1})})^{-1}}
\ar[drr]^{\bar{\Phi}(G \cou^{2}_{F R_{1}})} &
 &
\\
\Phi(G) \Phi(F R_{1}) \Phi(R_{1})^{\sharp}
\ar[d]_{(\com^{\Phi}_{F,R_{1}})^{-1}} \ar[rrrrr]_(.45){\com^{\Phi}_{GF,R_{1}}} &
& & & &
\Phi(G F R_{1}) \Phi(R_{1})^{\sharp}
\ar[d]^{(\com^{\Phi}_{GF,R_{1}})^{-1}}
\\
\Phi(G) \Phi(F) \Phi(R_{1}) \Phi(R_{1})^{\sharp}
\ar[d]_{(\uni^{\Phi(R_{1})})^{-1}} \ar[rrrrr]_(.45){\com^{\Phi}_{G,F}} &
& & & &
\Phi(G F) \Phi(R_{1}) \Phi(R_{1})^{\sharp}
\ar[d]^{(\uni^{\Phi(R_{1})})^{-1}}
\\
\Phi(G) \Phi(F)
\ar[rrrrr]_(.45){\com^{\Phi}_{G,F}} &
& & & &
\Phi(G F)    }$$
\end{landscape}
\oddsidemargin 0.3cm
\setlength{\textheight}{54.5cm}
\begin{landscape}
$$\xymatrix{
\Phi(R_{3}) \bar{\Phi}\left(\tilde{R}_{3*}^{\flat}(G R_{2})\right) \Phi(R_{2})^{\sharp}
\Phi(R_{2}) \bar{\Phi}\left(\tilde{R}_{2*}^{\flat}(F R_{1})\right) \Phi(R_{1})^{\sharp}
\ar[d]_{\com^{\bar{\Phi}}_{R_{3},\tilde{R}_{3*}^{\flat}(G R_{2})}} \ar[r]^(.57){\cou^{\Phi(R_{2})}} &
\Phi(R_{3}) \bar{\Phi}\left(\tilde{R}_{3*}^{\flat}(G R_{2})\right)
\bar{\Phi}\left(\tilde{R}_{2*}^{\flat}(F R_{1})\right) \Phi(R_{1})^{\sharp}
\ar[d]_{\com^{\bar{\Phi}}} \ar[r]^(.52){\com^{\bar{\Phi}}} &
\Phi(R_{3}) \bar{\Phi}\left(\tilde{R}_{3*}^{\flat}(G R_{2}) (\tilde{R}_{2*}^{\flat}(F R_{1})\right)
\Phi(R_{1})^{\sharp}
\ar[d]_{\com^{\bar{\Phi}}}  \ar[r]^(.45){\bar{\Phi}(\uni^{3})} &
\Phi(R_{3}) \bar{\Phi}\left((\tilde{R}_{3*}^{\flat} \tilde{R}_{3*})
\left[\tilde{R}_{3*}^{\flat}(G R_{2}) (\tilde{R}_{2*}^{\flat}(F R_{1})\right]\right)
\Phi(R_{1})^{\sharp}
\ar[d]_{\com^{\bar{\Phi}}} \ar[r]^(.55){\cou^{3}_{G R_{2}}} &
\Phi(R_{3}) \bar{\Phi}\left(
\tilde{R}_{3*}^{\flat}\left[G (\tilde{R}_{2*} \tilde{R}_{2*}^{\flat})(F R_{1})\right]\right)
\Phi(R_{1})^{\sharp}
\ar[d]_{\com^{\bar{\Phi}}} \ar[r]^(.57){\cou^{2}_{F R_{1}}} &
\Phi(R_{3}) \bar{\Phi}\left(\tilde{R}_{3*}^{\flat}(G F R_{1})\right) \Phi(R_{1})^{\sharp}
\ar[d]^{\com^{\bar{\Phi}}}
\\
\bar{\Phi}\left((\tilde{R}_{3*}\tilde{R}_{3*}^{\flat})(G R_{2})\right) \Phi(R_{2})^{\sharp}
\Phi(R_{2}) \bar{\Phi}\left(\tilde{R}_{2*}^{\flat}(F R_{1})\right) \Phi(R_{1})^{\sharp}
\ar[d]_{\bar{\Phi}(\cou^{3}_{G R_{2}})} \ar[r]^(.57){\cou^{\Phi(R_{2})}} &
\bar{\Phi}\left((\tilde{R}_{3*}\tilde{R}_{3*}^{\flat})(G R_{2})\right)
\bar{\Phi}\left(\tilde{R}_{2*}^{\flat}(F R_{1})\right) \Phi(R_{1})^{\sharp}
\ar[d]_{\bar{\Phi}(\cou^{3}_{G R_{2}})} \ar[r]^(.52){\com^{\bar{\Phi}}} &
\bar{\Phi}\left((\tilde{R}_{3*}\tilde{R}_{3*}^{\flat})(G R_{2})\tilde{R}_{2*}^{\flat}(F R_{1})\right)
\Phi(R_{1})^{\sharp}
\ar@{=}[dr] \ar[r]^(.45){\bar{\Phi}(R_{3} \uni^{3})} &
\bar{\Phi}\left((\tilde{R}_{3*}\tilde{R}_{3*}^{\flat}\tilde{R}_{3*})
\left[\tilde{R}_{3*}^{\flat}(G R_{2}) (\tilde{R}_{2*}^{\flat}(F R_{1})\right]\right)
\Phi(R_{1})^{\sharp}
\ar[d]^{\bar{\Phi}\left(\cou^{3}_{(\tilde{R}_{3*}\tilde{R}_{3*}^{\flat})(G R_{2}) \tilde{R}_{2*}^{\flat}(F R_{1})}\right)}
\ar[r]^(.53){\cou^{3}_{G R_{2}}} &
\bar{\Phi}\left((\tilde{R}_{3*}\tilde{R}_{3*}^{\flat})
\left[G (\tilde{R}_{2*}\tilde{R}_{2*}^{\flat})(F R_{1})\right]\right)
\Phi(R_{1})^{\sharp}
\ar[ddddl]^{\bar{\Phi}\left(\cou^{3}_{G (\tilde{R}_{2*}\tilde{R}_{2*}^{\flat})(F R_{1})}\right)}
\ar[r]^(.57){\cou^{2}_{F R_{1}}} &
\bar{\Phi}\left((\tilde{R}_{3*}\tilde{R}_{3*}^{\flat})(G F R_{1})\right) \Phi(R_{1})^{\sharp}
\ar[ddddd]^{\bar{\Phi}(\cou^{3}_{G F R_{1}})}
\\
\Phi(G R_{2}) \Phi(R_{2})^{\sharp} \Phi(R_{2}) \bar{\Phi}\left(\tilde{R}_{2*}^{\flat}(F R_{1})\right)
\Phi(R_{1})^{\sharp}
\ar[d]_{(\com^{\Phi}_{G,R_{2}})^{-1}} \ar[r]^(.57){\cou^{\Phi(R_{2})}} &
\Phi(G R_{2}) \bar{\Phi}\left(\tilde{R}_{2*}^{\flat}(F R_{1})\right) \Phi(R_{1})^{\sharp}
\ar[d]_{(\com^{\Phi}_{G,R_{2}})^{-1}} \ar@{=}[dr] &
&
\bar{\Phi}\left((\tilde{R}_{3*}\tilde{R}_{3*}^{\flat})(G R_{2}) \tilde{R}_{2*}^{\flat}(F R_{1})\right)
\Phi(R_{1})^{\sharp}
\ar[ddd]_{\bar{\Phi}\left(\cou^{3}_{G R_{2}}\tilde{R}_{2*}^{\flat}(F R_{1})\right)} & &
\\
\Phi(G) \Phi(R_{2}) \Phi(R_{2})^{\sharp} \Phi(R_{2}) \bar{\Phi}\left(\tilde{R}_{2*}^{\flat}(F R_{1})\right)
\Phi(R_{1})^{\sharp}
\ar[d]_{(\uni^{\Phi(R_{2})})^{-1}} \ar[r]^(.57){\cou^{\Phi(R_{2})}} &
\Phi(G) \Phi(R_{2}) \bar{\Phi}\left(\tilde{R}_{2*}^{\flat}(F R_{1})\right)
\Phi(R_{1})^{\sharp}
\ar[ddl]^{\com^{\Phi}_{R_{2},\tilde{R}_{2*}^{\flat}(F R_{1})}} \ar[r]_{\com^{\Phi}_{G,R_{2}}}  &
\Phi(G R_{2}) \bar{\Phi}\left(\tilde{R}_{2*}^{\flat}(F R_{1})\right)
\Phi(R_{1})^{\sharp}
\ar[ddr]^{\com^{\Phi}_{G R_{2},\tilde{R}_{2*}^{\flat}(F R_{1})}} &
& &
\\
\Phi(G) \Phi(R_{2}) \bar{\Phi}\left(\tilde{R}_{2*}^{\flat}(F R_{1})\right) \Phi(R_{1})^{\sharp}
\ar[d]_{\com^{\Phi}_{R_{2},\tilde{R}_{2*}^{\flat}(F R_{1})}} \ar@{=}[ur] &
& & & &
\\
\Phi(G) \bar{\Phi}\left((\tilde{R}_{2*}\tilde{R}_{2*}^{\flat})(F R_{1})\right) \Phi(R_{1})^{\sharp}
\ar[d]_{\bar{\Phi}(\cou^{2}_{F R_{1}})}  & & &
\bar{\Phi}\left(G (\tilde{R}_{2*}\tilde{R}_{2*}^{\flat})(F R_{1})\right) \Phi(R_{1})^{\sharp}
\ar[lll]^{(\com^{\Phi}_{G,(\tilde{R}_{2*}\tilde{R}_{2*}^{\flat})(F R_{1})})^{-1}}
\ar[drr]^{\bar{\Phi}(G \cou^{2}_{F R_{1}})} &
 &
\\
\Phi(G) \Phi(F R_{1}) \Phi(R_{1})^{\sharp}
\ar[d]_{(\com^{\Phi}_{F,R_{1}})^{-1}} \ar[rrrrr]_(.45){\com^{\Phi}_{GF,R_{1}}} &
& & & &
\Phi(G F R_{1}) \Phi(R_{1})^{\sharp}
\ar[d]^{(\com^{\Phi}_{GF,R_{1}})^{-1}}
\\
\Phi(G) \Phi(F) \Phi(R_{1}) \Phi(R_{1})^{\sharp}
\ar[d]_{(\uni^{\Phi(R_{1})})^{-1}} \ar[rrrrr]_(.45){\com^{\Phi}_{G,F}} &
& & & &
\Phi(G F) \Phi(R_{1}) \Phi(R_{1})^{\sharp}
\ar[d]^{(\uni^{\Phi(R_{1})})^{-1}}
\\
\Phi(G) \Phi(F)
\ar[rrrrr]_(.45){\com^{\Phi}_{G,F}} &
& & & &
\Phi(G F)    }$$
\end{landscape}
\setlength{\textheight}{23cm}
\oddsidemargin 1cm

\newpage

\subsection{La transformation pseudo-naturelle $\chi$}
\label{app-chi}

On commence par vérifier l'axiome de composition de $\chi$ pour deux morphismes composables
$F  \in \hmodqd(\mm_{1},\mm_{2})$ et $G \in \hmodqd(\mm_{2},\mm_{3})$
tels que $\mm_{1} = \mm_{1}^{p}$ et $\mm_{2} = \mm_{2}^{p}$.
Il existe $\bar{F} \in \modqd(\mm_{1}^{p},\mm_{2}^{p})$ et $\bar{G} \in \modqd(\mm_{2}^{p},\mm_{3}^{p})$ tels que
$F = \loc\bar{F} = \Gamma(\bar{F})$ et $G = \loc\bar{G}$.

Par définition de $\tpn^{\chi}$, on a le diagramme commutatif~:
$$\xymatrix{
\chi_{\mm_{3}}  \Psi(G)
\ar[d]_{\tpn^{\chi}_{G}} \ar[rr]^(.45){\chi_{\mm_{3}} \Psi(\uni^{3}_{G})}_(.4){\sim}  & &
\chi_{\mm_{3}} \Psi(\tilde{R}_{3*}^{\flat} \tilde{R}_{3*}(G))
\ar[d]^{\tpn^{\chi}_{\Gamma(R_{3} \bar{G})}} \ar@{=}[r] &
\chi_{\mm_{3}} \Psi(\Gamma(R_{3}\bar{G}))  \ar[d]^{\tpn^{\chi}_{R_{3} \bar{G}}}   \\
\Psi'(G) \chi_{\mm_{2}^{p}} \ar[rr]_(.45){\Psi'(\uni^{3}_{G}) \chi_{\mm_{2}^{p}}}^(.4){\sim}  & &
\Psi'(\tilde{R}_{3*}^{\flat} \tilde{R}_{3*}(G)) \chi_{\mm_{2}^{p}} \ar@{=}[r] &
\Psi'(\Gamma(R_{3}\bar{G})) \chi_{\mm_{2}^{p}}    }$$
Via l'ismorphisme $\uni^{3}_{G} :
G \stackrel{\sim}{\longrightarrow} (\tilde{R}_{3*}^{\flat} \tilde{R}_{3*})(G) = \Gamma(R_{3} \circ \bar{G})$,
le pentagone
$$\xymatrix{
\chi_{\mm_{3}} \Psi(G) \Psi(F)  \ar[d]_{Id \ast \com^{\Psi}_{G,F}} \ar[rr]^{\tpn^{\chi}_{G} \ast Id} & &
\Psi'(G) \chi_{\mm_{2}} \Psi(F)  \ar[rr]^{Id \ast \tpn^{\chi}_{F}} & &
\Psi'(G) \Psi'(F) \chi_{\mm_{1}} \ar[d]^{\com^{\Psi'}_{G,F} \ast Id}   \\
\chi_{\mm_{3}} \Psi(GF) \ar[rrrr]_{\tpn^{\chi}_{GF}}  & & & &  \Psi'(GF) \chi_{\mm_{1}}   }$$
est isomorphe au pentagone supérieur du diagramme
$$\xymatrix{
\chi_{\mm_{3}} \Psi(\Gamma(R_{3}\bar{G})) \Psi(\Gamma(\bar{F}))  \ar[d]_{Id \ast \com^{\Psi}_{\Gamma(R_{3}\bar{G}),\Gamma(\bar{F})}}
\ar[rr]^{\tpn^{\chi}_{\Gamma(R_{3}\bar{G})} \ast Id} & &
\Psi'(\Gamma(R_{3}\bar{G})) \chi_{\mm_{2}} \Psi(\Gamma(\bar{F}))  \ar[rr]^{Id \ast \tpn^{\chi}_{\Gamma(\bar{F})}} & &
\Psi'(\Gamma(R_{3}\bar{G})) \Psi'(\Gamma(\bar{F})) \chi_{\mm_{1}}
\ar[d]^{\com^{\Psi'}_{\Gamma(R_{3}\bar{G}),\Gamma(\bar{F})} \ast Id}   \\
\chi_{\mm_{3}} \Psi(\Gamma(R_{3}\bar{G})\Gamma(\bar{F}))
\ar[d]_{Id \ast \Psi(\com^{\Gamma}_{R_{3}\bar{G},\bar{F}})}
\ar[rrrr]_{\tpn^{\chi}_{\Gamma(R_{3}\bar{G})\Gamma(\bar{F})}}  & & & &
\Psi'(\Gamma(R_{3}\bar{G})\Gamma(\bar{F})) \chi_{\mm_{1}}
\ar[d]^{\Psi'(\com^{\Gamma}_{R_{3}\bar{G},\bar{F}}) \ast Id}    \\
\chi_{\mm_{3}} \Psi(\Gamma(R_{3}\bar{G}\bar{F}))
\ar[rrrr]_{\tpn^{\chi}_{\Gamma(R_{3}\bar{G}\bar{F})}}  & & & &
\Psi'(\Gamma(R_{3}\bar{G}\bar{F})) \chi_{\mm_{1}}       }$$
dont le carré inférieur est commutatif.
Le pentagone externe n'est autre que le pentagone commutatif correspondant à l'axiome de composition de
$\chi$ pour $R_{3} \circ \bar{G}$ et $\bar{F}$. On obtient ainsi l'axiome de composition de $\chi$ pour
deux morphismes composables de sources présentables.

Le cas général s'en déduit par le diagramme des pages suivantes, où l'on utilise également la commutativité
du pentagone suivant~:
$$\xymatrix{
\Psi(G) \Psi(F) \ar[rrrr]^{\com^{\Psi}_{G,F}} \ar[ddd]_{Id \ast \uni^{\Psi(R)} \ast Id} \ar[drr]_{Id \ast \Psi(\uni^{R}F)}
& & & &    \Psi(GF)  \ar[ddd]^{\Psi(G\uni^{R}F)}   \\
 & & \Psi(G) \Psi(RR^{-1}F) \ar[d]_{Id \ast (\com^{\Psi}_{G,R})^{-1}} \ar@/^1pc/[ddrr]^{\com^{\Psi}_{G,RR^{-1}F}} &     \\
 & & \Psi(G) \Psi(RR^{-1}) \Psi(F)  \ar[drr]_(.4){\com^{\Psi}_{G,RR^{-1},F}} & &  &     \\
\Psi(G) \Psi(R) \Psi(R^{-1}) \Psi(F) \ar[urr]^{Id \ast \com^{\Psi}_{R,R} \ast Id} \ar[rr]_(.55){\com^{\Psi}_{G,R} \ast \com^{\Psi}_{R^{-1},F}}
& & \Psi(GR) \Psi(R^{-1}F) \ar[rr]_(.55){\com^{\Psi}_{GR,R^{-1}F}} & & \Psi(GRR^{-1}F)   }$$

\newpage

\setlength{\textheight}{20.9cm}
\oddsidemargin -2cm
\begin{landscape}
{\tiny
$$\xymatrix{
\chi_{\mm_{3}} \Psi(G) \Psi(F)
\ar[dd]_{\uni^{\Psi\Gamma(R_{2})}} \ar[dr]_{\uni^{\Psi\Gamma(R_{1})}} \ar[rrrrr]^{\com^{\Psi}_{G,F}}  &
  &  &  &  &
\chi_{\mm_{3}} \Psi(GF)
\ar[d]^{\uni^{\Psi\Gamma(R_{1})}}
\\
 &
\chi_{\mm_{3}} \Psi(G) \Psi(F) \Psi\Gamma(R_{1}) \Psi(\Gamma(R_{1})^{\sharp})
\ar[d]_{\uni^{\Psi\Gamma(R_{2})}} \ar[rrrr]^{\com^{\Psi}_{G,F}} &
  &  &  &
\chi_{\mm_{3}} \Psi(GF) \Psi\Gamma(R_{1}) \Psi(\Gamma(R_{1})^{\sharp})
\ar[dl]^(.4){\uni^{\Gamma(R_{2})}}  \ar[d]^{\com^{\Psi}_{GF}\Gamma(R_{1})^{\sharp}}
\\
\chi_{\mm_{3}} \Psi(G) \Psi\Gamma(R_{2}) \Psi(\Gamma(R_{2})^{\sharp}) \Psi(F)
\ar[d]_{\com^{\Psi}_{G,\Gamma(R_{2})}} \ar[r]^(.43){\uni^{\Psi\Gamma(R_{1})}}  &
\chi_{\mm_{3}} \Psi(G) \Psi\Gamma(R_{2}) \Psi(\Gamma(R_{2})^{\sharp}) \Psi(F) \Psi\Gamma(R_{1}) \Psi(\Gamma(R_{1})^{\sharp})
\ar[d]_{\com^{\Psi}_{G,\Gamma(R_{2})}} \ar[r]^(.52){\com^{\Psi}_{} \ast \com^{\Psi}_{}}  &
\chi_{\mm_{3}} \Psi(G \Gamma(R_{2})) \Psi(\Gamma(R_{2})^{\sharp} F) \Psi\Gamma(R_{1}) \Psi(\Gamma(R_{1})^{\sharp})
\ar[d]_{\com^{\Psi}_{\Gamma(R_{2})^{\sharp} F,\Gamma(R_{1})}} \ar[rr]^{\com^{\Psi}_{G\Gamma(R_{2}),\Gamma(R_{2})^{\sharp}F}} &
  &
\chi_{\mm_{3}} \Psi(G \Gamma(R_{2})\Gamma(R_{2})^{\sharp} F) \Psi\Gamma(R_{1}) \Psi(\Gamma(R_{1})^{\sharp})
\ar[d]_{\com^{\Psi}_{G\Gamma(R_{2})\Gamma(R_{2})^{\sharp}F,\Gamma(R_{1})}} &
\chi_{\mm_{3}} \Psi(GF\Gamma(R_{1})) \Psi(\Gamma(R_{1})^{\sharp})
\ar[dl]^(.4){\uni^{\Gamma(R_{2})}} \ar[d]^{\tpn^{\chi}_{GF\Gamma(R_{1})}}
\\
\chi_{\mm_{3}} \Psi(G \Gamma(R_{2})) \Psi(\Gamma(R_{2})^{\sharp}) \Psi(F)
\ar[d]_{\tpn^{\chi}_{G \Gamma(R_{2})}} \ar[r]^(.43){\uni^{\Psi\Gamma(R_{1})}}  &
\chi_{\mm_{3}} \Psi(G \Gamma(R_{2})) \Psi(\Gamma(R_{2})^{\sharp}) \Psi(F) \Psi\Gamma(R_{1}) \Psi(\Gamma(R_{1})^{\sharp})
\ar[d]_{\tpn^{\chi}_{G \Gamma(R_{2})}} \ar[r]^{\com^{\Psi}_{}}  &
\chi_{\mm_{3}} \Psi(G \Gamma(R_{2})) \Psi(\Gamma(R_{2})^{\sharp} F \Gamma(R_{1})) \Psi(\Gamma(R_{1})^{\sharp})
\ar[d]_{\tpn^{\chi}_{G \Gamma(R_{2})}} \ar[rr]^{\com^{\Psi}_{G\Gamma(R_{2}),\Gamma(R_{2})^{\sharp}F\Gamma(R_{1})}}  &
  &
\chi_{\mm_{3}} \Psi(G \Gamma(R_{2}) \Gamma(R_{2})^{\sharp} F \Gamma(R_{1})) \Psi(\Gamma(R_{1})^{\sharp})
\ar[d]_{\tpn^{\chi}_{G\Gamma(R_{2})\Gamma(R_{2})^{\sharp} F \Gamma(R_{1})}}  &
\Psi'(GF \Gamma(R_{1})) \chi_{\mm_{1}^{p}} \Psi(\Gamma(R_{1})^{\sharp}) \ar@{=}[ddddddd] \ar[dl]^(.4){\uni^{\Gamma(R_{2})}}
\\
\Psi'(G \Gamma(R_{2})) \chi_{\mm_{2}^{p}} \Psi(\Gamma(R_{2})^{\sharp}) \Psi(F)
\ar[d]_{\left(\com^{\Psi'}_{G,\Gamma(R_{2})}\right)^{-1}} \ar[r]^(.43){\uni^{\Psi\Gamma(R_{1})}}  &
\Psi'(G \Gamma(R_{2})) \chi_{\mm_{2}^{p}} \Psi(\Gamma(R_{2})^{\sharp}) \Psi(F) \Psi\Gamma(R_{1}) \Psi(\Gamma(R_{1})^{\sharp})
\ar[d]_{\left(\com^{\Psi'}_{G,\Gamma(R_{2})}\right)^{-1}} \ar[r]^{\com^{\Psi}_{}}  &
\Psi'(G \Gamma(R_{2})) \chi_{\mm_{2}^{p}} \Psi(\Gamma(R_{2})^{\sharp} F \Gamma(R_{1})) \Psi(\Gamma(R_{1})^{\sharp})
\ar[d]_{\left(\com^{\Psi'}_{G,\Gamma(R_{2})}\right)^{-1}} \ar[r]^{\tpn^{\chi}_{}}  &
\Psi'(G \Gamma(R_{2})) \Psi'(\Gamma(R_{2})^{\sharp} F \Gamma(R_{1})) \chi_{\mm_{1}^{p}} \Psi(\Gamma(R_{1})^{\sharp})
\ar[d]_{\left(\com^{\Psi'}_{G,\Gamma(R_{2})}\right)^{-1}} \ar[r]^{\com^{\Psi'}_{}}  &
\Psi'(G \Gamma(R_{2}) \Gamma(R_{2})^{\sharp} F \Gamma(R_{1})) \chi_{\mm_{1}^{p}} \Psi(\Gamma(R_{1})^{\sharp})
\ar[d]_{\left(\com^{\Psi'}_{}\right)^{-1}}  &
\\
\Psi'(G) \Psi'\Gamma(R_{2}) \chi_{\mm_{2}^{p}} \Psi(\Gamma(R_{2})^{\sharp}) \Psi(F)
\ar[d]_{\left(\tpn^{\chi}_{\Gamma(R_{2})}\right)^{-1}} \ar[r]^(.43){\uni^{\Psi\Gamma(R_{1})}}  &
\Psi'(G) \Psi'\Gamma(R_{2}) \chi_{\mm_{2}^{p}} \Psi(\Gamma(R_{2})^{\sharp}) \Psi(F) \Psi\Gamma(R_{1}) \Psi(\Gamma(R_{1})^{\sharp})
\ar[d]_{\left(\tpn^{\chi}_{\Gamma(R_{2})}\right)^{-1}} \ar[r]^{\com^{\Psi}_{}}  &
\Psi'(G) \Psi'\Gamma(R_{2}) \chi_{\mm_{2}^{p}} \Psi(\Gamma(R_{2})^{\sharp} F \Gamma(R_{1})) \Psi(\Gamma(R_{1})^{\sharp})
\ar[d]_{\left(\tpn^{\chi}_{\Gamma(R_{2})}\right)^{-1}} \ar[r]^{\tpn^{\chi}_{}}  &
\Psi'(G) \Psi'\Gamma(R_{2}) \Psi'(\Gamma(R_{2})^{\sharp} F \Gamma(R_{1})) \chi_{\mm_{1}^{p}} \Psi(\Gamma(R_{1})^{\sharp})
\ar[r]^{\com^{\Psi'}_{}}  &
\Psi'(G) \Psi'(\Gamma(R_{2}) \Gamma(R_{2})^{\sharp} F \Gamma(R_{1})) \chi_{\mm_{1}^{p}} \Psi(\Gamma(R_{1})^{\sharp}) &
\\
\Psi'(G) \chi_{\mm_{2}} \Psi\Gamma(R_{2}) \Psi(\Gamma(R_{2})^{\sharp}) \Psi(F)
\ar[d]_{\left(\uni^{\Psi\Gamma(R_{2})}\right)^{-1}} \ar[r]^(.43){\uni^{\Psi\Gamma(R_{1})}}  &
\Psi'(G) \chi_{\mm_{2}} \Psi\Gamma(R_{2}) \Psi(\Gamma(R_{2})^{\sharp}) \Psi(F) \Psi\Gamma(R_{1}) \Psi(\Gamma(R_{1})^{\sharp})
\ar[d]_{\com^{\Psi}_{F,\Gamma(R_{1})}} \ar@/_1pc/[ddl]_(.55){\left(\uni^{\Psi\Gamma(R_{2})}\right)^{-1}}
\ar[r]^{\com^{\Psi}_{}}  &
\Psi'(G) \chi_{\mm_{2}} \Psi\Gamma(R_{2}) \Psi(\Gamma(R_{2})^{\sharp} F \Gamma(R_{1})) \Psi(\Gamma(R_{1})^{\sharp})
\ar[d]_{\com^{\Psi}_{}}  &
  &  &
\\
\Psi'(G) \chi_{\mm_{2}} \Psi(F)
\ar[d]_{\uni^{\Psi\Gamma(R_{1})}}  &
\Psi'(G) \chi_{\mm_{2}} \Psi\Gamma(R_{2}) \Psi(\Gamma(R_{2})^{\sharp}) \Psi(F \Gamma(R_{1})) \Psi(\Gamma(R_{1})^{\sharp})
\ar[ddl]_(.4){\left(\uni^{\Psi\Gamma(R_{2})}\right)^{-1}}  \ar[r]^{\com^{\Psi}_{}}  &
\Psi'(G) \chi_{\mm_{2}} \Psi(\Gamma(R_{2}) \Gamma(R_{2})^{\sharp} F \Gamma(R_{1})) \Psi(\Gamma(R_{1})^{\sharp})
\ar[uurr]^{\tpn^{\chi}_{\Gamma(R_{2}) \Gamma(R_{2})^{\sharp} F \Gamma(R_{1})}}  &
  &  &
\\
\Psi'(G) \chi_{\mm_{2}} \Psi(F) \Psi\Gamma(R_{1}) \Psi(\Gamma(R_{1})^{\sharp})
\ar[d]_{\com^{\Psi}_{F,\Gamma(R_{1})}}  &  &
\Psi'(G) \Psi'(F \Gamma(R_{1})) \chi_{\mm_{1}^{p}} \Psi(\Gamma(R_{1})^{\sharp})
\ar[rrd]_(.3){\Psi'\left(G\uni^{\Gamma(R_{2})}\right)} \ar[rr]^{\Psi'\left(\uni^{\Gamma(R_{2})}\right)\un^{\Psi'}}
\ar[uuurr]_{\Psi'\left(\uni^{\Gamma(R_{2})} F \Gamma(R_{1})\right)} &
  &
\Psi'(G) \Psi'(\Gamma(R_{2} \Gamma(R_{2})^{\sharp}) \Psi'(F \Gamma(R_{1})) \chi_{\mm_{1}^{p}} \Psi(\Gamma(R_{1})^{\sharp})
\ar[d]_{\com^{\Psi'}_{G, \Gamma(R_{2}) \Gamma(R_{2})^{\sharp}}}
\ar[uuu]^{\com^{\Psi'}_{\Gamma(R_{2}) \Gamma(R_{2})^{\sharp}, F \Gamma(R_{1})}}  &
\\
\Psi'(G) \chi_{\mm_{2}} \Psi(F \Gamma(R_{1})) \Psi(\Gamma(R_{1})^{\sharp})
\ar[d]_{\tpn^{\chi}_{F\Gamma(R_{1})}}  \ar[rr]_{\Psi'\left(G\uni^{\Gamma(R_{2})}\right)}
\ar[rru]_(.7){\tpn^{\chi}_{F\Gamma(R_{1})}}  \ar[uurr]^{\Psi\left(\uni^{\Gamma(R_{2})} F \Gamma(R_{1})\right)}  &
  &
\Psi'(G \Gamma(R_{2}) \Gamma(R_{2})^{\sharp}) \chi_{\mm_{2}} \Psi'(F \Gamma(R_{1})) \Psi(\Gamma(R_{1})^{\sharp})
\ar[rr]^(.43){\tpn^{\chi}_{F\Gamma(R_{1})}}  &
  &
\Psi'(G \Gamma(R_{2}) \Gamma(R_{2})^{\sharp}) \Psi'(F \Gamma(R_{1})) \chi_{\mm_{1}^{p}} \Psi(\Gamma(R_{1})^{\sharp})
\ar@/_4pc/@<-33ex>[uuuuu]^{\com^{\Psi'}_{G \Gamma(R_{2}) \Gamma(R_{2})^{\sharp}, F \Gamma(R_{1})}}
  &
\\
\Psi'(G) \Psi'(F \Gamma(R_{1})) \chi_{\mm_{1}^{p}} \Psi(\Gamma(R_{1})^{\sharp})
\ar[d]_{\left(\com^{\Psi}_{F,\Gamma(R_{1})}\right)^{-1}} \ar[rrrrr]_{\com^{\Psi'}_{G,F\Gamma(R_{1})}}
\ar[urrrr]_(.7){\Psi'\left(G\uni^{\Gamma(R_{2})}\right)}  &
  &  &  &  &
\Psi'(GF \Gamma(R_{1})) \chi_{\mm_{1}^{p}} \Psi(\Gamma(R_{1})^{\sharp})
\ar[d]^{\left(\com^{\Psi}_{GF,\Gamma(R_{1})}\right)^{-1}}
\\
\Psi'(G) \Psi'(F) \Psi'\Gamma(R_{1}) \chi_{\mm_{1}^{p}} \Psi(\Gamma(R_{1})^{\sharp})
\ar[d]_{\left(\tpn^{\chi}_{\Gamma(R_{1})}\right)^{-1}} \ar[rrrrr]_{\com^{\Psi'}_{G,F}}  &
  &  &  &  &
\Psi'(GF) \Psi'\Gamma(R_{1}) \chi_{\mm_{1}^{p}} \Psi(\Gamma(R_{1})^{\sharp})
\ar[d]^{\left(\tpn^{\chi}_{\Gamma(R_{1})}\right)^{-1}}
\\
\Psi'(G) \Psi'(F) \chi_{\mm_{1}} \Psi\Gamma(R_{1}) \Psi(\Gamma(R_{1})^{\sharp})
\ar[d]_{\left(\uni^{\Psi\Gamma(R_{1})}\right)^{-1}}  \ar[rrrrr]_{\com^{\Psi'}_{G,F}}
  &  &  &  &  &
\Psi'(GF) \chi_{\mm_{1}} \Psi\Gamma(R_{1}) \Psi(\Gamma(R_{1})^{\sharp})
\ar[d]^{\left(\uni^{\Psi\Gamma(R_{1})}\right)^{-1}}
\\
\Psi'(G) \Psi'(F) \chi_{\mm_{1}}
\ar[rrrrr]^(.55){\com^{\Psi'}_{G,F}}  &
  &  &  &  &
\Psi'(GF) \chi_{\mm_{1}}       }$$
}
\end{landscape}
\newpage
\setlength{\textheight}{65.4cm}
\oddsidemargin -1.7cm
\begin{landscape}
{\tiny
$$\xymatrix{
\chi_{\mm_{3}} \Psi(G) \Psi(F)
\ar[dd]_{\uni^{\Psi\Gamma(R_{2})}} \ar[dr]_{\uni^{\Psi\Gamma(R_{1})}} \ar[rrrrr]^{\com^{\Psi}_{G,F}}  &
  &  &  &  &
\chi_{\mm_{3}} \Psi(GF)
\ar[d]^{\uni^{\Psi\Gamma(R_{1})}}
\\
 &
\chi_{\mm_{3}} \Psi(G) \Psi(F) \Psi\Gamma(R_{1}) \Psi(\Gamma(R_{1})^{\sharp})
\ar[d]_{\uni^{\Psi\Gamma(R_{2})}} \ar[rrrr]^{\com^{\Psi}_{G,F}} &
  &  &  &
\chi_{\mm_{3}} \Psi(GF) \Psi\Gamma(R_{1}) \Psi(\Gamma(R_{1})^{\sharp})
\ar[dl]^(.4){\uni^{\Gamma(R_{2})}} \ar[d]^{\com^{\Psi}_{GF}\Gamma(R_{1})^{\sharp}}
\\
\chi_{\mm_{3}} \Psi(G) \Psi\Gamma(R_{2}) \Psi(\Gamma(R_{2})^{\sharp}) \Psi(F)
\ar[d]_{\com^{\Psi}_{G,\Gamma(R_{2})}} \ar[r]^(.43){\uni^{\Psi\Gamma(R_{1})}}  &
\chi_{\mm_{3}} \Psi(G) \Psi\Gamma(R_{2}) \Psi(\Gamma(R_{2})^{\sharp}) \Psi(F) \Psi\Gamma(R_{1}) \Psi(\Gamma(R_{1})^{\sharp})
\ar[d]_{\com^{\Psi}_{G,\Gamma(R_{2})}} \ar[r]^(.52){\com^{\Psi}_{} \ast \com^{\Psi}_{}}  &
\chi_{\mm_{3}} \Psi(G \Gamma(R_{2})) \Psi(\Gamma(R_{2})^{\sharp} F) \Psi\Gamma(R_{1}) \Psi(\Gamma(R_{1})^{\sharp})
\ar[d]_{\com^{\Psi}_{\Gamma(R_{2})^{\sharp} F,\Gamma(R_{1})}} \ar[rr]^{\com^{\Psi}_{G\Gamma(R_{2}),\Gamma(R_{2})^{\sharp}F}} &
  &
\chi_{\mm_{3}} \Psi(G \Gamma(R_{2})\Gamma(R_{2})^{\sharp} F) \Psi\Gamma(R_{1}) \Psi(\Gamma(R_{1})^{\sharp})
\ar[d]_{\com^{\Psi}_{G\Gamma(R_{2})\Gamma(R_{2})^{\sharp}F,\Gamma(R_{1})}} &
\chi_{\mm_{3}} \Psi(GF\Gamma(R_{1})) \Psi(\Gamma(R_{1})^{\sharp})
\ar[dl]^(.4){\uni^{\Gamma(R_{2})}} \ar[d]^{\tpn^{\chi}_{GF\Gamma(R_{1})}}
\\
\chi_{\mm_{3}} \Psi(G \Gamma(R_{2})) \Psi(\Gamma(R_{2})^{\sharp}) \Psi(F)
\ar[d]_{\tpn^{\chi}_{G \Gamma(R_{2})}} \ar[r]^(.43){\uni^{\Psi\Gamma(R_{1})}}  &
\chi_{\mm_{3}} \Psi(G \Gamma(R_{2})) \Psi(\Gamma(R_{2})^{\sharp}) \Psi(F) \Psi\Gamma(R_{1}) \Psi(\Gamma(R_{1})^{\sharp})
\ar[d]_{\tpn^{\chi}_{G \Gamma(R_{2})}} \ar[r]^{\com^{\Psi}_{}}  &
\chi_{\mm_{3}} \Psi(G \Gamma(R_{2})) \Psi(\Gamma(R_{2})^{\sharp} F \Gamma(R_{1})) \Psi(\Gamma(R_{1})^{\sharp})
\ar[d]_{\tpn^{\chi}_{G \Gamma(R_{2})}} \ar[rr]^{\com^{\Psi}_{G\Gamma(R_{2}),\Gamma(R_{2})^{\sharp}F\Gamma(R_{1})}}  &
  &
\chi_{\mm_{3}} \Psi(G \Gamma(R_{2}) \Gamma(R_{2})^{\sharp} F \Gamma(R_{1})) \Psi(\Gamma(R_{1})^{\sharp})
\ar[d]_{\tpn^{\chi}_{G\Gamma(R_{2})\Gamma(R_{2})^{\sharp} F \Gamma(R_{1})}}  &
\Psi'(GF \Gamma(R_{1})) \chi_{\mm_{1}^{p}} \Psi(\Gamma(R_{1})^{\sharp}) \ar@{=}[ddddddd] \ar[dl]^(.4){\uni^{\Gamma(R_{2})}}
\\
\Psi'(G \Gamma(R_{2})) \chi_{\mm_{2}^{p}} \Psi(\Gamma(R_{2})^{\sharp}) \Psi(F)
\ar[d]_{\left(\com^{\Psi'}_{G,\Gamma(R_{2})}\right)^{-1}} \ar[r]^(.43){\uni^{\Psi\Gamma(R_{1})}}  &
\Psi'(G \Gamma(R_{2})) \chi_{\mm_{2}^{p}} \Psi(\Gamma(R_{2})^{\sharp}) \Psi(F) \Psi\Gamma(R_{1}) \Psi(\Gamma(R_{1})^{\sharp})
\ar[d]_{\left(\com^{\Psi'}_{G,\Gamma(R_{2})}\right)^{-1}} \ar[r]^{\com^{\Psi}_{}}  &
\Psi'(G \Gamma(R_{2})) \chi_{\mm_{2}^{p}} \Psi(\Gamma(R_{2})^{\sharp} F \Gamma(R_{1})) \Psi(\Gamma(R_{1})^{\sharp})
\ar[d]_{\left(\com^{\Psi'}_{G,\Gamma(R_{2})}\right)^{-1}} \ar[r]^{\tpn^{\chi}_{}}  &
\Psi'(G \Gamma(R_{2})) \Psi'(\Gamma(R_{2})^{\sharp} F \Gamma(R_{1})) \chi_{\mm_{1}^{p}} \Psi(\Gamma(R_{1})^{\sharp})
\ar[d]_{\left(\com^{\Psi'}_{G,\Gamma(R_{2})}\right)^{-1}} \ar[r]^{\com^{\Psi'}_{}}  &
\Psi'(G \Gamma(R_{2}) \Gamma(R_{2})^{\sharp} F \Gamma(R_{1})) \chi_{\mm_{1}^{p}} \Psi(\Gamma(R_{1})^{\sharp})
\ar[d]_{\left(\com^{\Psi'}_{}\right)^{-1}}  &
\\
\Psi'(G) \Psi'\Gamma(R_{2}) \chi_{\mm_{2}^{p}} \Psi(\Gamma(R_{2})^{\sharp}) \Psi(F)
\ar[d]_{\left(\tpn^{\chi}_{\Gamma(R_{2})}\right)^{-1}} \ar[r]^(.43){\uni^{\Psi\Gamma(R_{1})}}  &
\Psi'(G) \Psi'\Gamma(R_{2}) \chi_{\mm_{2}^{p}} \Psi(\Gamma(R_{2})^{\sharp}) \Psi(F) \Psi\Gamma(R_{1}) \Psi(\Gamma(R_{1})^{\sharp})
\ar[d]_{\left(\tpn^{\chi}_{\Gamma(R_{2})}\right)^{-1}} \ar[r]^{\com^{\Psi}_{}}  &
\Psi'(G) \Psi'\Gamma(R_{2}) \chi_{\mm_{2}^{p}} \Psi(\Gamma(R_{2})^{\sharp} F \Gamma(R_{1})) \Psi(\Gamma(R_{1})^{\sharp})
\ar[d]_{\left(\tpn^{\chi}_{\Gamma(R_{2})}\right)^{-1}} \ar[r]^{\tpn^{\chi}_{}}  &
\Psi'(G) \Psi'\Gamma(R_{2}) \Psi'(\Gamma(R_{2})^{\sharp} F \Gamma(R_{1})) \chi_{\mm_{1}^{p}} \Psi(\Gamma(R_{1})^{\sharp})
\ar[r]^{\com^{\Psi'}_{}}  &
\Psi'(G) \Psi'(\Gamma(R_{2}) \Gamma(R_{2})^{\sharp} F \Gamma(R_{1})) \chi_{\mm_{1}^{p}} \Psi(\Gamma(R_{1})^{\sharp}) &
\\
\Psi'(G) \chi_{\mm_{2}} \Psi\Gamma(R_{2}) \Psi(\Gamma(R_{2})^{\sharp}) \Psi(F)
\ar[d]_{\left(\uni^{\Psi\Gamma(R_{2})}\right)^{-1}} \ar[r]^(.43){\uni^{\Psi\Gamma(R_{1})}}  &
\Psi'(G) \chi_{\mm_{2}} \Psi\Gamma(R_{2}) \Psi(\Gamma(R_{2})^{\sharp}) \Psi(F) \Psi\Gamma(R_{1}) \Psi(\Gamma(R_{1})^{\sharp})
\ar[d]_{\com^{\Psi}_{F,\Gamma(R_{1})}} \ar@/_1pc/[ddl]_(.55){\left(\uni^{\Psi\Gamma(R_{2})}\right)^{-1}}
\ar[r]^{\com^{\Psi}_{}}  &
\Psi'(G) \chi_{\mm_{2}} \Psi\Gamma(R_{2}) \Psi(\Gamma(R_{2})^{\sharp} F \Gamma(R_{1})) \Psi(\Gamma(R_{1})^{\sharp})
\ar[d]_{\com^{\Psi}_{}}  &
  &  &
\\
\Psi'(G) \chi_{\mm_{2}} \Psi(F)
\ar[d]_{\uni^{\Psi\Gamma(R_{1})}}  &
\Psi'(G) \chi_{\mm_{2}} \Psi\Gamma(R_{2}) \Psi(\Gamma(R_{2})^{\sharp}) \Psi(F \Gamma(R_{1})) \Psi(\Gamma(R_{1})^{\sharp})
\ar[ddl]_(.4){\left(\uni^{\Psi\Gamma(R_{2})}\right)^{-1}}  \ar[r]^{\com^{\Psi}_{}}  &
\Psi'(G) \chi_{\mm_{2}} \Psi(\Gamma(R_{2}) \Gamma(R_{2})^{\sharp} F \Gamma(R_{1})) \Psi(\Gamma(R_{1})^{\sharp})
\ar[uurr]^{\tpn^{\chi}_{\Gamma(R_{2}) \Gamma(R_{2})^{\sharp} F \Gamma(R_{1})}}  &
  &  &
\\
\Psi'(G) \chi_{\mm_{2}} \Psi(F) \Psi\Gamma(R_{1}) \Psi(\Gamma(R_{1})^{\sharp})
\ar[d]_{\com^{\Psi}_{F,\Gamma(R_{1})}}  &  &
\Psi'(G) \Psi'(F \Gamma(R_{1})) \chi_{\mm_{1}^{p}} \Psi(\Gamma(R_{1})^{\sharp})
\ar[rrd]_(.3){\Psi'\left(G\uni^{\Gamma(R_{2})}\right)} \ar[rr]^{\Psi'\left(\uni^{\Gamma(R_{2})}\right)\un^{\Psi'}}
\ar[uuurr]_{\Psi'\left(\uni^{\Gamma(R_{2})} F \Gamma(R_{1})\right)} &
  &
\Psi'(G) \Psi'(\Gamma(R_{2} \Gamma(R_{2})^{\sharp}) \Psi'(F \Gamma(R_{1})) \chi_{\mm_{1}^{p}} \Psi(\Gamma(R_{1})^{\sharp})
\ar[d]_{\com^{\Psi'}_{G, \Gamma(R_{2}) \Gamma(R_{2})^{\sharp}}}
\ar[uuu]^{\com^{\Psi'}_{\Gamma(R_{2}) \Gamma(R_{2})^{\sharp}, F \Gamma(R_{1})}}  &
\\
\Psi'(G) \chi_{\mm_{2}} \Psi(F \Gamma(R_{1})) \Psi(\Gamma(R_{1})^{\sharp})
\ar[d]_{\tpn^{\chi}_{F\Gamma(R_{1})}}  \ar[rr]_{\Psi'\left(G\uni^{\Gamma(R_{2})}\right)}
\ar[rru]_(.7){\tpn^{\chi}_{F\Gamma(R_{1})}}  \ar[uurr]^{\Psi\left(\uni^{\Gamma(R_{2})} F \Gamma(R_{1})\right)}  &
  &
\Psi'(G \Gamma(R_{2}) \Gamma(R_{2})^{\sharp}) \chi_{\mm_{2}} \Psi'(F \Gamma(R_{1})) \Psi(\Gamma(R_{1})^{\sharp})
\ar[rr]^(.43){\tpn^{\chi}_{F\Gamma(R_{1})}}  &
  &
\Psi'(G \Gamma(R_{2}) \Gamma(R_{2})^{\sharp}) \Psi'(F \Gamma(R_{1})) \chi_{\mm_{1}^{p}} \Psi(\Gamma(R_{1})^{\sharp})
\ar@/_4pc/@<-33ex>[uuuuu]^{\com^{\Psi'}_{G \Gamma(R_{2}) \Gamma(R_{2})^{\sharp}, F \Gamma(R_{1})}}
  &
\\
\Psi'(G) \Psi'(F \Gamma(R_{1})) \chi_{\mm_{1}^{p}} \Psi(\Gamma(R_{1})^{\sharp})
\ar[d]_{\left(\com^{\Psi}_{F,\Gamma(R_{1})}\right)^{-1}} \ar[rrrrr]_{\com^{\Psi'}_{G,F\Gamma(R_{1})}}
\ar[urrrr]_(.7){\Psi'\left(G\uni^{\Gamma(R_{2})}\right)}  &
  &  &  &  &
\Psi'(GF \Gamma(R_{1})) \chi_{\mm_{1}^{p}} \Psi(\Gamma(R_{1})^{\sharp})
\ar[d]^{\left(\com^{\Psi}_{GF,\Gamma(R_{1})}\right)^{-1}}
\\
\Psi'(G) \Psi'(F) \Psi'\Gamma(R_{1}) \chi_{\mm_{1}^{p}} \Psi(\Gamma(R_{1})^{\sharp})
\ar[d]_{\left(\tpn^{\chi}_{\Gamma(R_{1})}\right)^{-1}} \ar[rrrrr]_{\com^{\Psi'}_{G,F}}  &
  &  &  &  &
\Psi'(GF) \Psi'\Gamma(R_{1}) \chi_{\mm_{1}^{p}} \Psi(\Gamma(R_{1})^{\sharp})
\ar[d]^{\left(\tpn^{\chi}_{\Gamma(R_{1})}\right)^{-1}}
\\
\Psi'(G) \Psi'(F) \chi_{\mm_{1}} \Psi\Gamma(R_{1}) \Psi(\Gamma(R_{1})^{\sharp})
\ar[d]_{\left(\uni^{\Psi\Gamma(R_{1})}\right)^{-1}}  \ar[rrrrr]_{\com^{\Psi'}_{G,F}}
  &  &  &  &  &
\Psi'(GF) \chi_{\mm_{1}} \Psi\Gamma(R_{1}) \Psi(\Gamma(R_{1})^{\sharp})
\ar[d]^{\left(\uni^{\Psi\Gamma(R_{1})}\right)^{-1}}
\\
\Psi'(G) \Psi'(F) \chi_{\mm_{1}}
\ar[rrrrr]^(.55){\com^{\Psi'}_{G,F}}  &
  &  &  &  &
\Psi'(GF) \chi_{\mm_{1}}       }$$
}
\end{landscape}
\setlength{\textheight}{23cm}

\end{document}